\documentclass[envcountchap,chaprefs]{svmult}

\usepackage{graphicx}    
\usepackage{multicol}    

\begin{document}

\title{Random planar curves and Schramm-Loewner evolutions}
\titlerunning {Random planar curves}
\author{Wendelin Werner} 
\institute{Universit\'e Paris-Sud and IUF\\
Laboratoire de Math\'ematiques, Universit\'e Paris-Sud,\\
B\^at. 425, 91405 Orsay cedex, France\\
\textit{e-mail: wendelin.werner@math.u-psud.fr}
}
%
%

\null \vskip 2cm

\centerline {\Large Random planar curves and Schramm-Loewner evolutions}

\bigbreak

\centerline {Lecture Notes from the 2002 Saint-Flour summer school}
\centerline {(final version)}

\bigbreak

\centerline {\large Wendelin Werner}

\bigbreak

\centerline{Universit\'e Paris-Sud and IUF}

\bigbreak
\centerline {Laboratoire de Math\'ematiques, Universit\'e Paris-Sud,}
\centerline {
B\^at. 425, 91405 Orsay cedex, France}
\centerline {
\textit{e-mail: wendelin.werner@math.u-psud.fr}
}

\newtheorem{prop}[theorem]{Proposition}
\newcommand{\Prob} {{\bf P}}
\font \m=msbm10
\newcommand{\R}{{\hbox {\m R}}}
\newcommand{\C}{{\hbox {\m C}}}
\newcommand{\Z}{{\hbox {\m Z}}}
\newcommand{\N}{{\hbox {\m N}}}
\newcommand{\U}{{\hbox {\m U}}}
\def\H{{\hbox {\m H}}}
\def\P{{\bf P}}
\def \expect {{\bf E}}
\def \proof {{\medbreak
\noindent {\bf Proof. }}}
\def \eps {\varepsilon}
\def \gamm {{r}}
\def \bibc{{\section* {Bibliographical comments}}}
\def \ho {{\tau}}


\chapter* {Foreword and summary}

The goal of these lectures is to review some of the 
mathematical results that have been 
derived in the last years on conformal invariance,
scaling limits and properties of some two-dimensional random curves. 
The (distinguished)
audience of the Saint-Flour summer school consists mainly 
of probabilists and I therefore assume  
knowledge in stochastic calculus (It\^o's formula etc.), but no 
special background in basic complex analysis.
 
These lecture notes are neither a book nor a 
compilation of research papers. While preparing
 them, I realized
that it was hopeless to present all the recent results  
on this subject, or even to give the complete 
detailed proofs of a selected portion of them.
Maybe this 
will disappoint part of the audience but
the main goal of these lectures will be to try to transmit some
ideas and heuristics.
As a reader/part of an audience,
I often think that omitting details is dangerous, 
and that ideas are sometimes better understood when the 
complete proofs are given, but in the present case, partly because 
the technicalities often use complex analysis tools that 
the audience might not be so familiar with, partly also because of 
the limited number of lectures, 
I chose to 
focus on some selected results and on the main ideas of their 
proofs, sometimes omitting technical details, and giving 
references for those interested in full proofs or 
more results.
In the final chapter, I will briefly review what I omitted  
in these lectures, as well as work in progress or open
questions.

Of course, I would like to thank my coauthors 
Greg Lawler and Oded Schramm without which I would
not have been lecturing on this subject in Saint-Flour. 
Collaborating with them during these last years was 
a great pleasure and privilege.
Also, I would like to stress the fact that (almost)
none of the pictures in these notes are mine. Many thanks
to their authors Vincent Beffara, Tom Kennedy and
Oded Schramm. 
 I also take this opportunity to 
thank  Stas Smirnov,
 Rick Kenyon, as well as all my Orsay colleagues and
students who have directly or indirectly contributed to 
these lecture notes through
their work, comments and  
discussions. 
 
Finally, I owe many thanks to all participants of the 
summer school, as well as 
to all colleagues who have sent me their comments and remarks on the first 
draft of these notes that was distributed during the summer school and posted on the web
at that time. 

It has been a pleasure and a very rewarding experience to lecture in 
the studious, relaxed and enjoyable atmosphere of the 2002 St-Flour school.
I express my gratitude to all who have contributed to it,
my co-lecturers 
Jim Pitman and Boris Tsirelson, the Maison des Planchettes' staff, and last 
but not least,  Jean 
Picard, whose outstanding organization has been both efficient and discreet.

%

\bigbreak

Here is  a short description of these notes:
In the first introductory chapter, I will briefly describe two discrete 
models
(loop-erased random walks and critical percolation interfaces) 
that have now been proved to converge in their scaling
limit to SLE (Oded Schramm used these letters
 as shorthand for ``stochastic
Loewner Evolution'', but I will stick to 
Schramm-Loewner Evolution). Using these models,
I will try to show why it is natural to define this 
one-parameter family of  
random continuously growing 
processes based on Loewner's equation, and to 
introduce the difference between their chordal and 
radial versions.
 
The second chapter is a review of the necessary background 
on deterministic aspects of Loewner's equation in the 
upper half-plane, which is then used in Chapter 3 to 
define chordal SLE. Some 
first properties of this process are studied.  
In particular, some hitting probabilities are computed.

The fourth chapter is devoted to some special properties
of SLE that hold for some special values of the 
parameter $\kappa$: The locality property for 
 $SLE_6$, and the restriction property for $SLE_{8/3}$.
These are not surprising if one thinks of these processes as 
the respective scaling limits of critical percolation interfaces and 
self-avoiding walks, but somewhat 
surprising if one starts from the definition of SLE itself.
These properties are then used in Chapter 5, to make the link between 
the geometry of $SLE_{8/3}$, that of the outer boundary of a
planar Brownian motion and that of the outer 
boundary of $SLE_{6}$.

In Chapter 6, we define radial SLE which are processes 
defined in a similar way as chordal SLE except that they
are growing towards an interior point of the domain 
and not to a boundary point. We show in that chapter that 
radial and chordal SLE are very closely related, especially in 
the case $\kappa =6$. 

In Chapter 7, we show how to compute critical exponents
associated to SLE that describe the asymptotic decay
of certain probabilities (non-disconnection, non-intersection).
Using the relation between radial $SLE_6$, chordal $SLE_6$ and 
planar Brownian motion, we then use these computations in
Chapter 8 to 
determine the values of the critical exponents that describe the
decay of disconnection or non-intersection
probabilities for planar Brownian motions, which 
is one of the main goals of these lectures.
As already mentioned, it will not be possible to 
describe all proofs in detail, but I hope that all the 
main ideas and steps (that are spread over the first seven
chapters of these notes) are explained in sufficient
detail so that the reader can get an overview of the
proof. For simplicity, I will mainly focus on 
derivation of the disconnection exponent i.e. the
proof of the fact that the probability that a complex 
Brownian curve $Z[0,t]$ started from $Z_0=1$ disconnects
the origin from infinity decays like $t^{-1/8}$ when $t \to 
\infty$.
 
In Chapters 9 and 10, another important aspect
of SLE is discussed: The proofs that some curves arising 
in discrete models from statistical physics
converge to SLE in their scaling limit. The case of loop-erased
random walks and uniform spanning trees is treated in Chapter 9.
Chapter 10
is devoted to critical site percolation on the triangular
lattice, including a brief discussion 
of Stas Smirnov's proof of conformal invariance
and of its consequences.

A concluding chapter contains
a list of other results, work in progress 
and open questions.

\setcounter {chapter}0
\chapter {Introduction}

\section {General motivation}

One of the main aims of both statistical physics and probability 
theory is to study 
macroscopic systems consisting of many (i.e. in the limit when this number 
grows
to infinity)
small microscopic random inputs.
One may classify the results into two categories:
 In the limit, the behaviour
of the macroscopic system becomes deterministic 
(these are ``law of large number'' type 
of results, and large deviations can 
to some extent been used in this framework), or
random. The archetype for continuous random 
objects that appear as  scaling limit of finite systems
 is Brownian motion.
Note that it is the scaling limit of
 a large class of simple random walks, so 
that one might argue that Brownian motion is more
 universal than the 
discrete model (simple random walk) because there is no
 need to specify a 
lattice or a jump-distribution:
 it only captures the phenomenological properties of
the walks (mean zero, stationary increments etc.).

In two dimensions, Brownian motion has an important property which was
first observed by Paul L\'evy (\cite {Lev}, see e.g. \cite {LG,RY}
for ``modern'' proofs based on It\^o's formula) and  that can be 
heuristically related to the fact that it is the 
scaling  limit of simple random walks
on different lattices (which implies for 
instance invariance under rotations and under scaling):
 It is invariant under conformal transformations. 
Here is one way to state this property: Take a simply connected
open planar domain $D$ that contains the origin and is not equal to $\C$.
Consider planar Brownian motion $(B_t, t \in [0,\tau])$ started from $B_0=0$
up to its exit time $\tau=\tau_D$ of the domain $D$.
Suppose that $\Phi$ is a conformal map (that is, a one-to-one 
smooth map that preserves angles) from $D$ onto some 
other domain $D'$ with $\Phi (0) = 0$. Then,
there exists a (random) time change $A : [0, \sigma ] \to [0, \tau]$ so that 
$ ( \Phi (B_{A(s)}), s \in [0, \sigma ] )$
is planar Brownian motion started from $0$ and killed at
 its first exit time $\sigma$ of $D'$.
In other words, if
we forget about time-parametrization, the law of $\Phi (B)$ is 
again a Brownian motion.
As we shall see in these 
lectures, conformal invariance will turn out to be instrumental in 
the understanding of curves arising in more complicated setups.

\begin{figure}
\centerline{\includegraphics*[height=3in]{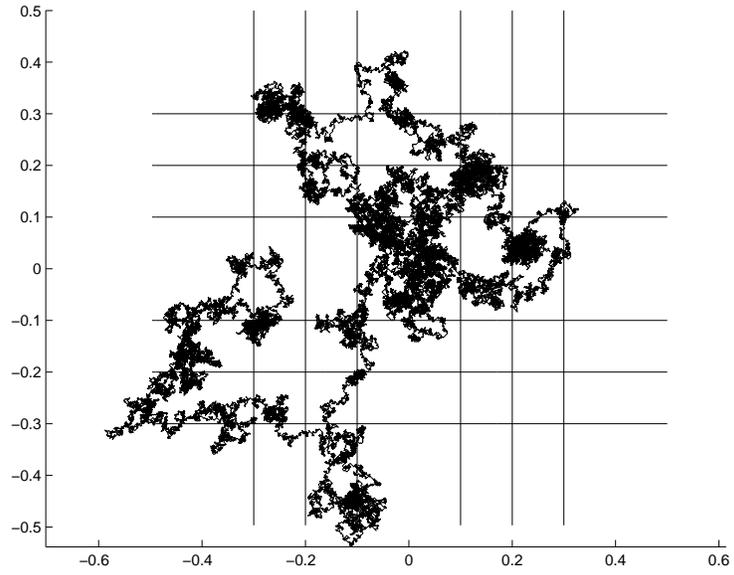}}
\caption {\label {f.rw1}Sample of a long simple random walk.}
\end{figure}

\begin{figure}
\centerline{\includegraphics*[height=3in]{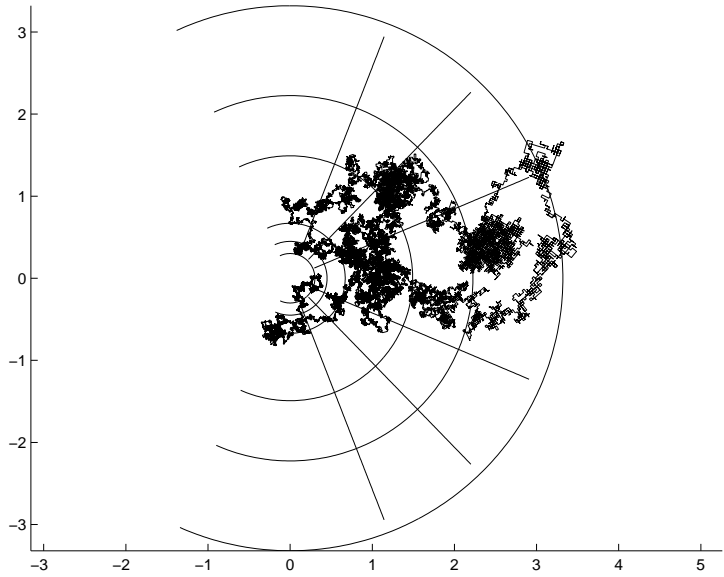}}
\caption{\label{f.rw2}The image of the previous sample under
an exponential mapping.}
\end{figure}

Actually, there exist only few 
known examples of probabilistic 
continuous models that 
are not directly related to Brownian motion.
For instance, under mild regularity conditions, continuous 
finite-dimensional Markov processes are solutions of 
stochastic differential equations
 and therefore constructed using Brownian motions.
If one looks for other types of continuous processes, one has therefore  
to give up the Markov property or the finite-dimensionality.
In many complex systems that we see around us
 and for which probability theory seems
a priori a well-suited tool (the shape of clouds, say), it is not possible
to explain the phenomena via Brownian motions, and there is still a long 
way to go for probabilists to understand their macroscopic behaviour.
  
In the present lectures, we shall focus on random planar curves. 
In two dimensions, (random) curves appear naturally as boundaries of 
domains, interfaces between two phases, level lines of random surfaces etc.
In all these cases, at least on microscopic level, the definition of the 
curve (say, as an interface) implies that it is a self-avoiding 
curve (or a simple closed loop). On the macroscopic scale, the 
continuous curves that we will be considering may have double-points
(in the scaling limit, simple curves may converge to curves with multiple
points), 
but self-crossings are forbidden.  
Of course, if $(\gamma_t, t \in [0, T])$ is such a random curve, we see 
that in general, this condition implies a strong correlation 
between $\gamma [0,t]$ and $\gamma [t, T]$, so that the Markov property is 
lost (if we look at these curves as living in the 
two-dimensional space). As we shall see, there is a
way to recover a 
Markov property for the random curves, using 
a coding of the curve in an infinite-dimensional  space of conformal 
maps.

\section {Loop-erased random walks}

In order to guide the intuition about the family of random curves that we will
be considering, it is helpful to have some discrete models in mind, for which 
one expects or can prove that they converge to this continuous object.
We therefore start these lectures with the description of one measure on 
discrete random curves that turns out to converge in the scaling limit.
This is actually the model that Oded Schramm considered when he invented these
random curves that he called $SLE$
 (for Stochastic Loewner Evolution, but we will 
replace this by Schramm-Loewner Evolution in these lectures).

For any ${\bf x} = (x_0, \ldots, x_m)$, we define the 
loop-erasure $L({\bf x})$ of ${\bf x}$ inductively as follows:
$L_0 = x_0$, and for all $j \ge 0$, we define
inductively  $n_j = \max \{ n \le m  \ : \ x_n = L_j \}$ and 
$$
L_{j+1} = X_{1 + n_j}  
$$
until $j = \sigma$ where $L_{\sigma} := x_{m}$.
In other words, we have
 erased the loops of ${\bf x}$ in chronological order. The number of 
steps $\sigma$ of $L$ is not fixed.

Suppose that $(X_n, n \ge 0)$
 is a recurrent Markov chain on a discrete state-space
${\cal S}$ started from $X_0 = x$.
 Suppose that $A \subset {\cal S}$ is non-empty,
and let $\tau_A$ denote the hitting time of $A$ by $X$.
Let $p(x,y)$ denote the transition probabilities  for the Markov chain $X$.
We define the loop-erasure $L=L(X[0, \tau_A])=L^A$ of $X$ up to its hitting 
time of $A$.
We
call  $\sigma$ the number of steps of $L^A$.
For $y \in A$ such that with positive probability $L^A (\sigma)
= X (\tau_A ) = y$, we call ${\cal L} (x,y; A)$ the law of 
$L^A$ conditioned on the event $\{ L^A (\sigma ) = y \}$.
In other words, it is the law of the loop-erasure of the
Markov chain $X$ conditioned to hit $A$ at $y$.

\begin {lemma}[Markovian property of LERW]
\label {l.lerw}
Consider $y_0, \ldots, y_j \in {\cal S}$ so that with positive
probability for ${\cal L}(x, y_0; A)$,
$$\{ L_\sigma = y_0, L_{\sigma -1} = y_1, \ldots, 
L_{\sigma -j} = y_j\}.$$
The conditional law of
$L[0,\sigma -j]$ given this event is 
 ${\cal L} ( x, y_j; A \cup \{ y_1, \ldots, y_j \} )$.
\end {lemma}

\proof
For each $A$ and $x \in {\cal A}$, we denote by 
$G(x, A)$ the expected number of visits by the Markov chain $X$
 before $\tau_A$ if 
$X_0  = x$.
Then, it is a simple exercise to check that for all $n \ge 1$, ${\bf w}
= (w_0, \ldots, w_n)$ with  $w_0=x$, 
$w_n \in A$ and $w_1, \ldots, w_{n-1} \in {\cal S} \setminus A$,
\begin {eqnarray*}
{\P
[ L^A = {\bf w} ]
} 
&=&
\sum_{{\bf x}  \ : \ L({\bf x}) = {\bf w} } \P [ X[0, \tau_A] = {\bf x} ] \\
&=&
 G( w_0, A) p(w_0, w_1) G(w_1, A \cup \{ w_0 \}) p(w_1, w_2)
 \cdots 
\\
&& \quad 
\times G( w_{n-1}, A \cup \{ w_0, w_1, \ldots, w_{n-2} \}) p( w_{n-1}, w_n).
\end {eqnarray*}
It is therefore natural to define the function 
$$
F( w_0, \ldots , w_{n-1} ; A) = 
\prod_{j=0}^{n-1} G ( w_j , A \cup \{ w_0, \ldots, w_{j-1} \} ).
$$
Again, it is a simple exercise on Markov chains to check that  
for all $A'$, $y$ and $y'$, 
$$
G ( y, A') G( y',A' \cup \{ y \}) = G(y',A') G(y, A' \cup \{ y' \}).
$$ 
It follows immediately that $F$ is in fact a 
symmetric function of its arguments. 
Hence,
\begin {eqnarray*}
\lefteqn {\P
[ L_0^A =w_0, \ldots , L_{\sigma}^A = w_n 
 | L_{\sigma} = w_n, L_{\sigma-1} = w_{n-1} ]
} \\
&=& 
\frac { 
p(w_{n-1}, w_n) 
G( w_{n-1}, A ) }{ \P [ L_\sigma = w_n , L_{\sigma-1} = w_{n-1} ] }
\\&&
\quad \times \prod_{j=0}^{n-2} 
p( w_j, w_{j+1}) 
G( w_{j}, (A \cup \{w_{n-1}\}) \cup \{ w_0, \ldots, w_{j-1} \} ) 
.\end {eqnarray*}
This readily implies the Lemma when $j=1$.
Iterating this $j$ times shows the Lemma.
\qed
\medbreak

This Lemma shows that it is in fact fairly natural to index the loop-erased 
path backwards
 (define $\gamma_j = L^A_{\sigma - j}$, so that $\gamma$ starts on $A$ and
goes back to $\gamma_\sigma= x$).
Then, the time-reversal of loop-erased (conditioned and stopped) Markov chains 
have themselves a Markovian-type property. 

Let us now
 come back to our two-dimensional setting: Suppose that $X$ is a simple
random walk 
on the grid $\delta \Z^2$ (we will then let the mesh $\delta$ of the 
lattice go to $0$) that is started from $0$. Let $D$ denote some 
simply connected domain $D$ with $0 \in D$ and $D \not= \C$, and let
$D_\delta = \delta \Z^2 \cap D$, 
$A=A_\delta = \delta Z^2 \setminus D$.
We are interested in the behaviour when $\delta \to 0$ of the law
of $\gamma^\delta$ which is defined as before as the time-reversed 
loop-erasure of $X[0, \tau_A]$. We now think on a heuristic 
level:
First, note that the law of $X_{\tau_A}$ converges to the harmonic measure 
on $\partial D$ from $0$, so that it is possible to study the behaviour of 
$\gamma^\delta$ conditional on the value of $\{ \gamma^\delta = y_0^\delta \}$ 
where $y_0^\delta \to y \in \partial D$ as $\delta \to 0$.
Second, one might argue that on the one hand, simple random walk converges 
to planar Brownian motion which is conformally invariant, and that on the 
other hand the chronological loop-erasing procedure is purely 
geometrical to conclude that when $
\delta  \to 0$, the law of $\gamma^\delta$ should converge to a conformal 
invariant curve that should be the loop-erasure of planar Brownian motion.

\begin{figure}
\centerline{\includegraphics*[height=2.7in]{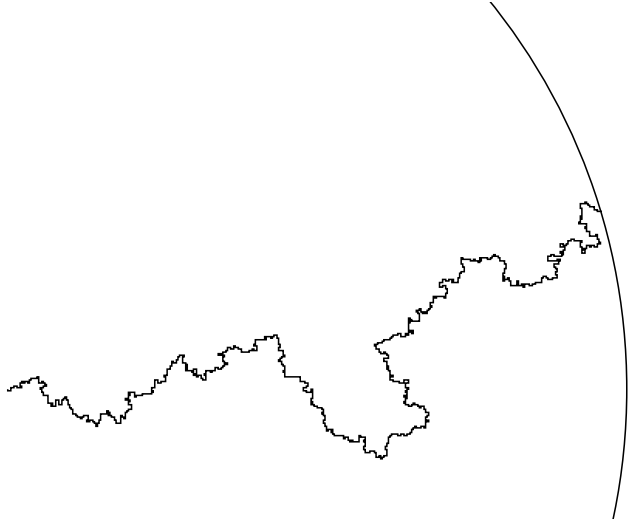}}
\caption{\label{f.lerw}A sample of the loop-erased
random walk.}
\end{figure}

Unfortunately (or fortunately!), the geometry of  planar Brownian curves is 
very complicated: It has points of any (even infinite) multiplicity
(see e.g. \cite {LG}), 
loops at any scale, so that there is no ``first'' loop to erase, and 
decisions about what small microscopic loops to erase first may propagate 
to the decisions about what macroscopic loops one should erase. In other words, there is 
no simple 
(even random) algorithm to loop-erase a Brownian path in chronological 
order.
Yet, the previous heuristic strongly suggests the law of $\gamma^\delta$ should
converge, 
and that the limiting law is invariant under conformal transformations:
The scaling limit of LERW in $D$ should be (modulo time-change) identical 
to the conformal image of the scaling limit of LERW in $D'$.
Furthermore, Lemma~\ref {l.lerw} should still be valid in the scaling limit.
We now show  that the combinations of these two 
properties 
in fact greatly reduce the family of possible scaling limits for LERW.

\section {Iterations of conformal maps and SLE}

We are therefore looking for the law of a random continuous curve $(\gamma_t, t \ge 0)$
with no self-crossings 
in the unit disc $\U$,  with $\gamma_0 = 1$, $\lim_{t \to \infty} \gamma_t = 0$ that could be 
the scaling limit of (time-reversed) loop-erased random walk on 
a grid approximation of $\U$ (conditioned to exit 
$\U$ near $1$).
Define for each $t \ge 0$, the conformal map $f_t$ from $\U \setminus \gamma [0,t]$
onto $\U$ which is normalized by $f_t (0)= 0$ and $f_t(\gamma_t)  = 1$
(actually, if $\gamma$ would have double-points, the domain of definition would 
be the connected component of $\U \setminus \gamma [0,t]$ that contains the 
origin, but let us a priori assume for convenience that $\gamma$ is a simple curve).

It is easy to check that $t \mapsto |f_t'(0)|$ is an increasing continuous 
function that goes to $\infty$ as $t \to \infty$ (see for instance 
\cite {A2}). Hence, it is possible to reparametrize $\gamma$ in such 
a way that 
\begin {equation}
\label {e.etime}
|f_t' (0)| = e^{t}.
\end {equation}
This is the natural parametrization in our context.
Indeed, let us now study the conditional law of $\gamma [t, \infty)$ given 
$\gamma[0,t]$. Lemma \ref {l.lerw} suggests that this law is the scaling limit of 
(time-reversed) LERW in the slit domain 
$\U
\setminus \gamma [0,t]$ conditioned to exit at $\gamma_t$, and
conformal invariance then says that this is the same (modulo time-reparametrization)
as the image under $z \mapsto f_t^{-1} ( z )$ of an independent 
copy  $\tilde \gamma$ of $\gamma$.
Note that if one composes conformal maps that preserve the origin, then the 
derivative at the origin multiply: This shows that in fact, no
time-change is necessary if we parametrize $\gamma$ (and $\tilde \gamma$) by
(\ref {e.etime}), in order for  the conditional law 
of $(\gamma_{t+s} , s \ge 0)$ given $\gamma [0,t]$ to be 
identical  
to that of 
$ (f_{t}^{-1} (\tilde \gamma_s ), s \ge 0)$.
In other words,
for all fixed $t \ge 0$,  
$$
(f_{t+s}, s \ge 0)  = (\tilde f_s \circ f_t , \ s \ge 0) 
\quad \hbox {in law}
$$
where $(\tilde f_s, s \ge 0)$ is an independent copy of $(f_s, s \ge 0)$.
In particular, $f_{2t}= \tilde f_t \circ f_t$
in law. Repeating this procedure, we see 
that for all $t \ge 0$ and all integer $n \ge 1$,
 $f_{nt}$ is the iteration of $n$ independent copies of $f_t$,
and that $f_t$ itself can be viewed as the iteration of $n$ independent 
copies of $f_{t/n}$. 
In other words, $(f_t, t \ge 0)$ is an ``infinitely divisible'' 
process of conformal maps, and $f_t$ is obtained by iterating 
infinitely many independent conformal  
maps that are infinitesimally close to the identity.

Back in the 1920's, Loewner observed that if $\gamma [0, \infty)$ 
is a simple continuous curve starting from $1$ in the unit disc, then it 
is naturally encoded via a continuous function $\zeta_t$ taking its values 
on the unit circle. Let us now describe briefly how it goes. 
Suppose, as in the previous section, that $\gamma (0) = 1$, 
$\lim_{t \to \infty} \gamma_t = 0$ and that $\gamma$ is parametrized in 
such a way that the modulus of the derivative at $0$ of the conformal map
$f_t$ from $U_t := \U \setminus \gamma [0,t]$
 into $\U$ that preserves the origin is 
$e^t$.
Define $\zeta_t = (f_t'(0)/|f_t'(0)|)^{-1}$. In other words, if $g_t$ denotes 
the conformal map from $U_t$ onto $\U$ such that $g_t (0) = 0$
and $g_t'(0) = e^t \in (0, \infty)$, then 
$$
\zeta_t = g_t ( \gamma_t)$$
and $g_t(z) = \zeta_t f_t(z)$.
One can note (see e.g. \cite {A2,Dur}) that for all $z \notin \gamma [0,t]$,
\begin {equation}
\label {e.radial}
\partial_t g_t (z) = - g_t (z) \frac { g_t(z) + \zeta_t }{g_t (z) - \zeta_t}
.\end {equation}
Hence, it is possible to recover $\gamma$ from $\zeta$ as follows:
For all $z \in \U$, define $g_t(z)$ as the unique solution to (\ref {e.radial})
starting from $z$. In case $g_t(z) = \zeta_t$ for some time $t$, then
define $\gamma_t = g_t^{-1} (\zeta_t)$
(we know already a priori that since $\gamma$ is a simple curve, the map $g_t^{-1}$ extends
continuously to the boundary). Note that if $g_t (z) = \zeta_t$, then $g_s (z)$ is not 
well-defined for $s \ge t$.

Hence, in order to define the random curve $\gamma$ that should be the scaling limit of 
loop-erased random walks, it suffices to define the random function $\zeta_t =\exp (iW_t)$,
where $(W_t, t \ge 0)$ is real-valued.
Our previous considerations suggest that the following conditions 
should be satisfied:
\begin {itemize}
\item
The process $W$ is almost surely continuous.
\item
The process $W$ has stationary increments (this is because 
$g_t$ is obtained by iterations of identically distributed 
conformal maps)
\item 
The laws of the processes $W$ and $-W$ are identical (this is because 
the law of $L$ and the law of the complex
conjugate $\overline L$ are 
identical).
\end {itemize}
The theory of Markov processes tells us that the 
only possible choices are:
$W_t = \beta_{\kappa t}$ where $\beta$ is standard
Brownian motion and $\kappa \ge 0$ a fixed constant.
In order to simplify some future notations, we will 
usually write 
$$
W_t = \sqrt {\kappa} B_t, t \ge 0
$$
where $(B_t, t \ge 0)$ is standard (one-dimensional) Brownian motion.

In summary, we have just seen that on a heuristic
level, if the scaling limit of loop-erased random walk
exists and is conformally invariant, then the 
scaling limit in the unit disk should be described as follows:
For some fixed constant $\kappa=\kappa_{LERW}$, define
$\zeta_t = \exp ( i \sqrt {\kappa} B_t), t \ge 0$, solve for each $z \in 
\U$, the equation (\ref {e.radial}) with $g_0 (z) = z$. This defines a 
conformal map $g_t$ from the subset $U_t$ of the unit disk 
onto $\U$. Then, one can construct $\gamma$ because
$$
U_t = \U \setminus \gamma [0,t] 
$$
and 
$$
\gamma_t = g_t^{-1} (\zeta_t).
$$
As we shall see later on in the lectures, this 
heuristic arguments can be 
made rigorous, and it will turn out that $\kappa_{LERW}= 2$.

\section {The critical percolation exploration process}

In the context
 of LERW, the  random curve  
joins a point in the inside of the domain to a point on the boundary of the 
domain. In statistical physics
models, one is often interested in ``interfaces''. Some of these interfaces 
appear to be random curves from one point on the boundary to another
point on the boundary. A natural setup
is to study curves from $0$ to infinity in 
the upper half-plane $\H := \{ x + iy \ : \ y >0 \}$.
Then, we look for random non-self-crossing curves $\gamma$ such that the 
law of $\gamma [t, \infty)$ given $\gamma [0,t]$ has the same law
than the conformal image of an independent copy $\tilde \gamma$
of $\gamma$ under a conformal map from $\H$ onto $\H \setminus \gamma [0,t]$
that maps $\infty$ onto itself and $0$ onto $\gamma_t$.

We now very briefly describe
 an important discrete model for which it has now 
also been
proved that it behaves in a conformally
invariant way in the scaling limit (more details on 
the model and its conformal covariance will
be given in Chapter 10): Critical site percolation on the
triangular lattice. Actually, it is more convenient to describe this 
in terms of cell-colouring  of the honeycombe lattice.
Suppose that a simply connected domain $D$ is fixed, as well as 
two distinct points $a$ and $b$ on $\partial D$. Let $D_\delta$ 
denote a suitably chosen
approximation of $D$ by a simply connected union of hexagonal
cells of size $\delta$. 
Let $a_\delta$ (resp. $b_\delta$) denote a vertex of the honeycombe 
lattice on $\partial D_\delta$ that is close to $a$ 
(resp. to $b$). Then, the cells on $\partial D_\delta$ 
can be divided into two ``arcs'' $B_\delta$ and $W_\delta$ in
such a way that $a_\delta$, $B_\delta$, $b_\delta$ and $W_\delta$
are oriented clockwise ``around'' $D_\delta$.
Decide that all hexagons in $B_\delta$ are colored in black and
that all hexagons in $W_\delta$ are colored in white.
On the other hand, all other cells in $D_\delta$ are chosen to 
be black or white with probability $1/2$ independently of
each other. Consider now the (random) path $\gamma_\delta$
from $a_\delta$ to $b_\delta$ that separates the cluster of 
black hexagons containing $B_\delta$ from the 
cluster of white hexagons containing $W_\delta$. 

\begin{figure}
\centerline{\includegraphics*[height=3in]{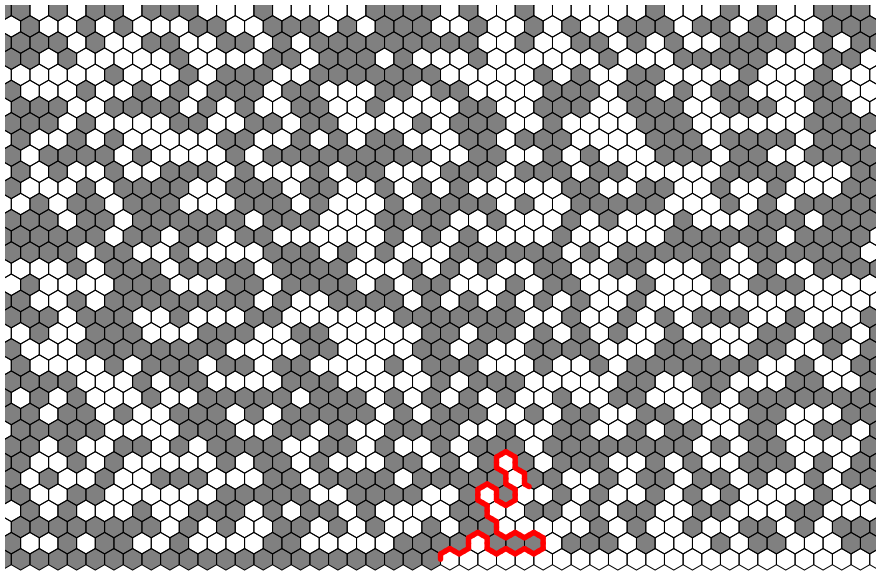}}
\caption{\label{f.exp1}
The beginning of the discrete exploration process.}
\end {figure}

For deep reasons that will be discussed later in these
lectures, it will turn out that when $\delta \to 0$, the law
of $\gamma_\delta$ converges towards that of a random curve $\gamma$
from $a$ to $b$ in $D$, and that the law
of that curve is conformally invariant: 
The law of $\Phi (\gamma)$ when $\Phi$ is a conformal map from 
$D$ onto $\Phi (D)$ is that of the corresponding path (i.e. of the 
scaling limit of percolation cluster interfaces) from 
$\Phi (a)$ to $\Phi (b)$ in $\Phi (D)$.
 
Again, on the discrete level, it is easy to see that $\gamma_\delta$
has the same type of Markovian property that LERW.  
More precisely, conditioning on the first steps of
$\gamma$
is equivalent to condition the percolation process to have black hexagons 
on the left-boundary of these steps and white hexagons on the right side.
Hence, the conditional law of the remaining steps is that of the 
percolation interface in the new domain obtained by slitting $D_\delta$
along the first steps of $\gamma_\delta$.
Figure \ref {f.exp1}
shows the beginning of the interface $\gamma_\delta$ in the 
case where $D$ is the upper half-plane.

Another equivalent way to define the interface $\gamma_\delta$ goes as follows:
It is a myopic self-avoiding walk. At each step $\gamma_\delta$ looks 
at its three neighbours (on the honeycomb lattice) and chooses at random 
one of the sites that it has not visited yet (there are one or two such sites 
since one site is anyway forbidden because it was the previous location of the 
walk).

This discrete walk in the upper half-plane is a very special discrete model  
that will turn out to converge to an SLE. The 
corresponding value of $\kappa$ is $6$.  Here, the starting point $a=0$ and
the end-point $b=\infty$ are both on the boundary of the domain, so that the 
previous definition of radial SLE is not well-suited anymore.

\begin{figure}
\centerline{\includegraphics*[height=3in]{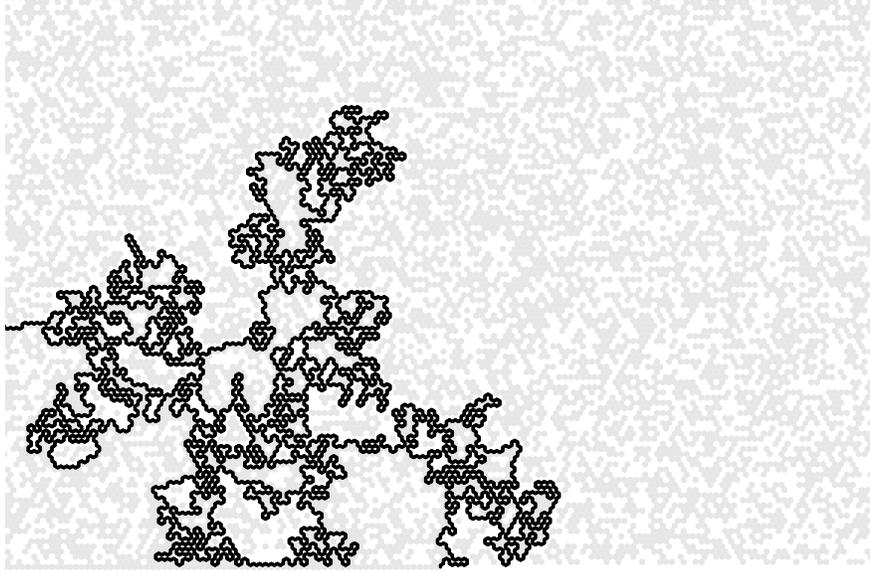}}
\caption{\label{f.longpcurve}The exploration process, proved to converge to $SLE_6$ (see Chapter 10)}
\end{figure}

\section {Chordal versus Radial}

The natural time-parametrization in the previous setup goes as follows:
Let $g_t$ denote the conformal map from $\H \setminus \gamma [0,t]$
onto $\H$ that is normalized at infinity in the sense that when $z \to \infty$,
$$
g_t (z) = z + \frac {a_t}{z} + o (1/z).
$$
It is easy to see that $a_t$ is positive, increasing and that it is natural 
to parametrize $g_t$ in such a way that $a_t$ is a multiple of $t$ (since the $a_t$
terms add up when one composes two such conformal maps). It is natural 
to choose $a_t = 2t$ (this is consistent with the chosen parametrization 
in the radial case).
Then, define $w_t = g_t ( \gamma_t)$, and observe that 
\begin {equation}
\label {e.chordal}
\partial_t g_t (z) = \frac {2}{g_t(z) - w_t }.
\end {equation} 
Hence, just as in the radial case, we observe that it is possible
to recover $\gamma$ using $w$, and that the only 
choice for $w$ that is consistent with the ``Markovian property''
is to take 
$w_t = \sqrt {\kappa} B_t$, where $B$ is ordinary 
one-dimensional Brownian motion.

Hence, one is lead to the following definition: 
Let $w_t = \sqrt {\kappa} B_t$, and define for all $z \in \overline \H$, the solution 
$g_t (z)$ of (\ref {e.chordal}) up to the (possibly infinite) time $T(z)$ 
at which $g_t (z)$ hits $w_t$. Then, define
$$ 
H_t = \{ z \in \H \ : \ T(z) > t \}
$$
and 
$$
K_t = \{ z \in \overline \H \ : \ T(z) \le t \}.
$$
Then, $g_t$ is the normalized conformal map from $H_t$ onto $\H$. 
We call $(K_t, t \ge 0)$ the chordal $SLE_\kappa$ in the upper half-plane.

It turns out that radial and chordal SLE's are rather closely related:
Consider for instance, the conformal image of radial $SLE_\kappa$ under the 
map that maps $\U$ onto $\H$, $1$ to $0$ and $0$ to $i$.
Consider both this process and chordal $SLE_\kappa$ up to their first hitting
of the circle of radius $1/2$ around zero say. Then, the laws of these
two processes are absolutely continuous with respect to each other \cite {LSW2}.
This justifies a posteriori the choice of time-parametrization in the
chordal 
case.

\section {Conclusion}

We have seen that if one considers a discrete model of random curves (or interfaces) that 
combine the two important features:
\begin {itemize}
\item
The Markovian type property in the discrete setting,
\item
Conformal invariance in the limit when the mesh of the lattice goes to zero,
\end {itemize}
then the good way to construct the possible candidates for the scaling limit of these
curves is to encode them via the corresponding conformal mappings. Then, these 
(random) conformal
mappings are themselves obtained by iterations of identically distributed random conformal 
maps. Loewner's theory shows that such families of conformal maps are themselves 
encoded by a one-dimensional function. If one knows this one-dimensional function, one can recover the family of conformal maps, and therefore also the two-dimensional curve.
The one-dimensional random function that generates the scaling limits of the discrete 
models must necessarily 
be a one-dimensional Brownian motion. The corresponding random two-dimensional curves 
are SLE processes.

\subsection* {Bibliographical comments}

Most of the intuition about how to define radial and 
chordal SLE (with LERW as a guide) was already 
present in the introduction of Oded Schramm's first paper \cite {S1} on SLE 
that he released in March 1999.
Our presentation of Lemma~\ref {l.lerw} is borrowed from Lawler
\cite {LLERW}, but there are other proofs of it
(it is for instance closely related to  
Wilson's algorithm \cite {WiD} that will be discussed in Chapter 9).

\chapter {Loewner chains}

This chapter does 
contain  background 
 material on conformal maps and on 
Loewner's equation
(no really new results will be presented here).
The setup is deterministic in this Chapter. SLE will be introduced in 
the next Chapter.

\section {Measuring the size of subsets of the half-plane}

We  
study increasing ``continuously growing''
compact subsets $(K_t, t \ge 0)$ of the upper half-plane. It will
turn out to be important to choose the good time-parametrization.  
We want to find the natural way to measure the
size $a(K)$ of a compact set $K$ and we will then choose
the
time-parametrization in such a way that $a(K_t) = t$.
We will use the following definition 
throughout the paper.

\medbreak
\noindent
{\bf Definition.}
We say that a compact subset $K$ of the closed
upper half-plane $\overline \H$,  
 such that 
$H := \H \setminus K$ is simply connected, is a 
hull.

\medbreak
Riemann's
mapping theorem asserts that there exist conformal maps
$\Phi$ from $H$ onto $\H$ with $\Phi (\infty) = \infty$.
Actually, 
 if $\Phi$ is such a map, the family of maps 
$b\Phi + b'$ for real $b'$ and positive $b$ 
is exactly the family of conformal maps from $H$
onto $\H$ that fix infinity.

Note that since $K$ is compact, the mapping 
$\Psi  : z \mapsto 1/ \Phi (1/z)$ is well-defined
on a neighbourhood of $0$ in $\H$. 
It is possible to 
extend this map $\Psi$ to a whole neighbourhood
 of $0$ in the plane by reflection
along the real axis (this is usually called Schwarz reflection) 
and to check that this extension is analytic. This 
implies that $\Phi$ can be expanded near infinity: There
exist $b_{1}, b_0, b_{-1}, \ldots $, such that 
$$
\Phi (z) 
= b_1 z + b_0 + b_{-1} z^{-1} + \cdots + b_{-n} z^{-n} + o (z^{-n})$$
when $z \to \infty $ in $\H$. Furthermore, since $\Phi$
preserves the
real axis near infinity, all coefficients $b_j$ are real.

Hence, for each $K$, there exists a unique conformal map $\Phi=\Phi_K$
from $H = \H \setminus K$ onto $\H$ 
such that:
$$
\Phi (z) = z + 0 + o (1/z) 
\hbox { when } 
z \to \infty.
$$
This is sometimes called the hydrodynamical normalization.
 In particular, there exists a real
 $a= a(K)$ such that 
$$
\Phi (z) = z+ \frac {2a }{z} + o (1/z) \hbox { when } z \to \infty.
$$
This number $a(K)$ is a way to measure the size of $K$. In a way, it tells
``how big $K$ is in $\H$, seen from infinity''.
It may a priori not be clear that $a$ is a non-negative increasing
function of the set $K$. There is 
a  simple probabilistic
interpretation of $a(K)$ that immediately implies these
facts: Suppose that 
$Z = X+iY$ is a complex Brownian motion started from
$Z_0 = iy$ (for some large $y$, so that $Z_0 \notin K$)
and stopped at its first exit time $\tau$ of $H$. The 
expansion
$\Phi (z) = z + o(1)$
near
infinity shows that $\Im ( \Phi(z) - z)$  
is a bounded harmonic function in $H$. The martingale
stopping theorem therefore shows that 
$$
\expect 
[ \Im ( \Phi (Z_\tau)) - Y_\tau ] = \Im ( \Phi (iy) - iy) 
= \frac {2a }{iy} + o (1/y).
$$
But $\Phi (Z_\tau)$ is real because of the definition of $\tau$.
Therefore
$$
2a= \lim_{y \to +\infty} y \ \expect [ \Im (Y_\tau)] 
.$$
In particular, $a \ge 0$.

One can also view $a$ as a function of the normalized conformal map $\Phi_K$ instead
of $K$. The chain rule for Taylor expansions then immediately shows 
that 
$$ a( \Phi^1 \circ \Phi^2 ) = a ( \Phi^1) + a (\Phi^2) 
$$
for any two normalized maps $\Phi^1$ and $\Phi^2$. In particular, this 
readily implies that $a(K) \le a(K')$ if $K \subset K'$ (because 
there exists a normalized conformal map from 
$\H \setminus \Phi_K ( K' \setminus K)$ onto $\H$).

Let us now observe two simple facts:
\begin {itemize}
\item
If $\lambda >0$, then $a ( \lambda K) = \lambda^2 a(K)$. This is 
simply due to the fact that 
$$
\Phi (z / \lambda) = \frac {z}{\lambda}
+ \frac {2 a(K) \lambda }{z} + o(\lambda /z)
$$
so that
\begin {equation}
\label {e.sc}
\Phi_{\lambda K} (z) = \lambda \Phi_K (z / \lambda) 
= z + \frac {2 a(K) \lambda^2}{z} + o \lambda (\lambda/z)
\end {equation}
when $z \to \infty$.

\item
When $K$ is the vertical slit $[0, iy]$, then 
$$
\Phi_K (z) = \sqrt { z^2 + y^2 }.
$$
In particular, we see that $a([0, iy]) = y^2 / 4$. Note that 
if $y$ is very small, the actual diameter of the vertical slit $[0, iy]$
is much larger than $a([0, iy])$.
\end {itemize}

Equation (\ref {e.sc}) shows that for all $K$ such that 
$a(K) =1$, one has $a(\sqrt {\lambda} K ) = \lambda$
and
\begin {equation}
\label {e.expa}
\lim_{\lambda \to 0}
\frac {\Phi_{  \sqrt \lambda K } (z) - \Phi_{\{0\}}  (z)}{\lambda}  
= 
\lim_{\lambda \to 0}
\frac {\Phi_{  \sqrt \lambda K } (z) - z}{\lambda}  
= 
\frac {2}{z}.
\end {equation}
Actually, it is not very difficult to prove that for all 
given $r$, there exists $C>0$ such that  
this convergence takes place uniformly over all $K$ of radius 
smaller than $r$ and $|z| >Cr$.
See Lemma 2.7 in \cite {LSW1}.

\section {Loewner chains}

Suppose that a continuous real function $w_t$ with $w_0=0$
is given. For each $z \in \overline \H$, define the function $g_t (z)$ as 
the solution to the ODE
\begin {equation}
\label {ode2}
\partial_t g_t (z) = \frac {2} {g_t (z) - w_t }
\end {equation}
with $g_0 (z) = z$. This is well-defined as long as $g_t (z) - w_t$
does not hit $0$, i.e., for all
$t < T(z)$, where
$$
T(z) := \sup \{ t \ge 0 \ : \ \min_{s \in [0,t] } | g_s (z) - w_s | > 0 
\} . $$
We define
\begin {eqnarray*} K_t &:= &\{ z \in \overline \H \ : \ T(z) \le t \}
\\
H_t &:= & \H \setminus K_t .
\end {eqnarray*}
Note for instance that if $w_t = 0$ for all $t$, then
$$ g_t (z) = \sqrt { z^2 + 4t }
$$
and
$K_t = [0, 2i \sqrt {t} ]$.  

It is very easy to check that $g_t$ is a bijection from $H_t$
onto 
$\H$ (in order to see that it is surjective, one can just look at the
ODE ``backwards in time'' to find which point $z$ is such that $g_t (z) =y$).
Moreover $K_t$ is bounded (because $w$ is continuous and 
bounded on $[0,t]$) and $H_t$ has 
a unique connected component (because $g_t^{-1}$ 
is continuous).  
Standard arguments from the theory of ordinary 
differential equations can be applied to check that 
$g_t$ is analytic and that one can formally differentiate the
ODE with respect to $z$, so that 
$$
\partial_t g_t'(z)
= \frac {-2 g_t'(z)}{(g_t(z) -w_t)^2}.
$$
So, $g_t$ is a conformal map
from $H_t$ onto $\H$.

Note also that $| \partial_t g_t (z)|$ is uniformly bounded when $z$ is large and
$t$ belongs to a given finite interval $[0, t_0]$. Hence, 
it follows that 
$g_t (z) = z + O (1)$ near infinity and uniformly over $t \in [0, t_0]$. Hence
(using the ODE yet again), 
$\partial_t g_t (z) = 2/z + o (1/z)$ uniformly over $t \in [0, t_0]$ so that 
finally,
for each $t$,
$$
g_t (z) = z + \frac {2t}{z} + o (1/z)
$$
when $z \to \infty$.
In other words, $a(K_t) =t$. The family 
$(K_t, t \ge 0)$ is called 
the Loewner chain associated to  
the driving function $(w_t, t \ge 0)$.

Loewner's original motivation was to control the 
behaviour of the coefficients
of the Taylor expansion of conformal maps and for this
goal,
it is sufficient to consider smooth slit domains (see e.g., 
\cite {A2, Dur}).
For this reason,   the following question was only addressed  
later (see \cite {Ku}): If the continuous function $(w_t , t \ge0)$ 
is given, what can be said about the family of 
compact sets $(K_t, t \ge 0)$?

In the introduction, we started with a continuous curve $\gamma$, then 
using $\gamma$, we defined $H_t$, the conformal maps $g_t$,
the function $w_t$ and argued that one could recover $\gamma$ from
$w_t$, using the fact that we a priori knew that $g_t^{-1}$
extends continuously to $w_t \in \partial \H$ and that 
$g_t^{-1} (w_t)$ was well-defined (and equal to $\gamma_t$)
because $\gamma$ is 
a continuous curve.
But if one starts with a general continuous function $w_t$, then
it can in fact happen that $g_t^{-1}$ does not
extend continuously to $w_t$. 

Before making general considerations, let us
exhibit a simple example to show that $(K_t, t \ge 0)$
does not need to be a simple curve.
For $\theta \in [0, \pi)$, let 
$\eta (\theta) = \exp ( i \theta) -1 $. Define $t(\theta)
= a ( \eta [0, \theta])$ the ``size'' of the arc $\eta [0, \theta]$.
Finally, define the reparametrization $\gamma$ of $\eta$ in 
such a way that $a ( \gamma [0, t] ) = t$.
$\gamma$ is defined for all $t< T := \lim_{\theta \to \pi-}
a ( \eta [0, \theta])$. 
It is simple to see that there exists a continuous function 
$(w_t, t <T)$ such that the normalized conformal maps $g_t$ 
from $\H \setminus \gamma [0,t]$ onto $\H$ satisfy the 
equation 
(\ref {ode2}). Furthermore, when $t \to T-$, $w_t$ converges 
to a finite limit $w_T$.  
At time $T$, the curve $\gamma [0,T]$ disconnects
the inside of the semi-circle from the outside.
Just before $T$, because $g_t$ is normalized ``from infinity'',
the inside of the semi-circle
is mapped onto a small region which is very close to $w_t
= g_t ( \gamma_t)$.  
When $t \to T-$, all points inside the semi-circle 
are hitting $w_T$. In other words, 
$
K_T$ is the whole semi-disc,  $H_T$ is the 
complement of the semi-disc, and $g_T$ is 
the normalized map from the simply connected 
domain $H_T$ onto $\H$. 
  
Let us now give a couple of general definitions:
\begin {itemize}
\item
We say that $(K_t, t \ge 0)$ is a simple curve if there exists a 
simple continuous curve $\gamma$ such that $K_t = \gamma [0,t]$.

\item
We say that $(K_t, t \ge 0)$ is generated by a curve if there exists a 
continuous curve $\gamma$ with no self-crossings, such that for all $t \ge 0$,
$H_t = \H \setminus K_t$ is the 
unbounded connected component of $\H \setminus 
\gamma [0,t]$. In other words, 
$K_t$ is the union of $\gamma$ and of the inside of the loops that $\gamma$ 
creates.

\item 
We say that $(K_t, t \ge 0)$ is pathological if it is not generated by a curve.
\end{itemize}

In each of these three cases, 
one can find (deterministic) continuous functions $w_t$ such
that the family $(K_t, t \ge 0)$ that it constructs falls into this 
category: For the first case, 
consider for instance $w_t=0$ as before, for the second case, one can use 
the example with the semi-circle.
For the more intricate third case, let us 
mention the following example 
(due to Don Marshal and Steffen Rohde, see \cite {MR}): Let $\gamma$ denote a simple curve in 
$\H$ started from 
$\gamma_0 = 0$ that 
spirals clockwise around the segment $[i, 2i]$ an infinite 
number of times, and then unwinds itself. Then at the ``time'' 
at which it winds around the segment an infinite number of times, $\gamma$
is not continuous i.e. $K_t \setminus K_{t-}$ is the whole
segment. However, this 
Loewner chain corresponds to a continuous function $w_t$.
Such pathologies could arise at any scale.

We now characterize the families $(K_t, t \ge 0)$ 
of compact sets 
that are Loewner chains:

\begin {prop}
\label {p.charact}
The following two conditions are equivalent:
\begin {enumerate}
\item
$(K_t, t \ge 0)$ is a Loewner chain associated to 
a continuous driving function $(w_t, t \ge 0)$.
\item
For all $t \ge 0$, $a(K_t) = t$, and for all $T>0$,
and $\eps >0$, there exists $\delta >0$ such that 
for all $t \le T$, there exists a bounded connected
set $S \subset \H \setminus K_t$
with diameter not larger than $\eps$ such that $S$ 
disconnects $K_{t+\delta} \setminus K_t$
from infinity in $\H \setminus K_t$. 
\end {enumerate}
\end {prop}

\noindent
{\bf Sketch of the proof.}
Let us now 
 prove that 2. implies 1. (the fact that 1. implies 2.
is very easy):   
2. implies that for all $t \ge 0$, the diameter 
of the sets 
$g_{t} ( K_{t+\delta} \setminus K_t ) $ decrease 
towards $0$ when $\delta \to 0$. Hence, one can 
simply define $w_t$ by
$$\{ w_t\}  = \lim_{\delta \to 0} \overline { 
g_t ( K_{t+ \delta} \setminus K_t )}.
$$
Then, one uses 2. to show  that $t \mapsto w_t$ is uniformly continuous.
It then only remains to check that indeed
$$
\lim_{\delta \to 0}
\frac { g_{t+\delta} (z) -g_t (z) }{\delta}
 = \frac {2}{g_t (z) - w_t}.
$$
This is achieved by applying the uniform version of (\ref {e.expa}).
\qed 

\medbreak

Suppose now that $K_t$ is the Loewner chain 
$$ K_t  = [ 0, c \sqrt {t}]$$
for some $c = c(\theta) \exp ( i \theta) \in \H$. Here, $\theta \not= 0$ is given, and 
then the positive real $c (\theta)$ is chosen in such a way that $a(K_1) = 1 $.
Scaling immediately shows that $a(K_t) = t $ for all $t >0$, so that 
there exists therefore a continuous driving function $w$ that generates these
slits.
Again, scaling (because $K_{\lambda t} = \sqrt {\lambda} K_t$)
shows that necessarily, this function $w$ must be of the type 
$$ w_t = c_1 \sqrt {t}$$
for some real constant $c_1 = c_1 (\theta)$.
Let $g_t^{\theta}$ denote the corresponding family of conformal maps.

Let us now choose a new driving function $w$ as follows: $w_t = 0$
when $t < 1$ and for $t \ge 1$:
$$ w_t = c_1 \sqrt {t-1}.$$
When $t < 1$, then $K_t$ is just the straight slit.
 In particular, $g_1 (z) = \sqrt { z^2 + 4 }$.
When $t > 1$, then  $K_t$ is obtained by mapping the angled slit 
$K_{t-1}^\theta$ back by $g_1^{-1}$.
In particular, we see that the curve $\gamma$ generated by this function 
$w$ is not differentiable at $t=1$.
This is one simple hint to the fact that H\"older-$1/2$ regularity may be 
critical (note that at $t=1$, $w$ is just H\"older $1/2$).

The general relation between smoothness of the driving function and 
regularity of the slit has also recently been investigated (in the
deterministic setting) by Marshall-Rohde \cite {MR}.
In this paper, it is shown that H\"older-$1/2$ 
is in a sense a ``critical regularity'' for the 
driving function  $w_t$: Loosely speaking (their results are more precise than 
that), if $w$ is better than H\"older-$1/2$, then it defines a ``smooth'' (in some appropriate sense) slit, 
but nasty ``pathological'' phenomena can occur for H\"older-$1/2$ driving functions.
See \cite {MR} and the references therein.

\bibc
For general background on complex analysis, Riemann's mapping
theorem, there are plenty of 
 good references, see for instance \cite {A1, Ru}.
Loewner introduced his equation (in the radial setting)  
in 1923 \cite {Lo}.
For general information about Loewner's equation, and in particular how
Loewner used it to prove that $|a_3| \le 3$ for univalent functions 
$z+ \sum_{n \ge 2} a_n z^n$ on $\U$ as well
as other applications, see for instance \cite {A2, Dur}.
For how it is used in de Branges' proof of the Bieberbach conjecture, 
a good self-contained reference is Hayman's book \cite {Hay}.
For basics on hypergeometric functions, see e.g., \cite {Le}.

Proposition~\ref {p.charact} 
is derived in \cite {LSW1}, see also 
\cite {Po}.
Carleson and Makarov \cite {CM1,CM2} have used Loewner's (radial)
 equation in the context
of Diffusion Limited Aggregation.

\chapter {Chordal SLE}

\section {Definition}

Chordal $SLE_\kappa$ is the Loewner chain 
$(K_t, t \ge 0)$ that is obtained when 
the driving function $$
w_t = W_t := \sqrt {\kappa}  B_t 
$$
is $\sqrt \kappa$ times a  standard real-valued Brownian motion
$(B_t, t \ge 0)$ with
$B_0=0$.
Let us now list a couple of consequences of the simple properties
of Brownian motion:

\begin {itemize}
\item
Brownian motion is a strong Markov process with independent increments. This
implies that for any stopping time $T$ (with respect to the
natural filtration $({\cal F}_t, t \ge 0)$ of $B$), the
process
$$
( g_{T+t} ( K_{T+t} \setminus K_T ) - W_T , t \ge 0)
$$
is independent of ${\cal F}_T$ and that its law is identical to 
that of $(K_t, t \ge 0)$. Note that one has to shift by $W_T$ in order to 
obtain a process starting at the origin.

\item
Brownian motion is scale-invariant:   
For each $\lambda >0$, the 
process $W_t^\lambda := 
W_{\lambda t} / \sqrt {\lambda}, t \ge 0$ has the same law 
than $W$. But 
$$
\partial_t ( g_{\lambda t} (  \sqrt {\lambda} z )) 
= 
\frac {2 \lambda }{g_t ( \sqrt {\lambda}  z) - W_{\lambda t} }.
$$
In particular, if 
$$
g_t^\lambda (z) := g_{\lambda t } ( z \sqrt {\lambda}) / \sqrt {\lambda},$$
then 
$$
\partial_t g_t^\lambda (z) 
= \frac {2} { g_t^\lambda (z)  - W_t^\lambda }
$$
and $g_0^\lambda (z) = z$. In other words, $(K_{\lambda t}, t \ge 0)$ 
and 
$( \sqrt {\lambda} K_{t}, t \ge 0)$ have the same law: Chordal $SLE_\kappa$
is scale-invariant.

\item
Brownian motion is symmetric ($W$ and $-W$ have the
same law). Hence, the law of $(K_t , t \ge 0)$ is symmetric 
with respect to the imaginary axis.

\end {itemize}

It is actually possible to prove the following result:

\begin {prop}
For all $\kappa \ge 0$, chordal $SLE_\kappa$ is
almost surely not
pathological. 
When $\kappa \le 4$, it is a.s. a simple curve $\gamma$, when
$\kappa > 4$, it is a.s. generated by a (non-simple) curve
$\gamma$.
\end {prop}

This result is due to Rohde-Schramm \cite {RS} (see 
\cite {LSWlesl} for the critical case $\kappa =8$). It is not
an easy result, especially for the values $\kappa >4$.
Actually, while this fact is important 
and useful in order to understand 
heuristically the behaviour and the properties of $SLE_\kappa$,
it turns out that one can derive many of them without 
knowing that the $SLE_\kappa$ is generated by a continuous 
curve.
We therefore omit the proof in these lectures, and we will 
call $(K_t, t \ge 0)$ the SLE. 
In some cases that we will 
focus on ($\kappa=2,8/3,6,8$), the fact that
$SLE_\kappa$ is a.s. generated
by a curve will actually follow from other considerations.

It is however easy to see that $\kappa=4$ is a critical value:
Consider chordal $SLE_\kappa$, and define 
$$
X_t = \frac { g_t (1) - W_t}{\sqrt {\kappa}} 
.$$
Note that $X$ hits zero if and only if the chordal SLE absorbs the boundary point 
$1$. 
But $X$ satisfies  
\begin {equation}
\label {bessel}
dX_t = dB_t + \frac { 2 }{\kappa X_t} dt .
\end {equation}
It is a $1+ (  4 / \kappa) $ dimensional Bessel process, and it is well-known 
(see e.g. \cite {RY}) that such a process a.s. hits zero if and only if $\kappa > 4$.
This can for instance be viewed as a consequence of the fact that 
if $X$ is a Bessel process of dimension $d$ started away from zero, then 
if $d \not= 2$, $X^{2-d}$ is a local martingale, and when $d=2$, 
$\log X$ is a local martingale.

It follows that:
\begin {prop}
\label {p.phases}
\begin {itemize}
\item
If $\kappa \le 4$, then almost surely $\cup_{t \ge 0} K_t 
\cap \R = \{ 0 \}$.
\item
If $\kappa > 4$, then almost surely, $\R \subset \cup_{t \ge 0} K_t $.
\end {itemize}
\end {prop}

 Assuming that the SLE is generated by a curve, this readily shows that
the SLE curve is simple if and only if $\kappa \le 4$.

If one defines, for all $z \in \H$,
the solution $X_t^z$ to (\ref {bessel}) started from $X_0^z = z / \sqrt \kappa$
(up to the stopping time $T(z)$).
Then, we see that $SLE_\kappa$ can be interpreted in terms
of the flow of a complex Bessel process: 
For each $t>0$, $K_t$ is the set of starting points such that 
$X_t^z$ has hit $0$ before time $t$.

\section {A first computation}

We now compute the probability of some simple
events involving the chordal Schramm-Loewner evolution. 
Suppose that $a<0<c$. Let $\kappa>0$ be fixed. Define the event
$E_{a,c}$ that 
the chordal $SLE_\kappa$ hits $[c, \infty)$ before $(- \infty, a]$.
For the reasons that we just discussed, this makes sense only
if $\kappa > 4$ (otherwise, it never hits these intervals).
The goal of this section is to compute the probability of 
$E_{a,c}$.
The scaling property of chordal SLE shows that this is a function of 
the ration $c/a$ only. We can therefore define
$F= F_\kappa$ on the interval  $(0,1)$ 
by
$$
\P [ E_{a,c} ] = F ( -a / (c-a)).
$$

\begin {prop}
\label {p.cardy}
For all $\kappa >4$ and $z \in (0,1)$,
$$
F(z)  =
 c(\kappa) \int_0^z \frac {du}{u^{4/\kappa} (1-u)^{
4/\kappa}}
$$
where $c(\kappa) = ( \int_0^1 u^{-4/\kappa} (1-u)^{-4/\kappa} du )^{-1}$ is
chosen so that $F(1) =1$.
\end {prop}

Note that this Proposition is in fact a property of the 
real Bessel flow:
$E_{a,c}$ is the event  that $X^{c}$ hits $0$ before
$X^{a}$ does.

\medbreak
\noindent
{\bf Proof.}
Suppose that ${\cal F}_t = \sigma ( B_s, s \le t)$ is the natural
filtration associated to the Brownian motion, and define 
$T_a=T(a)$ and $T_c=T(c)$ as before (the times at which $a$ and $c$ are
respectively absorbed by $K_t$).
For $t< T_a$ and $t< T_c$ respectively, define
$$ 
A_t := g_t (a)  \hbox { and } C_t := g_t (c)
.$$
Suppose that $t< \min (T_a, T_c) $, and define
$$
K_{t,s} = g_t ( K_{t+s} \setminus K_t ) - W_t.
$$
The strong Markov property shows that 
$
(K_{t,s}, s \ge 0)$ is also chordal $SLE_\kappa$, and that it 
is independent from ${\cal F}_t$.
Also, if $t< \min (T_a, T_c)$, 
the event $E_{a,c}$ corresponds to the event that 
$(K_{t, s}, s \ge 0)$ hits $[C_t-W_t, \infty)$ 
before 
$(-\infty, A_t - W_t]$. 
Hence, if $t < \min (T_a, T_c)$,
$$
\P [ E_{a,c} \mid {\cal F}_t ] 
= F \left( \frac { W_t - A_t }{C_t - A_t} \right).
$$
In particular, 
this shows that the right-hand side of the previous 
identity is a (bounded) martingale.
We know that $W_t = \sqrt {\kappa }  B_t$, and
that 
$$ \partial_t A_t = \frac {2}{A_t - W_t}, \ 
\partial_t C_t = \frac {2} {C_t - W_t}.
$$
Hence, if we put $Z_t := (W_t -A_t)/ (C_t - A_t)$, 
stochastic calculus yields
$$
d Z_t  
=
\frac {\sqrt { \kappa}  dB_t}{C_t - A_t}
+ \frac { 2 dt}{ ( C_t - A_t)^2}
\left( \frac {1}{Z_t} - \frac {1}{1- Z_t} \right)
.$$
One can now also introduce the natural time-change
$$
s = s(t) : = \int_0^t \frac {du}{(C_u-A_u)^2} 
$$
and define $\tilde Z$ in such a way that 
$\tilde Z_{s(t)} = Z_{t}$.
Then, 
$$
\tilde Z_s = \sqrt {\kappa} d \tilde B_s +2 \left( 
 \frac {1}{Z_s} - \frac {1}{1- Z_s} \right) ds
$$ 
where $(\tilde B_s, s \ge 0)$ is a standard Brownian 
motion.
 
But $K_t$ hits $(-\infty, a)$ if and only $Z_t$ hits $0$,
and  $K_t$ hits $(c, \infty)$ if and only if $Z_t$
hits $1$.
Hence, 
$F(z)$ is the probability that the diffusion $\tilde Z$
started from $\tilde Z_0 = z$ hits $1$ before $0$.
One can invoke (for instance) the general theory of diffusions 
to argue that the function $F$ is therefore smooth on 
$(0,1)$. It\^o's formula (since $F(\tilde Z_s)$ is a martingale)
then implies that 
\begin {equation}
\label {e.ode1}
\frac {\kappa}{4} F'' (z) +  \left(
 \frac {1}{z} - \frac {1}{1- z} \right)  F' (z) = 0
. \end {equation}
Furthermore, the boundary values of $F$ are 
simple to work out: When $\kappa >4 $, one can see
(for instance comparing $\tilde Z$ with a Bessel
process) that
$$
\lim_{z \to 0} F(z) = 0 
\hbox { and } \lim_{z \to 1} F(z) = 1.
$$
Hence, 
$F$ is the only solution to the ODE (\ref {e.ode1})
with boundary values $F(0) =0$ and $F(1)=1$. 
This immediately proves the Proposition.
\qed

\medbreak

Note that when $z \to 0$,
$$
F(z)  \sim \frac {c(\kappa)}{1- 4/\kappa}
  z^{1- 4/\kappa}
.$$
In particular, for $\kappa=6$,  we get the 
exponent $1/3$.

Exactly in the same way, it is possible (for $\kappa > 4$) to compute the 
probability that chordal $SLE_\kappa$ (started from $0$) hits the 
interval $[a,c]$ before $[c, \infty)$ when $0< a <c$.
This is a function $\tilde F$ of the ratio $a/c$, satisfying a 
linear second-order differential equation, with the boundary conditions
$$
\tilde F (1) =0 
\hbox { and } \tilde F (0) = 1
.$$

\section {Chordal $SLE_\kappa$ in other domains}

Suppose that $D$ is some given
non-empty open simply connected
subset of the complex plane with $D \not= \C$.
We do not impose any regularity 
condition on $\partial D$. Riemann's mapping
theorem shows that there exist 
(many) conformal maps $\Phi$ from 
the upper half-plane $\H$ onto $D$.
Even if the boundary of $\partial D$ is
not smooth, one can define a general notion that 
coincides with that of boundary  
 points when it is smooth: For each $x \in \overline \R$,
we  
say that (if some map $\Phi$ is given) $\Phi(x)$ is a 
prime end of $D$ (see e.g. \cite {P2} for a more precise 
and correct definition). 

Suppose that $O$ and $U$ are two distinct prime ends in $D$.
Then, there exists a conformal map $\Phi$ from $\H$ onto $D$
such that $\Phi (0) =O$ and $\Phi (\infty) = U$.
Actually, this only 
characterizes $\Phi ( \cdot)$ up to a multiplicative 
factor (because 
$\Phi (\lambda \cdot)$ would then also do).

Suppose that $(K_t, t \ge 0)$ is chordal $SLE_\kappa$ 
in $\H$ as defined before.
We define $SLE_\kappa$ in $D$ from $O$ to $U$ as
the image of the process $(K_t, t \ge 0)$
under $\Phi$.
Recall that $\Phi$ is  defined
up to a multiplicative constant.
However, the scaling property of 
$SLE_\kappa$ in $\H$ shows that 
the law of $(\Phi (K_t), t \ge 0)$ is 
invariant (modulo linear time-change) if we replace
$\Phi (\cdot)$ by $\Phi (\lambda \cdot)$.
 
To illustrate this definition, consider the following setup: 
Suppose that $\kappa=6$ and that 
$OAC$ is an equilateral triangle.
Let $\Phi$ denote the conformal map from $\H$ onto the 
triangle defined in such a way that 
$$
\Phi (a) = A, \Phi (0) = O , \Phi (c) = C
$$
where $a<0<c$ are given.
This conformal map can be easily described explicitly 
using the Schwarz-Christoffel transformations
\cite {A1, Ru}.
Note that $U= \Phi( \infty)$ is on the interval $AC$.
It turns out that 
$$ 
\frac {AU}{AC}  = F ( z)
$$
where $z = -a / (c-a)$ and
$F = F_{\kappa=6}$ is precisely the same hypergeometric function 
as in Proposition \ref {p.cardy}.
Hence, the probability that chordal $SLE_6$ from
$O$ to $U$ in the equilateral
triangle $OAC$ hits $AU$ before $UC$ is simply
the ratio $AU/AC$.

\medbreak

Suppose now that $\kappa \in (4, 8)$. Just as for the hypergeometric function 
$F$, the functions $\tilde F$ that were defined at the end of the last subsection
 have a nice interpretation in terms of 
conformal mappings onto triangles:
Consider an isocele triangle ${\cal T} = OAU$ with $OA=AU=1$ and angle $\pi (1- 4 / \kappa)$ 
at the vertices $O$ and $U$. The angle at the vertex $A$ is therefore 
$\pi ( 8/ \kappa - 1)$.
Consider now a chordal $SLE_\kappa$ from $O$ to $U$ in the triangle ${\cal T}$.
Let $X$ denote the random point at which it first hits the segment $AU$. 

\begin {prop}
\label {p.unif}
The law of $X$ is the uniform distribution on $AU$.
\end {prop}

This is a direct consequence of the explicit computation of ${\tilde F}$ and of the 
explicit Schwarz-Christoffel mapping from the upper half-plane onto ${\cal T}$: For 
each $C \in AU$, one can compute the probability that $X \in [AC]$ via the 
function $\tilde F$.
\qed

\medbreak
This gives a first justification to the fact that the only possible 
conformally invariant scaling limit of the critical percolation exploration process
is $SLE_6$ (see more on this in Chapter 10).
Indeed, suppose that the critical percolation exploration 
process is conformally invariant. We have argued in the first chapter that the 
scaling limit is one of the $SLE$s. Suppose that it is $SLE_\kappa$ for 
a given value of $\kappa$, and consider the corresponding triangle ${\cal T}$.
 
Clearly in the discrete case (for a fixed small meshsize), up to the 
first time at which it hits the edge $AU$, the critical exploration process 
from $O$ to $U$ and the critical exploration process from $O$ to $A$ 
in ${\cal T}$ coincide. Hence, the hitting distributions on 
$AU$ for chordal $SLE_\kappa$ from $O$ to $U$ and for chordal 
$SLE_\kappa$ from $O$ to $A$ coincide. In particular, the uniform distribution 
on $AU$ must be invariant under the anti-conformal map from 
${\cal T}$ onto itself that maps $O$ onto itself and interchanges the 
vertices $A$ and $U$. This is only true when the triangle is symmetric 
(i.e. the angles at $U$ and $A$ are identical), in other words 
when $\alpha= \pi /3$ or $\kappa = 6$.

We shall see in the next chapter that indeed, for $SLE_6$, the 
whole paths from $O$ to $A$ and from $O$ to $U$ coincide up to their 
first hitting of $AU$. This is the so-called locality property of $SLE_6$.

\section {Transience}

We conclude this chapter with the following fact
(assuming the fact that the SLE is a.s.  
a simple curve $\gamma_t = K_t \setminus K_{t-}$ for 
$\kappa < 4$). This is also to illustrate the type of
techniques that is used to derive such properties of SLE:

\begin {prop}
\label {p.transience}
For $\kappa<4$, almost surely,
$\lim_{t \to \infty} \gamma_t = \infty$.
\end {prop}
Loosely speaking, the SLE is transient.
Actually (see \cite {RS}), this result is in fact valid for all 
$\kappa$, but the proof is (a little bit) more involved.

\proof
Let $\delta \in (0, 1/4)$, $x > 1$, and suppose that 
$$
t_\delta : = \inf \{ t > 0 \ : \ d( \gamma_t, [1,x] ) \le \delta \}
$$
is finite.
Let $z_\delta = \gamma_{t_\delta}$.
Clearly, $g_{t_\delta} (z_\delta) = W_{t_\delta}$.
Note that 
$g_{t_\delta} (1/2) - W_{t_\delta}$ is (up to a 
multiplicative constant)  the limit when $y \to + \infty$ of 
$y$ times the probability 
that a planar Brownian motion started from $iy$ exits $\H$ in the 
interval $[W_{t_\delta}, g_{t_\delta} (1/2)]$.   
By conformal invariance, this is the same as the limit of $y$
times the probability that a planar Brownian motion started from 
$iy$ exits $H_{t_\delta}$ through the boundary of $H_{t_\delta}$ which is 
``between''    
$z_\delta$ and $1/2$.
But in order to achieve this, the planar Brownian motion 
has in particular to hit the vertical segment joining 
$z_\delta$ to the real line before exiting $\H$.
 This segment has length at most
$\delta$.  
Hence, 
$$
|g_{t_\delta} (1/2) - W_{t_\delta}| \le O (\delta).
$$
On the other hand, $\lim_{t \to \infty}
(g_t(1/2) - W_t) = \infty$ because $\kappa < 4$ (and the
corresponding Bessel process is transient).
It follows that a.s.,
$$
d ( \gamma [ 0, \infty]  , [1,x] )  > 0 .
$$
By the scaling property and monotonicity,
 it follows that almost surely, for all $0<x_1<x_2$, 
the distance $d( \gamma [0, \infty] , [x_1,x_2]) $ is almost surely
strictly positive.

Let $\tau$ denote the hitting time of the unit
circle by the SLE.
Since $\R \cap \gamma [0,\tau] = \{ 0 \}$, 
it follows that $0 \in \partial H_{\tau}$.
For all $\eps>0$, there exists $ 0 < x_1 < x_2$ 
 such that with probability at least 
$1-\eps$  the two images
of $0$ under $g_\tau$ are in $[W_\tau - x_2, W_\tau -x_1] \cup 
[W_\tau + x_1, W_\tau + x_2]$.
It follows from the strong Markov property and from the 
previous result that with probability at least $1- \eps$, 
$$ d ( g_\tau ( \gamma [\tau, \infty)) - W_\tau , [-x_2, -x_1] \cup
[x_1, x_2] ) >0 .
$$
Hence, it follows that in fact, almost surely
$$
d ( 0, \gamma [ \tau, \infty) ) > 0
$$
and the Lemma readily follows (for 
instance using the scaling property once again).
\qed

 \section* {Bibliographical comments}

Again, many of the ideas in this chapter were contained 
or follow readily from Schramm's first paper \cite {S1}.
Rohde-Schramm \cite {RS} have derived various almost 
sure properties of SLE (H\"older boundary, generated by a continuous 
path, transience). Proposition~\ref {p.cardy} is derived (in a more general 
setting) in \cite {LSW1}.
It was Carleson who first noted that Cardy's formula (which 
Cardy predicted for crossing probabilities for critical percolation)
has a simple interpretation in an equilateral triangle. The interpretation of the 
functions ${\tilde F}$ in terms of isocele triangles was pointed out 
by Dub\'edat \cite {Dub}. 
Another justification to the fact that $\kappa = 6$ is the unique possible scaling limit 
of critical percolation exploration processes (for site percolation 
on the triangular lattice, or for bond percolation on the square lattice) uses the fact that 
for these models the probability of existence of a left-right crossing of a square must be $1/2$ (see \cite {S1}).
For references on Bessel processes,
stochastic calculus, see e.g.~\cite {IW,RY}.

\chapter {Chordal SLE and restriction}

\section {Image of SLE under conformal maps}

Suppose now that $(K_t, t \ge 0)$ is chordal $SLE_\kappa$
in the upper half-plane $\H$.

\medbreak
\noindent
{\bf Definition.} We say that a hull $A$ that is at positive
distance of the origin is a Hull (with capital H).
When $A$ is such a Hull, we define $\Phi_A$ the normalized
conformal map from $\H \setminus A$ onto $\H$
as before. We also 
define $\Psi_A$ the conformal map from
$\H \setminus A$ onto $\H$ such that $\Psi (z) \sim z$
when $z \to \infty$ and $\Psi(0)=0$. Note that 
$\Psi (z)= \Phi (z) - \Phi (0)$.

\medbreak
Let $A \subset \overline \H$ denote a Hull.
Define $T = \inf \{ t \ : \ K_t \cap A \not= \emptyset \}$ and for all
$t < T$,
$$ \tilde K_t := \Phi (K_t ).
$$
Let us immediately emphasize that the time-parametrization of 
$K_t$ and therefore also of $\tilde K_t$ is given in terms of 
the ``size'' of $K_t = \Phi^{-1} (\tilde K_t)$ in $\H$
and not in terms of the ``size'' of $\tilde K_t$ itself in $\H$.
One of the goals of this section is to study the 
evolution of $\tilde K_t$ and to compare it with
that of $K_t$.

For $t<T$, we also define the conformal map
$h_t$ from $g_t ( H_t \cap H)$ onto $\H$ (where
 $H = \H \setminus A$).
Note that $h_0 = \Phi$. Since $g_t (A)$ is
at positive distance of $W_t$ for $t<T$, we can define
$$
\tilde W_t = h_t (W_t).
$$
Define finally also the normalized  conformal map $\tilde g_t$ 
from $\Phi (H_t \cap H)$ onto $\H$. Note that (as long as $t<T$), 
$$
h_t \circ g_t = \tilde g_t \circ h_0.
$$
In short, all these maps are normalized, $h_0=\Phi$ removes $A$
and $\tilde g_t$ removes $\tilde K_t$, while $g_t$ removes
$K_t$ and $h_t$ removes $g_t (A)$.

The family $(\tilde K_t, t < T)$ is a ``continuously'' growing family
of  subsets of $\H$ satisfying Proposition~\ref {p.charact} except that 
a time-change is required in order to 
parametrize it as a Loewner chain.
We therefore define
the function
$$
a(t) := a ( A \cup K_t) = a(A) + a ( \tilde K_t).
$$
A simple time-change shows that 
$$
\partial \tilde g_t (z) 
= \frac {2 \partial_t a }{\tilde g_t (z) - \tilde W_t } 
.$$
Hence, in order to understand the evolution of $\tilde K_t$, we
have to understand the evolutions of $\tilde W_t$
and of $a (t)$.
 
The scaling rule $a(\lambda \cdot) = \sqrt {\lambda}
a (\cdot)$ shows that
$$
\partial_t a(t) =
 h_t' (W_t)^2 .
$$
On the other hand,
$$
h_t = \tilde g_t \circ \Phi \circ g_t^{-1} 
$$
and
$$
\partial_t ( g_t^{-1} (z)) 
= - 2 \frac { (g_t^{-1})' (z) } {z - W_t}
$$
so that putting the pieces together, we see that
\begin {equation}
\label {e.h}
\partial_t h_t (z) 
= \frac {2 h_t' (W_t)^2 }{h_t (z) -\tilde W_t }
- \frac {2 h_t' (z) }{z- W_t}.
\end {equation}
Recall that $\tilde W_t = h_t (W_t)$.
The previous formula is valid for all $z \in \H \setminus g_t (A)$.
In fact, one can even extend it to $z= W_t$:
$$
(\partial_t h_t) (W_t)
= \lim_{z \to W_t} 
\left(  \frac {2 h_t' (W_t)^2 }{h_t (z) -\tilde W_t }
- \frac {2 h_t' (z) }{z- W_t}  
\right)
= - 3 h_t'' (W_t)
$$
(note that $h_t$ is smooth near $W_t$ because of Schwarz reflection).
It\^o's formula (this is not the classical formula since 
$h_t$ is random, but it is 
adapted with respect to the 
filtration of $W_t$, it is $C^1$ with respect to $t$, so that It\^o's 
formula still holds, see e.g., exercise IV.3.12 in \cite {RY})
can be applied:
$$
d  \tilde W_t =
(\partial_t h_t ) (W_t) dt + h_t' (W_t) dW_t + \frac {\kappa}{2} h_t'' (W_t) dt.
$$
Hence,
$$
d \tilde W_t = h_t'(W_t) dW_t + [ ( \kappa/2) - 3 ] h_t'' (W_t).
$$
Clearly, the value $\kappa= 6$ will play a special role here.
The next section is devoted to this case.

\section {Locality for $SLE_6$}

Throughout this section, we will assume that $\kappa=6$.
Then, 
$$\tilde W_t = \int_0^t h_s'(W_s) dW_s.$$
Recall also that 
 $a_t  - a_0 = \int_0^t h_s'(W_s)^2 ds = \langle \tilde W 
\rangle_t$. Hence, if we define
$(\hat W_a, a \ge 0)$ in such a way that 
$$\tilde W_t = \hat W_{a(t) - a(0)},$$
 then 
$\hat W - \hat W_0$ and $W$ have the same law.
If we define $\hat g_a$ in such a way
that $\tilde g_t = \hat g_{a(t)}$,
then
$$
\partial_a \hat g_a (z) = \frac { 2 }{ \hat g_a (z) - \hat W_a }.
$$
Hence, modulo time-change, the evolution of $\tilde K_t - \hat W_0$ up to $t=T$
is that of chordal $SLE_6$.
Suppose that $\tilde T$ is the first time at which $K_t$ hits
$\Phi (\partial A)$.
We have just proved $SLE_6$'s locality property:
 
\begin {theorem}
\label {t.loc}
Modulo time-reparametrization, 
the processes $ (\tilde K_t - \Phi (0), t < T ) $ and 
$ ( K_t, t < \tilde T)$
have the same law.
\end {theorem} 

We now discuss some consequences of this result.
Suppose first
that 
$$A=A_\eps=\{ e^{i\theta} \ :  \ \theta \in [0, \pi - \eps] \}.
$$
Recall that $\Phi=\Phi_\eps$ is the normalized map from 
$\H \setminus A$ onto $\H$.
Let 
$$\psi_\eps (z)= \frac {\Phi_\eps (z) }{ \Phi_{\eps}' (0)}.$$
It is easy to see that when $\eps \to 0$, the mappings
$\psi_\eps$ converge uniformly on any set $
V_\delta := \{ z \in \H \ : \ |z| < 1-\delta\}$
towards the conformal map $\psi $ from $V:=\{ z \in \H \ : \ |z| <1 \}$
onto $\H$ such that $\psi(0)=0$, $\psi'(0) = 1$ and $\psi(-1) = \infty$.
Theorem \ref {t.loc} shows that for each $\eps>0$, the law
of the process $\psi_\eps (K_t)$ up to its hitting time of $\psi_\eps
(A_\eps)$
is a time-change of chordal $SLE_6$. In particular, letting 
$\eps \to 0$ for each fixed $\delta >0$ shows readily 
that:

\begin {corollary}
Let $(K_t, t \ge 0)$ denote the law of chordal $SLE_6$
from $0$ to $-1$ in $V$. Let $T$ the first time at which
$K_t$ hits the unit circle. Then, the law of 
$(K_t, t<T)$ is identical (modulo time-change)
to that of chordal $SLE_6$ in $\H$ (from $0$ to $\infty$)
up to its first hitting time of the unit circle.
\end {corollary}

The same reasoning can be applied to
$\{ e^{i \theta} \ : \ \theta \in [\eps , \pi] \}$ 
instead of $A_\eps$.
It shows that the law described in the corollary 
is also identical to that of chordal $SLE_6$ 
from $0$ to $+1$ in $V$ (up to the hitting time of the
unit circle).
By mapping the set $V$ onto any other simply connected domain,
we get the following splitting property:

\begin {corollary}
\label {c.splitting}
Let $D \subset \H$ denote a simply connected subset of $\H$
such that the boundary of $\partial D$ is a continuous Jordan curve.
Let $a, b, b'$ denote three distinct points on $\partial D$
and call $\partial$ the connected component of $\partial D
\setminus \{ b,b'\}$ that does not contain $a$.
Then: up to their first hitting times of $\partial$ and
modulo time-change, the laws of chordal $SLE_6$ 
from $a$ to $b$ and from $a$ to $b'$ in $D$ are identical.
\end {corollary}



 
Note that these properties of chordal $SLE_6$ 
are not surprising if one thinks of $SLE_6$ as the
scaling limit of critical percolation interfaces.
They generalize the properties of hitting probabilities 
for $SLE_6$ that 
we derived in the previous chapter.

\section {Restriction for $SLE_{8/3}$}

We now apply the same technique as in the first subsection to 
understand how $h_t' (W_t)$ evolves.
Recall that $h_t$ is smooth in the neighbourhood of $W_t$
by Schwarz reflection.
Hence $h_t'(W_t)$ is a positive real (as long as $t<T$).
Differentiating Equation (\ref {e.h}) with respect to 
$z$ (this is licit as long as $t<T$) gives
$$
\partial_t h_t' (z) = \frac {-2 h_t'(W_t)^2 h_t' (z)}{(h_t (z) - \tilde W_t)^2}
+ \frac {2 h_t' (z)}{(z- W_t)^2} - \frac {2 h_t'' (z) }{z- W_t}.
$$
If we take the limit when $z \to W_t$, we get that
$$
( \partial_t h_t') (W_t) = \frac {h_t''(W_t)^2}{2 h_t' (W_t)} - \frac {4}{3} h_t'''(W_t)
.$$
Hence, 
It\^o's formula (in its random version as before) shows that
$$
d [ h_t'(W_t) ] 
= 
h_t''(W_t) dW_t + \left[ \frac {h_t''(W_t)^2}{2 h_t'(W_t)} 
+ ( \kappa/2 - 4/3) h_t''' (W_t) \right] dt.
$$
This time, it is the value $\kappa =8/3$ that plays a special role.
Let us in this section from now on suppose that $\kappa = 8/3$.
Then, we see that 
$$
d [ h_t'(W_t)^{5/8} ] 
= \frac {5 h_t'' (W_t)}{8 h_t'(W_t)^{3/8}} dW_t.
$$
The important feature is that the drift term disappear so that:
$(h_t' (W_t)^{5/8} , t < T)$
is a local martingale.
This has the following important consequence: 

\begin {prop}
\label {p.restriction}
Consider chordal $SLE_{8/3}$ 
in $\H$. Then, 
for any Hull $A$,
$$
\P [ \forall t \ge 0 , \ K_t \cap A = \emptyset  ] = \Phi_A '(0)^{5/8}
.$$
\end {prop}

\proof
The quantity $M_t := h_t' (W_t)^{5/8}$ is a local martingale.
Recall that $h_t$ is a normalized map from a subset of $\H$
onto $\H$. Hence, for all $t<T$, $M_t \le 1$
and $M$ is a bounded martingale. We have to understand the
behaviour of $M_t$ when $t \to T$ in the two
cases $T<\infty$ and $T= \infty$.
When $T=\infty$, one can use the transience of the SLE: 
Define for each $R$, the hitting time $\tau_R$ of the 
circle of radius $R$. Then, simple considerations using harmonic 
measure for instance show that 
$$ 
\lim_{R \to \infty} h_{\tau_R}' (W_{\tau_R})
= 1.
$$
In the case where $T<\infty$, one can for instance first
approximate $A$ by a Hull with a smooth boundary, and show 
that in this case, $\lim_{t \to T} h_t'(W_t)=0$ for any
path $\gamma$ in the upper half-plane that hits $A$ away
from the real line. See \cite {LSWrest} for details.

Finally, since $M_t$ converges in $L^1$ and almost surely
when $t \to T$, we get
that $\P [ T= \infty ] = \expect [  M_T ] = E [M_0 ] = \Phi' (0)^{5/8}$.
\qed
\medbreak

Let us now define the random set 
$$
K_\infty = \cup_{t >0} K_t 
.$$ 

\begin {corollary}
\label {c.rest}
Suppose that $A_0$ is a Hull, then
the conditional law of $K_\infty$ given $K_\infty \cap A_0 =
\emptyset$ is identical to the law of $\Psi^{-1}_{A_0}
(K_\infty)$.
\end {corollary}
 
\proof
Note that  $K_\infty$ 
is a closed set because of the transience of 
$(K_t, t \ge 0)$.
The law of such a 
random set is characterized by the value of 
$\P [ K_\infty \cap A = \emptyset ]$ 
for all Hulls $A$ (this set of events is a generating
$\pi$-system of the $\sigma$-field on which we define $K_\infty$).
Suppose now that the Hull $A_0$ is 
fixed. By Proposition~\ref{p.restriction}, 
$K_\infty$ avoids $A_0$ with positive probability.
Suppose that $A$ is another Hull.  
Then
\begin {eqnarray*}
\lefteqn {\P [ \Psi_{A_0} ( K_\infty) \cap A = \emptyset |
K_\infty \cap A_0 = \emptyset ]
} 
\\
&=&
\frac { \P [ K_\infty \cap (\H \setminus (\Psi_{A_0}^{-1} \circ \Psi_A^{-1} (\H)) = \emptyset  ] }
{\P [ K_\infty \cap A_0 = \emptyset ] }
\\
&=&\left( \frac { \Psi_{A_0}'(0) \Psi_A'(0) }{\Psi_{A_0}'(0)}
\right)^{5/8}
\\
&=&
\P [ K_\infty \cap A = \emptyset ].
\end {eqnarray*}
Since this is true for all Hull $A$, 
it follows that the 
 the law
of $\Psi_{A_0} (K_\infty)$ given $\{ K_\infty \cap A_0 = \emptyset \}$
is identical to the law of $K_\infty$. 
\qed
\medbreak

\begin{figure}
\centerline{\includegraphics*[height=4in]{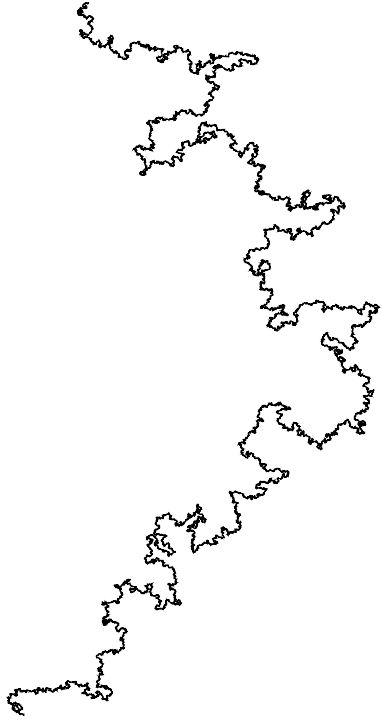}}
\caption{\label{f.saw}
Sample of the beginning of a half-plane walk (conjectured to 
converge to chordal $SLE_{8/3}$).}
\end{figure}

This striking property of $SLE_{8/3}$ has many nice consequences. 
It will enable us to relate it to the Brownian frontier in the 
next chapter. It also shows that it is the natural candidate for 
the scaling limit of planar self-avoiding walks. More precisely, one 
can show that when $n \to \infty$, the uniform measure 
on self-avoiding walks of 
length $n$ in the upper
half-plane  $\N \times \Z$ started from the origin 
converges to a law of infinite self-avoiding walks.
The conjecture is that the scaling limit of this infinite 
self-avoiding walk is $SLE_{8/3}$. See~\cite {LSWSAW} for more
on this. Note that there exist
algorithms to simulate
 half-plane self-avoiding walks (see \cite {Kennalgo,MS}; Figure~\ref {f.saw}
is due to Tom Kennedy). The 
conjecture that the half-plane SAW scaling limit is chordal $SLE_{8/3}$
has recently been comforted by simulations \cite {Kenn}.

Let us briefly conclude this chapter by mentioning the 
following characterization of $SLE_{8/3}$ that does not use
explicitly  Loewner's equation (even though its proof does):

\begin {theorem}
\label {t.uniq}
Chordal $SLE_{8/3}$ is the unique measure on continuous
simple  curves 
$\gamma$ from $0$ to $\infty$ in $\H$ such that for all Hull,
the law of $\gamma $ conditioned to avoid $A$ is identical to 
the law of $\Psi^{-1} (\gamma)$.
\end {theorem}

The proof of this Theorem uses the complete description of all 
measures on simply connected closed sets
(not necessarily curves)  joining $0$ to $\infty$
in $\H$ that satisfy this condition. These measures (called restriction 
measures in \cite {LSWrest})
are  constructed using $SLE_\kappa$ (in fact, by adding Brownian
bubbles to the $SLE_\kappa$ paths) for other values
of $\kappa$ (in fact for $\kappa \in (0, 8/3]$) and it
turns out that the only measure with these
properties that is supported on simple curves  is $SLE_{8/3}$.

\section*{Bibliographical comments}
All the material of this chapter is borrowed from \cite {LSWrest},
to which we refer for further details. 
The locality  property for $SLE_6$ was first proved in  
 \cite {LSW2}, using a different method. 
Restriction properties are closely related to 
conformal field theory \cite {BPZ0,BPZ,Ca1,Ca2,Cabook,CaChuo},
as pointed out in \cite {FW1,FW2}.
They have also interpretations in terms of highest-weight representations of the 
Lie algebra of polynomial vector fields on the unit circle. In fact, Theorem~\ref {t.uniq} 
corresponds to the fact that the unique such representation that is degenerate at level 2 has 
its highest weight equal to $5/8$. See \cite {FW1,FW2}.

\chapter {SLE and the Brownian frontier}

\section {A reflected Brownian motion}

In this section, we introduce a  two-dimensional Brownian motion
with  a certain oblique reflection on the boundary
of a domain, and
we will  
relate its outer boundary to that of  $SLE_6$. 

Let us first define this 
reflected Brownian motion in the upper half-plane
$\H$.
Define for any $x \in \R$, 
the vector $u(x) = \exp (i \pi /3)$ if $x \ge0$
and $u(x) = \exp (2 i \pi /3)$ if $x <0$.
It is the vector field 
with angle $2 \pi /3$ pointing ``away from the origin''.
Suppose that $Z^*_t =  X_t^* + i Y_t^* $ is an ordinary planar Brownian path 
started from 
$0$. Then, there exists a unique pair $(Z_t, \ell_t)$ 
of continuous processes such that $Z_t$ takes its values
in $\overline \H$, $\ell_t$ is a non-decreasing real-valued
continuous function with $\ell_0=0$
 that increases only when $Z_t \in \R$, and
$$
Z_t =  Z_t^*  + \int_0^t u( Z_s) d\ell_s.
$$
The process $(Z_t, t \ge 0)$ is called the reflected 
Brownian motion in $\H$ with reflection vector field $u(\cdot)$.
Note that the process $Z$ in fact only depends on the direction
of $u(\cdot)$ and not on its modulus. For instance $Z$ is also 
the reflected Brownian motion in $\H$ with reflection vector field
$2u (\cdot)$ (just change $\ell$ into $\ell/2$).

An equivalent way to define this process is to first
define the reflected
(one-dimensional) Brownian motion  
$$ Y_t =  Y_t^*  - \min_{s \in [0,t]}  Y_s^* 
.$$
The local time at $0$ of $Y$ is simply 
$l_t = -\min_{[0,t]}  Y^*$.
Then, define $X$ in such a way that 
$$
X_t =  X_t^* + \int_0^t {\hbox {sgn}}(X_s) \frac {1}{\sqrt {3}} dl_s 
$$
and verify that  $Z_t = X_t + i Y_t$ satisfy the 
required conditions.

Brownian motion with oblique reflection on domains have been
extensively studied, and this is not the proper place to 
review all results. We just mention that the general theory
of such processes (e.g., \cite {VW})
 ensures that the previously defined
process $Z^*$  
exists.

Reflected planar Brownian motion 
(even with oblique reflection)
are also 
invariant under conformal transformations.
Suppose for instance that $\phi$ is a conformal transformation
from a smooth subset $V$ (such that $[-1,1] \subset \partial V$)
 of $\H$ 
onto a smooth domain $D$.
Recall that 
$$Z_t =  Z_t^* + \int_0^t  u(Z_s) d\ell_s.
$$
Define 
$$
\sigma_V := \inf \{ t > 0 \ : \ \partial V \setminus (-1, 1) \}.
$$
Taylor-expanding 
each term in the sum 
\begin {eqnarray*}
{\phi ( Z_t) - \phi (0)}
= \sum_{j=1}^{n}
(\phi (Z_{jt/n}) - \phi (Z_{(j-1)t/n}))
\end {eqnarray*}
just as in the proof of It\^o's formula (letting $n \to \infty$),
it follows (using the fact that the real and imaginary parts of
$\phi$ are harmonic) that for all $t \le \sigma_V$, 
$$
\phi (Z_t) = \int_0^t \phi'(Z_s) d Z_s^* 
+  
+  \int_0^t u(Z_s) \phi' (Z_s) d\ell_s  
. $$
Hence, if one time-changes $\phi(Z)$ using the clock 
$u(t) = \int_0^t |\phi' (Z_s)|^2 ds$, we see
that 
$\phi (Z_u)$ is also a stopped  reflected Brownian motion in $D$ 
with the reflection 
vector field
$(\phi'(\phi^{-1} (\cdot)) \times u (\phi^{-1} (\cdot))$
on $\partial D$.

\medbreak

This has the following useful consequences:
Suppose that $V \subset \H$ and $\sigma_V$ are as before.
Note that 
$\sigma_\H$ is the first time at which $Z_t$ hits 
$\R \setminus (-1,1)$.
There exists a unique conformal map $\phi$ from $V$ onto $\H$
such that $\phi (-1)=-1$, $\phi(0)=0$ and $\phi (1) =1$.
\begin {lemma}
\label {l.cirbm}
Modulo time-change, the laws of $(\phi (Z_t), t \le \sigma_V)$
and of $(Z_t , t \le \sigma_{\H})$ are identical.
\end {lemma}
In other words, The reflected Brownian motion $Z$ satisfies the 
same locality property as $SLE_6$.

\medbreak
A slight modification of the above proof of conformal 
invariance for reflected Brownian motions shows that the image of 
$Z$ under the conformal map 
$z \mapsto z^{1/3}$ from $\H$ onto the 
wedge 
$${\cal W} := \{ r e^{i \theta} \ : \ r >0 , \theta \in (0, \pi/3) \}$$
is reflected Brownian motion in that wedge, started from the 
origin, with reflection vector field
$u(x)= e^{i \pi / 3}$ on $\R_+$ and $u(x) = 1$ on $e^{i \pi/ 3} \R_+$.
We use this observation to give a simple proof of the following fact
on hitting probabilities for $Z$:
\begin {lemma}
\label {l.cardybm}
Suppose that $\Phi$ is the conformal transformation from $\H$ 
onto an equilateral triangle $OAC$ such that 
$\Phi (0) = O$, $\Phi(-1) = A$ and $\Phi (1) = C$.
Then, the law of $\Phi (Z_{\sigma_\H})$ 
is uniform on $AC$.
\end {lemma}
\proof
One elementary convincing proof uses discrete approximations. 
Here is a brief outline of this proof:
Define $\omega = \exp ( i \pi /3 )$.
Consider a triangular grid in the wedge ${\cal W}$ i.e.
$ \{ m + m' \omega  \ : \ m, m' \ge 0 \}$.
Let $(S_n, n \ge 0)$
denote simple random walk on this grid that is started from 
$0$. In the inside of ${\cal W}$, its transition 
probabilities are that of simple random walk (with probability $1/6$ to 
jump to each of its neighbours).
When $S$ hits the (positive) real line at $x$, it has the 
following transition probabilities:
$p(x , x+1) =  1/3$ and
$$
p (x, x -1) = p (x, x+\omega )= p (x, 
x + \omega^2 ) = p (x,x) = \frac 16 
.$$
and the symmetric ones on $\omega \N$:
$p(x,x+ \omega) =1/3$ and
$$
p(x, x +1) = p (x, x +1/ \omega) = p (x, x+  1/ \omega^2)
= p(x,x)= \frac 16.
$$
Finally, at the origin, $p(0,1)=p(0,\omega)=1/2$.
It is not difficult to see that in the scaling limit, such a 
random walk converges to reflected Brownian motion in ${\cal W}$ 
with the reflection vector field $u(\cdot)$ on $\partial {\cal W}$.
This is due to the fact that the bias of the simple random walk
when it hits $\partial {\cal W}$ is proportional to $u$.
Moreover, it is easy to check that if $S_0 = 0$, then
if one writes $S_n = e^{i\pi /6} r_n + \omega^2 s_n$, 
then the conditional law of $s_n$ given
  $(r_j, j \le n)$ is the uniform distribution among the permitted
values of $s$ given $r_n$. In other words, the ``uniform distribution 
of $s$ is preserved, independently from $r$''. In particular,
the hitting distribution of the simple random walk
$S$ on the segment $N + \omega^2 [0,N]$,  
is simply the uniform distribution on $\{ N, N+ \omega^2 ,
N + 2 \omega^2 , \ldots, N + N \omega^2 \}$.
The Lemma follows, letting $N \to \infty$. 
\qed
\medbreak

\begin{figure}
\centerline{\includegraphics*[height=2.3in]{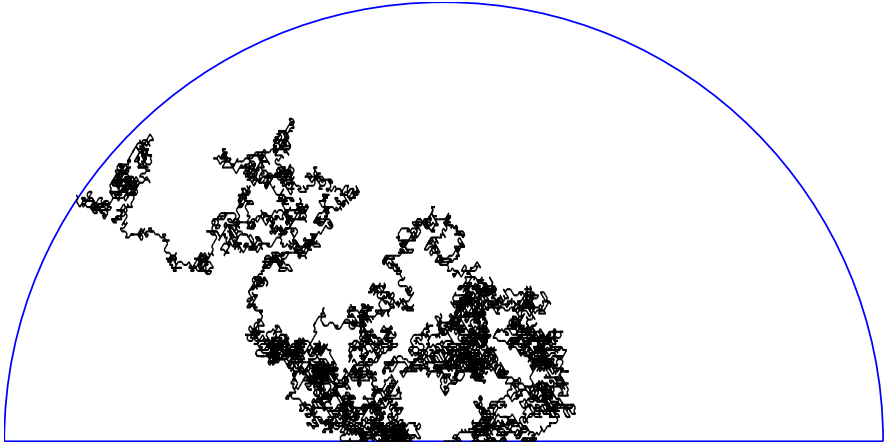}}
\caption{\label{f.rbm}
The reflected Brownian motion stopped at its 
hitting time of the unit circle}
\end{figure}

We are now ready to state and prove the following result:
\begin {theorem}
Define the following two sets:
\begin {itemize}
\item
Consider chordal $SLE_6$ $(K_t, t \ge 0)$
in $\H$ (or in $V$) up to its first hitting
time $T$ of $\R \setminus (-1,1)$.
Let $e$ denote the point at which the SLE hits 
$\R \setminus (-1,1)$, and let $E := 
\{e \} \cup \cup_{t<T} K_t$.

\item 
Consider the set of points $F$ in $\overline \H$ that 
are disconnected (in $
\H$) from $\R \setminus (-1,1)$
by $Z [ 0, \sigma_{\H}]$.
\end {itemize}
Then, the laws of $E$ and of $F$ are identical.
\end {theorem}

\proof
Note that Lemma~\ref {l.cardybm}, Lemma~\ref {l.cirbm}, 
Theorem~\ref {t.loc} and Proposition \ref {p.unif}
show that $E$ and $F$ both have the following properties:
\begin {itemize}
\item 
They are random compact sets that intersect 
$\R \setminus (-1,1)$ at just one point $x$ and the law
of $\Phi (x)$ is uniform on $AC$.
\item
Their complement in $\overline \H$ consists of two 
connected components (one unbounded, one bounded).
\item
For all $V$ as before, the 
probability that $E \subset V$ is   
identical to the probability that $\sigma_V = \sigma_\H$
(and the corresponding result for $F$).
\end {itemize}
If we combine these two properties, we see that for all
such 
$V$,
$$
\P [ E \subset V ] = \P [F \subset V ] 
= \frac {\hbox {length} (\Phi \circ  \phi ( \partial V \setminus \R))}{AC}
$$
(this is because the law of the image under $\Phi \circ \phi$ of the 
``hitting point'' of $\partial V \setminus (-1,1)$ is uniform
on $AC$.
But this determines completely the laws of $E$ and of $F$
and therefore implies that they are equal.
\qed
\medbreak

\begin{figure}
\centerline{\includegraphics*[height=3in]{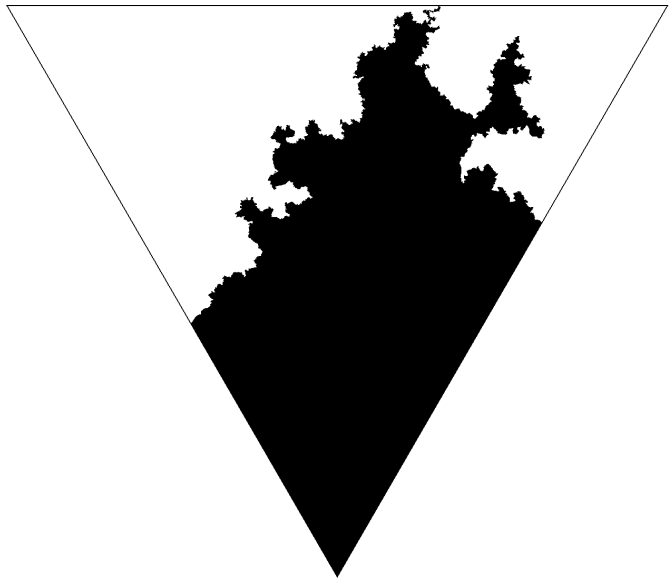}}
\caption{\label{f.tri}
The 
filling of RBM (or of the $SLE_6$ curve)
in a triangle}
\end{figure}

Using
 conformal invariance, the previous result can be 
adapted in any domain. For instance, Figure \ref {f.tri}
could represent both the filling of a reflected Brownian motion
(or of a $SLE_6$ curve), 
started at the bottom of 
the triangle stopped at their first hitting of the top
segment. Recall that the law of this hitting point is uniformly
distributed.

\section {Brownian excursions and $SLE_{8/3}$}

We now describe a probability measure on 
Brownian excursions from $0$ to infinity in $\H$ (which is 
closely related to the measures on excursions 
that were considered in \cite {LW2}).
One can view this measure on 
paths as the law of planar Brownian motion $W$
(not to be confused with the $\sqrt {\kappa} B$
in the previous chapters) started 
from $0$ and conditioned to stay in $\H$ at all 
positive times.

Let $X$ and  $Y$ denote two independent processes such that 
$X$ is standard one-dimensional Brownian motion and $Y$ is 
a three-dimensional Bessel process (see e.g., \cite {RY}
for background  on three-dimensional Bessel processes, its
relation to Brownian motion conditioned to stay positive
and stochastic differential equations) that are both started 
from $0$.
Let us briefly recall that a three-dimensional Bessel process is 
the modulus of a three-dimensional Brownian motion, and that it can be 
defined as the solution to the stochastic differential equation
$$ 
dY_t = dw_t + \frac {1}{Y_t} dt
$$
(where $w$ is one-dimensional 
standard Brownian motion). It is very easy to see that
$(1/Y, t \ge t_0)$ is a local martingale for all $t_0>0$, and 
that if $T_r$ denotes the hitting time of $r$ by $Y$, 
then the law of $(Y_{T_r + t} , t < T_R - T_r )$ is 
identical to that of a Brownian motion started from $r$ and
conditioned to hit $R$ before $0$ (if $0<r<R$).
Loosely speaking $Y$ is a Brownian motion started from $0$
and conditioned to stay forever positive. Note that almost
surely $\lim_{t \to \infty} Y_t = \infty$.
 
We now define $W = X+iY$.
In other words,
 $W$ has the same law as the solution to the 
following stochastic differential equation: 
\begin {equation}
\label {e.sde}
dW_t = d\beta_t + i \frac {1}{ \Im (W_t)} dt
\end {equation}
with $W_0 = 0$,
where $\beta$ is a complex-valued Brownian motion.
Note that $W$ is a strong Markov process.
Let $T_r$ denote the hitting time of the line $
\R + ir$
by this process $W$ (i.e., the hitting time of $r$ by $X$).
Let $S$ denote a random variable with the same 
law as $W_{T_1}$.
Then,
 scaling and
the relation between one-dimensional Brownian motion 
conditioned to stay positive and the three-dimensional Bessel 
process shows immediately that for all 
$0<r<R$, the law of $W [ T_r, T_R]$ is the law of a Brownian motion
started with the same law as $rS$, stopped at its first hitting
of $iR + \R$, and conditioned to stay in the upper 
half-plane up to that time.
Note that the probability of this event is $r/R$.

\begin{figure}
\centerline{\includegraphics*[height=3in]{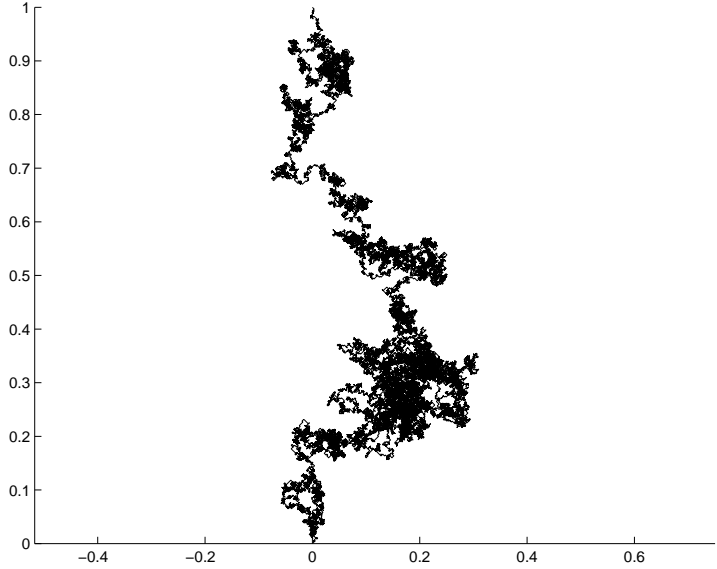}}
\caption{\label{f.ex}
An excursion from $0$ to $i$
in the strip $\R \times [0,1]$} 
\end{figure}

By mapping conformally $\H$ onto any other simply
connected domain $D$ ($D \not= \C$), and looking at the 
image of the Brownian excursion in $\H$ under this map, one
gets the law of a Brownian excursion in $D$ from the image of 
$0$ to the image of $\infty$.
As for SLE, this law is well-defined up to linear 
time-change. One can also directly define this excursion in $D$
as the solution to a stochastic differential equation ``forcing the 
Brownian motion to hit $\partial D$ at the image of infinity.''

The following result was observed by B\'alint Vir\'ag \cite {V} (see also 
\cite {LW2,LSWrest}):

\begin{lemma} Suppose $A$ is a Hull
and $W$ is a Brownian excursion in $\H$
from $0$ to $\infty$.  Then
$\P [W [0, \infty) \cap  A = \emptyset ]  
  = \Phi_A'(0)$.
\end{lemma}

\proof
Suppose   that 
$W$ is a  solution to (\ref {e.sde})
started from $z \in \Phi^{-1} (\H)$.
Let $Z$ denote a planar Brownian motion 
started from $z$.
Let $\tau_R(V)$  denote the hitting time of $iR + \R$ by a process 
$V$.
When $\Im (z) \to \infty$, $\Im (\Phi (z)) = \Im (z) + o(1)$,
and it therefore follows easily  
from the strong Markov property  
of planar Brownian motion that when $R \to \infty$,
$$
\P [ \Phi (Z) [ 0, \tau_R (Z)]  \subset \H ]
\sim
\P [ \Phi (Z) [ 0, \tau_R (\Phi (Z))] \subset \H]
.$$
But since $\Phi (Z)$ is a time-changed Brownian motion,
the right-hand probability is equal to $
\Im ( \Phi(z)) / R$, so that 
$$\P [ W [ 0, \tau_R (W) ] \subset \Phi^{-1} (\H) ] 
=
\frac { \P [ Z [ 0, \tau_R (Z)] \subset \Phi^{-1} (\H ) ] }
{ \P [ Z [ 0, \tau_R (Z) ] \subset \H ] } 
=
\frac { \Im \Phi (z)}{\Im (z) } + o(1) 
$$
when $R \to \infty$. In the limit $R \to \infty$, we get 
\begin {equation}
\label {e.Im}
\P [ W[0, \infty)  \subset \Phi^{-1} (\H) \mid W_0 = z ] 
= \frac { \Im \Phi(z)}{\Im (z)}
.\end {equation}
When $z \to 0$, 
$ \Phi (z) = z  \Phi' (0) + O ( |z|^2) $  so that 
\begin {eqnarray*}
\P [ 
W[0, \infty) \subset \Phi^{-1} ( \H ) ] 
&=& \lim_{r \to 0} \P [ W [ T_r, \infty ) \subset \Phi^{-1} (\H) ] 
\\&=& \lim_{r \to 0} \expect [ \Im (\Phi (rA)) / \Im (rA) ] 
\\&=& \Phi'(0) 
\end {eqnarray*}
(one can use dominated convergence here since $\Im ( \Phi(z)) \le 
\Im (z)$ for all $z$).
\qed
\medbreak

We now define the filling ${\cal H}$ of $W[0, \infty)$ as the set of points in 
$\overline {\H}$ that are disconnected from $\R$ by 
$W [ 0, \infty)$. This set is obtained by filling in 
all the bounded connected components of the complement of 
the curve $W$.
Then, ${\cal H}$ is a closed unbounded set and $\H \setminus 
{\cal H}$ consists of two open connected components 
(with $[0,\infty)$ and $(- \infty, 0]$ on their respective boundaries).
The law of such a random set is characterized by the values of 
$\P [ {\cal H} \cap A = \emptyset ]$, where $A$ spans all 
Hulls, because this family of
events turn out to generate the $\sigma$-field on which ${\cal H}$
is defined, and to be stable under finite intersections. 
Hence, as in the case of $K_\infty$ for $SLE_{8/3}$, 
the fact that 
\begin {equation}
\label {e.ber}
\P [ {\cal H} \cap A = \emptyset ] = \Phi'(0)
\end {equation}
characterizes the law of ${\cal H}$ and yields that ${\cal H}$
also satisfies Corollary~\ref {c.rest}.

\begin {theorem}
\label {5=8}
Suppose that ${\cal H}_8$ denotes the filling of the union of 
$8$ independent chordal $SLE_{8/3}$'s. Suppose that ${\cal H}_5$ denotes 
the filling of the union of $5$ independent Brownian excursions.
Then, ${\cal H}_5$ and ${\cal H}_8$ have the same law.
\end {theorem}

\proof
This is simply due to the fact that for all Hull $A$ 
$$
\P [ {\cal H}_5 \cap A = \emptyset 
] = \P [ {\cal H}_8 \cap A = \emptyset 
] = \Phi_A '(0)^5
$$
and that this characterizes these laws.
\qed
\medbreak

This has various nice consequences (see \cite {LSWrest}),
 some of which  we now
heuristically 
describe:  
First, since the boundary of ${\cal H}_8$ consists 
of the union of some parts of the $SLE_{8/3}$
curves, it follows that ``locally'', the outer boundary of 
a Brownian excursion (and therefore also of a Brownian 
motion) looks like one $SLE_{8/3}$ path. 
In the previous section, we did see that the outer
boundaries of reflected Brownian motion  
and of $SLE_6$ are the same. 
Hence, ``locally'', the outer frontiers of $SLE_6$ and of 
planar Brownian motion look like an $SLE_{8/3}$ curve.
Furthermore, since $SLE_{8/3}$ is symmetric, this shows that 
one cannot distinguish the inside from the outside of a 
planar Brownian curve by only seeing a part of its frontier. 
Since $SLE_{8/3}$ is conjectured to be the scaling limit 
of self-avoiding walks, this would also show that the 
Brownian frontier looks locally like the scaling limit
of long self-avoiding curves (see \cite {LSWSAW}).

\bibc
The idea that conformal invariance and restriction defines 
measures on random sets and makes
it possible to understand the Brownian frontier
in terms of other models (or the corresponding exponents) first 
appears in \cite {LW2}.
Most of the material of this chapter is borrowed from
\cite {LSWrest}.

A discussion of the conjectured relation between
$SLE_{8/3}$ and planar self-avoiding walks is discussed in \cite {LSWSAW};
one can in particular recover the predictions of Nienhuis \cite {N2}
on the critical exponents for self-avoiding walks using SLE arguments.

The fact that the Brownian frontier had the same dimension as the
scaling limit of self-avoiding walks was first observed visually by 
Mandelbrot \cite {Ma}.
  
\chapter {Radial SLE}

\section {Definitions} 

Motivated by the example of LERW (among others)
given in 
the introductory chapter, 
we now want to find a nice way to encode 
growing families of compact subsets $(K_t, t \ge 0)$
of the closed unit disk that are growing from the 
boundary point $1$ towards $0$. 
As in the chordal case, we are in fact going to 
focus on the conformal geometry of the complement $H_t$ of 
$K_t$ in the unit disc $\U$. One first has to 
find a natural time-parametrization.
It turns out to be convenient to 
define the conformal map $g_t$ from
$H_t$ onto $\U$ that is normalised
by
$$ 
g_t (0) = 0 \hbox { and } g_t'(0) > 0 .
$$
Note that $g_t'(0) \ge 1$.
This can be for instance derived using the fact that
$\log g_t' (0)$ is the limit when $\eps \to 0$ 
of $\log (1/ \eps)$ times
the probability that a planar Brownian 
motion started from $\eps$ hits the circle of
radius $\eps^2$ before exiting $H_t$ (an analyst
would find this justification very strange, for sure).

Then (and this is simply because 
with obvious notation, $(\tilde g_s \circ g_t ) (0)
= \tilde g_s (0) \circ g_t'(0) $),  
one measures the ``size'' $a(K_t)$
of $K_t$ via the 
derivative of $g_t$ at the origin:
$$
g_t'(0) = \exp (a(t)).
$$
Hence, we will consider growing families of 
compact sets such that $a(K_t) = t$.

Suppose now that $(\zeta_t, t \ge 0)$ 
is a continuous 
function on the unit circle $\partial \U$. 
Define for all $z \in \overline \U$,
the solution $g_t(z)$ to the ODE
\begin {equation}
\label {rODE}
\partial_t g_t (z) = - g_t (z) \frac { g_t (z) + \zeta_t}{g_t (z) - 
\zeta_t}
\end {equation}
such that 
$g_0 (z) = z$.
This solution is well-defined up to the 
(possibly infinite) time $T(z)$ defined by
$$
T(z) = \sup \{ t> 0 \ : \ 
\min_{s \in [0,t)} |g_s (z)- \zeta_s| 
>0 \}.
$$
We then define
$$
K_t := \{ z \in \overline \U \ : \ T(z) \le t \}
$$
and 
$$ U_t := \U \setminus K_t .$$
The family $(K_t, t\ge 0)$ is called the (radial) Loewner chain
associated to the driving function $\zeta$.

The general statements that we described in the chordal case 
are also valid in this radial case.
One can add one feature that has no analog in the 
chordal case:
It is possible to estimate the Euclidean distance $d_t$ from $0$
to $K_t$ in terms of $a(t)=t$.
Indeed, since $U_t$ contains the disc $d_t \times  \U$, it is clear
that $g_t' (0) \le 1/ d_t$.
On the other hand, a classical result of the theory of 
conformal mappings known as Koebe's $1/4$ Theorem states that 
(if $a(K_t) = t$)  
$1/d_t \le 4 g_t' (0)$.
This is loosely speaking due to the fact that the best $K_t$ 
can do to get as close to $0$ 
in ``time $t$'' is to shoot straight i.e. to choose $\zeta =1 $.
Hence, for all $t \ge 0$,
\begin {equation}
\label {e.koebe}
 e^{-t}/4 \le d(0, K_t) \le e^{-t}.
\end {equation}
This will be quite useful later on.

Radial $SLE_\kappa$ is then simply the random family of 
sets $(K_t, t \ge 0)$ that is obtained when 
$$
\zeta_t = \exp ( i \sqrt {\kappa} B_t )
$$
where $\kappa>0$ is fixed and $(B_t, t \ge 0)$ is 
standard one-dimensional Brownian motion.

As in the chordal case, one can then define radial $SLE$
from $a \in \partial D$
to $b \in D$ in any open simply connected domain $D$  
by taking the image of radial $SLE$ in $\U$ under
the conformal map $\Phi$ from $\U$ onto $D$ such 
that $\Phi(1)=a$ and $\Phi (0)=b$.
Note that this time, the
time-parametrization is also 
well-defined since there exists 
only one such conformal map (recall that in the
chordal case, one had to invoke the
scaling property to make
sure that chordal SLE in other domains than the 
half-space was properly defined).

\section {Relation between radial and chordal SLE}

In this section, we show that 
 chordal SLE and radial SLE 
are very closely related. 
Let us start with the special case $\kappa =6$.

\begin{theorem}
\label{t.equiv}
Suppose that $x \in (0, 2 \pi)$.
Let 
$(K_t, t \ge 0)$ be a radial $SLE_6$ process.
Set 
$$
T := \inf \{t\ge 0 \  : \  \exp (ix) \in{K_t}\}.
$$
Let $(\tilde K_u , u \ge 0)$ be a chordal $SLE_6$
 process in $\U$ starting
also at $1$ and growing towards $\exp (ix)$, and let 
$$
\tilde T := \inf \{u\ge 0 \ : \  0\in\tilde K_u\}.
$$
Then, up to a random time change, the process $t\mapsto K_t$
restricted to $[0,T)$ has the same law as the process
$u\mapsto\tilde K_u$ restricted to $[0,\tilde T)$.
\end{theorem}
Note that $T$ (resp.\ $\tilde T$) is the first time where 
$K_t$ (resp.\ $\tilde K_u$) disconnects 0 from 1.

When $\kappa \not= 6$, a weaker form of equivalence holds:

\begin {prop}
\label {p.equiv}
Let $(K_t, t \ge 0)$, $(\tilde K_u, u \ge 0)$, $T$ and $\tilde T$
be defined just as in Theorem \ref {t.equiv}, except that they
are SLE with general $\kappa >0$.
There exist two nondecreasing families 
of stopping times $(T_n , n \ge 1)$ and $(\tilde T_n , n \ge 1)$
such that almost surely,
$T_n \to T$ and $\tilde T_n \to \tilde T$ when $n \to \infty$,
and such that for each $n \ge 1$, the laws of 
$ (K_t , t \in [0, T_n])$ 
and $(\tilde K_u , u \in [0, \tilde T_n])$
are equivalent (in the sense that they have a positive density
with respect to each other) modulo increasing time change.
\end {prop}

These results imply that the properties of chordal SLE 
such as ``being generated by a continuous curve'' are also 
valid for radial SLE. 

We prove both results simultaneously:

\proof
Let us first briefly recall 
how $\tilde K_u$ is defined.
For convenience, we will restrict ourselves
to $x=\pi$ (the proof
in the general case is almost identical). 
Define the conformal map
$$ \psi (z) = i \frac {1-z}{1+z}$$
from $\U$ onto $\H$ that satisfies
 $\psi(-1)=\infty$, $\psi(1)=0$, and  $\psi(0)=i$.
Suppose that $u \mapsto \tilde B_u$ is a real-valued 
Brownian motion such that $\tilde B_0 =0$.
For all $z \in \U$, define the function $\tilde g_u = \tilde g_u (z)$
such that $\tilde g_0 (z) = \psi (z)$ and 
$$
\partial_u \tilde g_u = \frac {2} { \tilde g_u - \sqrt {\kappa} \tilde
B_u}.
$$
This function is defined up to the (possibly infinite)
 time $\tilde T_z$ where 
$\tilde g_u(z)$ hits $\sqrt {\kappa} \tilde B_u$.
Then, $\tilde K_u$ is defined by $\tilde K_u = \{ z \in \U \ : \ \tilde T_z 
\le u\}$, so that $\tilde g_u$ is a conformal map from 
$\U \setminus \tilde K_u$ onto the upper half-plane.
This defines
the process
$( \tilde K_u , u \ge 0)$.
 
We are now going to compare it to  radial SLE.
Let $g_t:\U\setminus K_t\to\U$ be the conformal map 
normalized by $g_t(0)=0$ and $g'_t(0)>0$.
Recall that
\begin{equation}
\label{e.loewd}
\partial_t g_t (z) =  g_t (z) \frac{\zeta_t +g_t (z)}{\zeta_t -g_t (z)}\,, 
\end{equation}
where $\zeta_t= \exp (i \sqrt {\kappa}  B_t)$,
and $B$ is Brownian motion
on $\R$   with $B_0=0$.
Let $\psi$ be the same conformal map as before,
and define
 \begin {eqnarray*}
e_t &:=& g_t (-1),
\\
f_t(z) & := & \psi\bigl(g_t(z)/ e_t \bigr),
\\
\gamm_t&  := &
 \psi\bigl(\zeta_t / e_t \bigr).
\end {eqnarray*}
These are well defined, as long as $t<T$.
Note that 
$f_t$ is a conformal map from $\U \setminus K_t$
onto the upper half-plane, $f_t (1) = \infty$,
and $\gamm_t \in \R$.
From (\ref{e.loewd}) it follows that
$$
\partial_t f =-\frac {(1+\gamm^2)(1+f^2)}{ 2(\gamm-f)}.
$$
Let 
$$\phi_t(z)=a(t)z+b(t)$$
 where
$$
a(0)=1,\qquad \partial_ta=-(1+\gamm^2)a/2
$$ 
and
$$
\qquad b(0)=0,\qquad
\partial_tb=-(1+\gamm^2)a\gamm /2.
$$
Set 
\begin {eqnarray*}
h_t
&:= &  \phi_t\circ f_t\,,\\
\beta _t
& := &
\phi_t\bigl( \gamm (t)\bigr)
.\end {eqnarray*}
Then
(and this is the reason for the choice of the 
functions $a$ and $b$)
$$
\partial_t  h =
-(a/2)
\frac {(1+\gamm^2)^2}{ \gamm- f}
=
-\frac {(1+\gamm^2)^2a^2/2}{  \beta - h}\,.
$$
$h_t$ is also a conformal map from $\U \setminus K_t$
onto the upper half-plane with $h_t (1) = \infty$. Note also
that $h_0 (z) = \psi (z)$.
We introduce a new time parameter $u = u(t)$ by
setting
$$
\partial_t u = {(1+\gamm^2)^2a^2/4},\qquad u(0)=0\,.
$$
Then
$$
\frac
{\partial h}{\partial u} = \frac {-2}
{   \beta  - 
h}.$$
Since this is the equation defining the chordal SLE process,
it remains to show that
$ u  \mapsto   \beta _{t(u)} / \sqrt {\kappa}$ 
is related to 
Brownian motion (stopped at some random time).
This is a direct but  tedious
application of  It\^o's formula:
$$
d\gamm_t = \frac {(1+\gamm^2) \sqrt {\kappa}}{2}\, dB_t
+ \frac {\gamm (1+\gamm^2)}{2} 
\left( \frac {\kappa}{2} - 1  \right)dt 
$$
and 
$$
d  \beta_t 
=
\frac {  (1+ \gamm^2) a}{2}
\left( \sqrt {\kappa}\, dB_t + ( -3 + \frac {\kappa}{2}) \gamm \,dt 
\right).
$$
When $\kappa=6$, the drift term disappears and this proves
Theorem~\ref{t.equiv}.
When $\kappa \not=6$, the drift term does not
disappear. However,
the law of $u\mapsto  \beta_{t(u)}$
 is absolutely continuous with
respect to that of $\sqrt {\kappa}$ times a Brownian motion,
as long as $\gamm$ and $u$ remain bounded.
More precisely:
It suffices to take 
$$ T_n = \min \Bigl\{ n,\, \inf \{ t > 0 \ : \ |\zeta_t - e_t | < 1/ n \}
\Bigr\}. 
$$
Before $T_n$, $|\gamm|$ remains bounded,
$a$ is bounded away from $0$ (note also that $a \le 1$ always), so that 
$t/u$ is bounded and bounded away from $0$.
Hence, $u( T_n)$ is also bounded (since $T_n \le n$).

It now follows directly from Girsanov's Theorem (see e.g., \cite {RY})
that the law of $\bigl(\beta(u) / \sqrt {\kappa}\bigr)_{ u \le  u(T_n)}$
is equivalent to that of 
 Brownian motion up to some (bounded)
stopping time, and Proposition~\ref {p.equiv} follows.
\qed

\section {Radial $SLE_6$ and reflected Brownian motion}

If one combines the radial-chordal equivalence for $SLE_6$
with the locality property for chordal $SLE_6$, 
one gets immediately a locality property for
radial $SLE_6$, and the relation between 
fillings of radial $SLE_6$ and 
of reflected Brownian motion.
We do not state the locality property here (and leave it to 
the interested reader), but we state the relation 
between fillings of radial $SLE$ and 
of reflected Brownian motions that we will use 
in the next chapters.

Before that, we have to say some words about how this
reflected Brownian motion is defined in the unit disc.
Suppose that $(Z_t, t \ge 0)$ is the reflected Brownian
motion in the upper half-plane with reflection 
angle $2 \pi /3$ away from the origin as in the previous 
chapter.
Let us now define
$$ \tilde Z_t 
:= \exp (- i Z_t )$$
so that $\tilde Z$ takes its values in the unit disk
and is started from $\tilde Z_0=1$.
Clearly, since $\tilde Z_t \not= 0 $ for all $t$, one can define the
continuous version of its argument $(\theta_t, t \ge 0)$.
Conformal invariance of planar Brownian motion shows that 
$\tilde Z_t$ behaves like (time-changed) Brownian motion
as long as it stays away from the unit circle, and when 
it hits the unit circle, then it is reflected 
with angle $2\pi /3$ in the direction that 
``increases'' $|\theta|$.
Define 
$$
\tilde \sigma (r) 
: = \inf \{ t >0 \ : \ |\tilde Z_t | = r\}
$$
which is also the first time at which the imaginary part
of $Z$ hits $\log (1/r)$. 
 
\begin {theorem}
\label {t.hulleq}
Suppose that $r<1$. Define the two following random hulls:
\begin {itemize}
\item
Suppose that $(K_t, t \ge 0)$ 
is radial $SLE_6$ as before. 
Let $\tau_r$ denote the first time
at which radial $K_t$ intersects the 
circle $\{ |z| = r \}$.
Define the event ${\cal H} (x , \tau_r)$ that
 $K_{\tau_r}$ does
not disconnect $0$ from $\exp (ix)$.

\item
On the event $\tilde {\cal H} (x, \tilde \sigma_r)$ that 
$\tilde Z [0, \tilde \sigma_r]$ does not disconnect $0$ from
$\exp (ix)$, define the connected component 
$H$ of $\U \setminus \tilde Z$ that contains $0$, 
and the hull $\tilde K_{\sigma_r} = \overline {\U \setminus H}$.
\end {itemize}
Then, the two random sets
$1_{{\cal H} (x, \tau_r)} K_{\tau_r } $ and
$1_{\tilde {\cal H}_r ( x ,\tilde \sigma_r )} 
\tilde K_{\tilde \sigma_r}$
have the same law.
\end {theorem}

In particular, 
$$
\P [ {\cal H}(x, \tau_r) ] = \P [ \tilde {\cal H} (x , \tilde \sigma_r))].
$$
This shows that one can compute non-disconnection
probabilities for reflecting 
Brownian motions using radial $SLE_6$.

\section* {Bibliographical comments}
For basic results on Loewner's equation, 
and basic complex analysis, we refer again to 
\cite {A1, A2, Dur, Hay}.
The radial-chordal equivalence for $SLE_6$ has been derived in \cite {LSW2}.

\chapter {Some critical exponents for SLE}

\section {Disconnection exponents}

In this section, we fix $\kappa >4$,
and we consider radial $SLE_\kappa$
in the unit disc started from $1$.
Our goal will be to estimate 
probabilities
of events
like 
$$
{\cal H} (x,t) = \{ \exp (ix) \in \partial H_t \}
$$
that $K_t$ has not swallowed the point $\exp (ix) \in
\partial \U$ from $0$ at time $t$. 
Let us define
 the numbers
$$
q_0 = q_0(\kappa) 
:=  1  - \frac 4 \kappa
$$
and
$$
\lambda_0 = \lambda_0 (\kappa) 
:= \frac {\kappa }{8} - \frac 12.
$$

\begin {prop}
\label {p.disco}
There exists a constant $c$ such that for all $t \ge 1$
and for all $x \in (0, 2 \pi)$,
$$
e^{-\lambda_0 t} (\sin (x/2))^{q_0}
\le \P [ {\cal H}(x,t) ] \le c e^{-\lambda_0 t} (\sin (x/2))^{q_0}.
$$
\end {prop}

\proof
We will use the notation
$$
f(x,t) = \P [ {\cal H} (x,t) ] 
.$$
Let
$\zeta_t = \exp ( i \sqrt {\kappa} B_t )$
be the driving process of the radial $SLE_\kappa$, with $B_0=0$. 
For all $x \in (0, 2 \pi)$,
let $Y^x_t$ be the continuous real-valued function
of $t$ which satisfies 
$$ g_t (e^{ix} ) = \zeta_t \exp ( i Y_t^x ) $$ 
and $Y_0^x=x$.
The function $Y_t^x$ is defined on the set of pairs
$(x,t)$  such that ${\cal H}(x,t)$ holds.
Since
$g_t$ satisfies Loewner's differential equation
\begin{equation}\label{e.loew}
\partial_t g_t(z) = g_t(z) \frac{\zeta_t+g_t(z)}{\zeta_t-g_t(z)}\,,
\end{equation}
we find that
\begin{equation}\label{Ysde}
d Y_t^x = \sqrt {\kappa} \  dB_t 
+ \cot(Y_t^x/2)\ dt.
\end{equation}
Let 
$$
\tau^x 
:=
\inf\Bigl\{t \ge 0 \ : \  Y_t^x\in\{0,2\pi\}\Bigr\}\,
$$
denote the time at which $\exp (ix)$ is absorbed
by ${K_t}$, so that
$$
f (x,t) = \P [ \tau_x > t ].
$$
We therefore want to estimate the 
probability that the 
diffusion $Y^x$ (started from $x$) 
has not hit $\{0, 2 \pi \}$ before time $t$
as $t \to \infty$.
This is a standard problem.
The general theory of diffusion processes can be used to argue that 
	$f(x,t)$ is smooth on $(0,2\pi) \times \R_+$, and It\^o's 
	formula immediately shows that 
	\begin {equation}
	\label {e.od}
	\frac {\kappa}{2} \partial_{x}^2 f + \cot (x/2) 
	\partial_x f = \partial_t f.
	\end {equation}
	Moreover, for instance comparing $Y$ with Bessel processes
	when $Y$ is small, 
	one can easily see that (here we use that $\kappa >4$)
	for all $t >0$,
	\begin {equation}
	\label {bv}
	\lim_{x \to 0+} f(x,t) 
	= \lim_{x \to 2 \pi - } f(x,t) = 0.
	\end {equation}
	Hence, $f$ is solution to (\ref {e.od})
	with boundary values (\ref {bv}) and $f(x,0) = 1$.
	This in fact characterizes $f$, and its long-time
	behaviour is described in terms of the 
	first eigenvalue of the operator
	$\kappa \partial_x^2 / 2 + \cot (x/2) \partial_x$.
	More precisely, define  
	$$ 
	F(x,t) = \expect[   1_{{\cal H}(x,t)} \sin (Y_t^x/2)^{q_0} ]. $$
	Then, it is easy to see that $F$ also solves (\ref {e.od})
	with boundary values (\ref {bv}) but this time with initial 
	data $F(x,0) = \sin (x/2)^{q_0}$.
	One can for instance invoke the maximum principle 
	to construct a handcraft proof (as in \cite {LSW1})
	of the fact that this characterizes $F$.
	Since $e^{-\lambda_0 t} \sin (x/2)^{q_0}$
	also satisfies these conditions, it follows that 
	$$
	F(x,t) = e^{-\lambda_0 t} \sin (x/2)^{q_0}.
	$$
	Hence, 
	$$
	f(x,t) = \P [ {\cal H} (x,t)] 
	\ge \expect [ 1_{ {\cal H} (x,t)} \sin (Y_t^x /2)^{q_0}]
	=  e^{-\lambda_0 t} \sin (x/2)^{q_0}.
	$$
	To prove the other inequality, one can for instance use 
	an argument based on Harnack-type considerations:
	For instance, one can see 
	that (uniformly in $x$) a positive fraction of 
	the paths $(Y_t^x , t \in [0, 1])$ 
	such that $\tau_x > 1$ satisfy 
	$Y_1^x \in [\pi/2, 3 \pi /2]$.
	This then implies readily (using the Markov
	property at time $t-1$) that
	for all $t \ge 1$,
	$$
	f(x,t)  
	\le c_0 \P [ \tau_x > t 
	\hbox { and } Y_t^x \in [ \pi/2, 3\pi /2] ]  
	\le c_1 
	 F(x,t) =  
	c_1 e^{-\lambda_0 t} \sin (x/2)^{q_0}.$$
	\qed

	\section {Derivative exponents}

	The previous argument can be generalized in order to derive the 
	value of other exponents that will be very useful
	later on:
	 We
	will focus on the moments of the  derivative of $g_t$ at 
	$\exp (ix)$ on the event ${\cal H}(x,t)$.
	Note that on a heuristic level, $|g_t'(e^{ix})|$ measures
	how ``far'' $e^{ix}$ is from the origin in $H_t$.

	More precisely, we fix $b \ge 0$, and we 
	define
	$$
	f(x,t)
	:= 
	\expect \Bigl[\left|g_t'\bigl(\exp(i x)\bigr)\right|^b\,
	1_{{\cal H}(x,t)}\Bigr] .
	$$
	We also define the numbers 
	\begin {eqnarray*}
	\nonumber 
	q=q(\kappa,b)
	&
	:= &  \frac{\kappa-4+\sqrt{(\kappa -4)^2+16 b\kappa }}{2\kappa}
	\\
	\lambda  = \lambda (\kappa, b)
	&:=&
	\frac{8b+\kappa-4+\sqrt{(\kappa -4)^2+16 b \kappa }}{16}
	.\end {eqnarray*}
	The main result of this Section is the following
	generalization of Proposition~\ref {p.disco}:
	\begin {prop}
	\label {p.exp}
	There is a constant $c>0$ such that
	for all $t \ge 1$, for all $x \in (0, 2\pi)$,
	$$
	 e^{-\lambda t} \bigl(\sin(x/2)\bigr)^q
	\le f(x,t)\le c  e^{-\lambda t} \bigl(\sin(x/2)\bigr)^q
	$$
	\end {prop}

	\proof
	We can assume that $b>0$ since the case $b=0$ was treated in 
	the previous section.
	Let $Y^x_t$ be as before
	and define  for all $t < \tau^x$
$$
\Phi_t^x :=
\left| g_t'\bigl(\exp(ix)\bigr)\right|
\,
.$$
On $t\ge \tau^x$ set $\Phi_t^x := 0$.
	Note that on $t<\tau^x$
	$$
	\Phi_t^x =\partial_x Y_t^x .
	$$
	By differentiating (\ref{e.loew}) with respect to $z$, we find
	that for $t<\tau^x$
	\begin{equation}\label{e.logder}
	\partial_t \log \Phi_t^x = -
	\frac 1{2\sin^2(Y_t^x/2)}
	\end{equation}
	and hence (since $\Phi_0^x = 1$),
	\begin{equation}\label{Phiis}
	(\Phi_t^x)^b =
	\exp
	\left( - \frac {b}{2} 
	\int_0^t \frac {ds} {\sin^{2}(Y_s^x/2)} \right)
	\ ,
	\end{equation}
	for $t<\tau^x$.
	So, we can rewrite
	$$ 
	f(x,t) = \expect \left[ 1_{{\cal H}(x,t)}
	\exp
	\left( - \frac {b}{2}
	\int_0^t \frac {ds} {\sin^{2}(Y_s^x/2)} \right)
	\right]
	.$$
	Again, it is not difficult to see 
	that 
	 the right hand side of~(\ref{Phiis}) is $0$
	when $t=\tau^x$
	and that
	\begin{equation}\label{sidebd}
	\lim_{x\to 0} f(x,t)=\lim_{x\to 2\pi} f(x,t)=0
	\end{equation}
	holds for all fixed $t>0$.

	Let $F:[0,2\pi]\to\R$ be a continuous function with
$F(0)=F(2\pi)=0$, which is smooth in $(0,2\pi)$, and set
$$
h(x,t)=h_F(x,t):=
\expect \Bigl[(\Phi_t^x)^b\,F(Y_t^x)\Bigr].
$$
By~(\ref{Phiis}) and the general
theory of diffusion Markov processes,
we know that $h$ is smooth in $(0,2\pi)\times\R_+$.
From the Markov property for $Y_t^x$ and~(\ref{Phiis}), it
follows that $h(Y_t^x,t'-t) (\Phi_t^x)^b$ is a local
martingale on $t<\min\{\tau^x,t'\}$.  Consequently,
the drift term of the stochastic
differential $d\bigl(h(Y_t^x,t'-t) (\Phi_t^x)^b\bigr)$
is zero at $t=0$.  By It\^o's formula, this means
\begin{equation}\label{ekol}
\partial_t h = \Lambda h\,,
\end{equation}
where
$$
\Lambda h:=
\frac \kappa 2\, \partial_x^2 h
+\cot(x/2)\, \partial_x h
-\frac b {2 \sin^2 (x/2)}\, h
\,.
$$
We therefore {choose} 
$$
F(x):= \bigl(\sin(x/2)\bigr)^q, 
$$
and note that
$
F(x) e^{-\lambda t}
= h_F$ because both satisfy~(\ref{ekol}) on
$(0,2\pi)\times [0,\infty)$,
 and have the same boundary values.
Finally, one can conclude using the same type of argument 
as in Proposition \ref {p.disco}.
\qed

\section {First consequences}

Recall that  for 
all $t \ge 0$, $d(0,K_t) e^t \in [1/4,1]$. 
Hence, if $\tau_r$ denotes the hitting time of the 
circle of radius $r<1$ by the radial $SLE_\kappa$, then 
$ r e^{\tau_r} \in [1/4,1]$.  
Combining this with Propositions~\ref{p.disco} and 
\ref {p.exp} 
then implies that for all fixed $\kappa >4$, all $b\ge 0$,
 if $\lambda$, $q$ are defined as before, 
there exists two positive finite constants $c_1$ and $c_2$ such that 
for all $r<r_0$,
\begin {equation}
\label {e.exphit}
c_1 r^{\lambda}  (\sin (x/2))^q 
\le \expect \left[ 1_{{\cal H} (x, \tau_r)} 
|g_{\tau_r}'(\exp (ix)) |^b \right]
\le 
c_2 r^{\lambda} (\sin (x/2))^q
\end {equation}
(we used also the fact that 
$| g_t' (\exp (ix)|$ is an decreasing function of $t$).

When $b=1$, one can note that 
$$
l_t := \int_0^{2\pi} dx |g_t' (e^{ix})| 1_{{\cal H}(x,t)}
$$
is simply the length of 
the image under 
$g_t$ of the arc 
$A_t := \partial H_t \cap \partial \U$ 
on the unit circle that have not yet been swallowed by $K_t$.
In particular, if one starts a planar Brownian motion from 
$0$, it has a probability 
$l_t / 2 \pi$ to hit the unit circle on the
arc $g_t (A_t)$.
By conformal invariance of planar Brownian motion, we see 
that 
$l_t / 2\pi$ is also the probability that a planar Brownian motion 
started from $0$ hits the unit circle before hitting $K_t$.
Let $Z$ denote planar Brownian motion, stopped at 
its hitting time $\sigma$ of the unit circle.
Integrating Proposition~\ref {p.exp} for $x \in [0,2\pi]$ 
therefore shows that there exist constants $c_1'$ and $c_2'$
such that (if $K_t$ is radial $SLE_6$) 
\begin {equation}
\label {e.5/4}
c_1' r^{5/4} \le \P [ Z[0, \sigma] \cap K_{\tau_r} = \emptyset ] 
\le c_2' r^{5/4}.
\end {equation}

Combining  these results  with 
Theorem~\ref {t.hulleq}, we see
that these estimates are also valid for
reflected Brownian motions.
In particular, let us now define a reflected 
Brownian motion $\tilde Z$ 
in the unit disc 
as in Theorem~\ref {t.hulleq} (reflected on 
$\partial \U$ with angle $2\pi /3$
``away'' from $\tilde Z_0=1$).
Let $\tilde \sigma_r$ denote its hitting
time of the circle $r \partial \U$.
Then  
there exist constants $c_1$ and $c_2$ such that 
for all $r<1/2$,
\begin {equation}
\label {e.drbm}
c_1 r^{1/4} \le \P [ \tilde Z [ 0, \tilde \sigma_r] 
\hbox {does not disconnect } 0 \hbox { from } -1 ] 
\le c_2 r^{1/4}.
\end {equation}
Similarly, (\ref {e.5/4}) holds if one replaces $K_{\tau_r}$ 
by $\tilde Z [ 0, \tilde \sigma_r]$.
 
We will see in the next chapter that this 
also yields the corresponding 
estimates for (non-reflected) Brownian motions.


\bibc
The material of this chapter is borrowed from \cite {LSW2}, in which the reader can 
find more detailed proofs.
It is possible to compute analogous exponents for chordal SLE. These
``half-plane exponents'' are determined in \cite {LSW1,LSW3}.

Other important exponents are derived in \cite {RS,LSW5,Be3}. As in 
this chapter, the exponents
appear always as leading eigenvalues of some differential operators. 

\chapter {Brownian exponents}

\section {Introduction}

The goal of this chapter is to relate the previous 
computations to the exponents associated to
planar Brownian motion itself (not only to 
reflected Brownian motion).

Suppose that a planar Brownian motion $Z$ is started 
 from $1$.
Let
$\sigma_r$ denote its hitting time of the circle of 
radius $r>0$, and 
let 
$$ p_r := \P [{\cal D}( Z [ 0, \sigma_r ] 
)]
,$$
where ${\cal D} (K)$ denotes the event 
that $K$ does not disconnect 
the origin from infinity.
 Note that by inversion, $p_R = p_{1/R}$ for all $R>1$
(one can map the disk $\{ |z| < R \}$ conformally
on $\{|z| > 1/R \}$ by $z \mapsto 1/z$ and use
conformal invariance of planar Brownian motion).

The strong Markov property and
the scaling property of planar Brownian motion 
imply readily that for all $R,R'>1$,
$$
p_{RR'} 
\le  \P [ {\cal D} (Z[0, \sigma_R]) \hbox { and } 
{\cal D}(Z[\sigma_{R}, \sigma_{RR'} ])
]
 \le  
 p_R p_{R'}
.$$
On the other hand, it is not difficult to see
that 
$$
p_R \ge \P [ Z [ 0, \sigma_R] \cap [-R,0]
= \emptyset ] \ge  c R^{-1/2}
$$
for all $R >1$ and some constant $c$.
Hence, a standard subadditivity argument
implies that 
there exists a constant $\eta \le 1/2$ 
such that 
$$
p_R \approx R^{-\eta}
$$
when $R \to \infty$, where this notation means that 
$\log p_R \sim  - \eta \log R$.
It turned out that there seems to be no direct way to
determine the value of this exponent $\eta$.

Similarly, if $Z^1$ and $Z^2$ denote two independent Brownian motions
started uniformly on the unit circle, then subadditivity implies
the existence of a positive constant $\xi$ such that 
$$
\P [ Z^1 [ 0, \sigma_R^1 ] \cap Z^2 [ 0, \sigma_R^2 ] 
= \emptyset ] \approx R^{-\xi}.
$$

The exponents $\eta$ and $\xi$ are respectively called the 
disconnection exponent and the intersection exponent  
for planar Brownian motion.

\section {Brownian crossings}

We now make some considerations that will help us 
relating the results on reflected Brownian motions 
derived in the previous chapter to the exponents 
$\eta$ and $\xi$.
For simplicity, we first focus on the disconnection 
exponent $\eta$.

Suppose that 
$Z$ denotes a planar Brownian motion that is started from $1$,
and define the random times:
\begin {eqnarray*}
\sigma_r & := & \inf \{ t > 0 \ : \ |Z_t | = r \} 
\\
\sigma_r^\# & := &
\max \{ t <  \sigma_r \ : \ |Z_t | = 1 \} \\
\sigma_r^* &=&
\inf \{ t > \sigma_r^\# \ : \ | Z_t | = 1/2 \}
. \end {eqnarray*}

It is a fairly standard application of the decomposition of 
the path $\Re ( \log Z)$ into excursions away from
the origin to see
that 
\begin {itemize}
\item
The paths
$P^1_r :=(Z_t, t \in  [0, \sigma_r^\#])$ and
$P^2_r := (Z_{t+ \sigma_r^\#}/ Z_{\sigma_r^\#}, t \in [0
, \sigma_r - \sigma_r^\#] )
$ are independent.
\item
The law of $P^3_r :=
( Z_{t + \sigma_r^*} / Z_{\sigma_r^*} ,
t \in [0, \sigma_r - \sigma_r^* ] )$
is identical to the
conditional law  of $(Z_t, t \le \sigma_{2r})$
on the event $E_r:= \{Z_{[0, \sigma_{2r}]} \subset 2 \U\}$.
\end {itemize}
Note also that $\P [ E_r] = \log 2 / \log r$ because
$\log |Z|$ is a local martingale.
We will call $P^2_r$ a
Brownian crossing of the annulus ${\cal A}_r :=
\{ 1 > |z| > r \}$.

When $r'<r$, one can construct a Brownian crossing of the
annulus ${\cal A}_{r'}$ starting from a crossing
$P_r^2$ of the 
annulus ${\cal A}_r$ as follows:
Attach to the endpoint $e_r:=
Z_{\sigma_r}/ Z_{\sigma_r^\#}$  of $P_r^2$ a Brownian motion
started from $e_r$, that is conditioned to hit the circle
of radius $r'$ before the unit circle, and stop it at that
hitting time of the circle of radius $r'$
(note that this event has probability $\log (1/r) /
\log (1/r')$).

We now define the probability $p_r^*$ 
that the crossing does not
disconnect the origin from infinity:
$$
p_r^*
:= {\P [ {\cal D} (P^2_r)]} .
$$
Since a crossing is a subpath of a stopped Brownian motion, it follows
from the a priori lower 
bound for $p_r$ that $p_r^* \ge c r^{1/2}$ for some absolute constant $c$.  

We now define for $\delta>0$,
$$
p_r^* (\delta)
:= 
{\P [ {\cal D} (P_r^2 \cup {\cal B}(1, \delta) \cup
{\cal B} (e_r, \delta) )] },
$$
where ${\cal B }(z,r)$ stands for the ball of radius $r$ around $z$.

The following observations will be useful:

\begin {lemma}
\label {l.separation}
There exists $\delta >0$ and $\eps >0$ such that for all integer $n$,
then for at least $99\%$ of the integers $j \in \{1, \ldots, n\}$,
one has
$$
p_{r_j}^* (\delta) > \eps p_{r_j}^* 
$$
where $r_j = 2^{-j}$.
\end {lemma}

\proof
We only sketch the main ideas of the proof.
First, notice that 
$j \mapsto p_{r_j}^*$ is decreasing in $j$ so that the a priori 
lower bound for $p_{r_j}^*$  
implies that there exists $\eps$ such that for all $n$, then 
for at least $99\%$ of the values of $j$ in $\{ 1, \ldots, n-1\}$,
$2\eps p_{r_j}^* \le  p_{r_{j+1}}^*$
(otherwise $p_{r_n}$ would be too small).
On the other hand, it is easy to see that there exists 
$\delta>0$ such that 
$$
p_{r_{j+1}}^* \le \eps (p_{r_j}^*
- p_{r_j}^* (\delta))  + p_{r_j}^* (\delta)  
.$$
This is due to the fact that one can construct a sample of $P_{r_{j+1}}^2$
by extending the crossing $P_{r_j}^2$ into a crossing of 
$\{ \sqrt {2} > |z| > r / \sqrt {2} \}$
by attaching conditioned Brownian motions to both ends
(and then rescale this into a crossing of ${\cal A}_{r_{j+1}}$). And
if $\delta$ is sufficiently small, then each of the 
attached parts disconnect the ball of radius $\delta$
around their starting point with very high probability.
It therefore follows that ``for $99\%$ of the values of $j$'',
$$
p_{r_j}^* (\delta ) \ge p_{r_{j+1}}^* - \eps
p_{r_j}^* \ge \eps p_{r_j}^* .
$$
\qed
\medbreak

\begin {lemma}
\label {l.start}
For all fixed $\delta$, for some constant $c=c(\delta)$,
$$\P [ P_r^1 \subset {\cal B} (1,\delta /2) ]
\ge \frac {c}{\log (1/r)}.
$$
\end {lemma}
\proof
With positive probability, $Z$ hits the circle of radius
$1- \delta/4$ around $0$ before $\partial {
\cal B} (1, \delta/2)$. Then, if this 
is the case,  with probability $
\log (1/(1-\delta/4)) / \log (1/r)$ it hits the circle
of radius $r$ before going back to the unit circle.
\qed

\section {Disconnection exponent}

We now use combine these considerations with the 
computation of the exponents for reflected Brownian motion
to prove the following result:

\begin {theorem}
\label {t.bmdisco}
One has $\eta = 1/4$. Furthermore, there exist two
constants $c_1$ and $c_2$ such that for all $R>1$,
$$
c_1 R^{-1/4} \le p_R \le c_2 R^{-1/4} .
$$
\end {theorem}

As we shall see later, it is important to 
have estimates  ``up-to-constants'' as in this Theorem
(rather than $\approx$)
 in order to make the link with Hausdorff dimensions.

\proof
By inversion, this is equivalent
to
corresponding result for small $r$ i.e., that
for all $r <1$,
\begin {equation}
\label {e.24}
c_1 r^{1/4} \le p_r \le c_2 r^{1/4}.
\end {equation}
In order to compare $p_r$ to $\tilde p_r$ (this is the
non-disconnection probability for reflected Brownian motion
that was defined at the end of the previous chapter where
we proved that it is close to $r^{1/4}$), 
we will in fact compare both to $p_r^*$.

First, one can notice using the previous lemma that
\begin {eqnarray*}
p_r 
& \ge & 
\P [ {\cal D} ( P_r^2 \cup {\cal B} (1, \delta) ) \cap \{ P_r^1 \subset 
{\cal B} (1, \delta/2 ) \} ]
\\
& \ge &
p_r^* (\delta) \times \frac { c}{\log (1/r)}
\end {eqnarray*}
for some constant $c$ which is independent of $r<1/2$.
The same argument can be adapted to the reflected
Brownian motion $\tilde Z$. Hence, ``for $99\%$ of $j$'s'', 
$$
p_{r_j} \ge \frac {c p_{r_j}^*}{j}
\hbox { and }
\tilde p_{r_j} \ge \frac {c p_{r_j}^*}{j} $$
for some universal constant $c$.

On the other hand, 
let us now define 
inductively the stopping times:
$\rho_0=0$ and for all $n \ge 0$,
\begin {eqnarray*}
\tau_n &:=& \inf \{ t > \rho_{n} \ : \ |Z_t| = 1/2 \} \\
\rho_{n+1} &:=& \inf \{ t >  \ho_n \ : \ | Z_t | =1 \}
\end {eqnarray*}
the successive times of downcrossings and upcrossings
between the two circles $\{ |z| =1 \}$ and $\{ |z| = 1/2\}$.
Let $N_r$ denote the number of upcrossings before $\sigma_r$. In other
words,
$$
N = N(r) := \max \{ n \ge 0 \ : \ \rho_n < \sigma_r \}.
$$
Note that the probability that a Brownian motion started 
on the circle $\{ |z| = 1/2 \}$ hits $\{ |z| = r \}$ 
before the unit circle is $c_r :=\log 2 / \log (1/r)$, because 
$\log |Z|$ is a local martingale.
Hence, $\P [ N_r \ge n ] = (1-c_r)^n$.
For each $n \ge 0$, the probability that 
$Z [  \rho_n , \tau_n ]$ disconnects $0$ from the unit circle 
and does not hit the circle of radius $1/4$ is strictly positive
(and independent from $n$). 
Note that if $Z [ 0, \sigma_r]$ does not disconnect the
origin from the unit circle, then
for all $n \le N$, $Z[\rho_n ,\tau_n ] $
does not disconnect the origin from the unit circle, and
$Z [ \tau_{n} , \sigma_{n,r} ]$ doesn't either,
where
$$ 
\sigma_{n,r} = \inf \{ t > \tau_{n}  \ : \ |Z_t| = r \}
.$$
It follows that for some absolute constant $c>0$, 
\begin {eqnarray*}
p_r & \le &
 \sum_{n \ge 0} (1-c)^n (1-c_r)^n \P [ {\cal D}(Z [  \tau_{n}
, \sigma_{n, r} ])
] \\
& \le &
 \frac {p_{2r}^*} {c c_r} \\
&\le &
\frac { \log 2 }{c} \times \frac {p_{2r}^*}{\log {1/r}} 
\end {eqnarray*}
A close inspection at the proof  
actually shows that the very same proof goes through if one
replaces the Brownian motion $Z$ by the reflected Brownian motion
$\tilde Z$. Hence, 
for some 
absolute constant $c'$,
$$
 p_r 
\le c' \frac {p_{2r}^*}{\log (1/r)}
\hbox { and }
\tilde p_r                   
\le c' \frac {p_{2r}^*}{\log (1/r)}.
$$
Putting the pieces together, we see that ``for $98 \%$ of $j's$'',
$$
\tilde p_{r_j} 
\le c_1 \frac {p_{2r_j}^*}{\log (1/r_j)}
\le c_2 p_{2r_j}
\le c_3 \frac {p_{2r_j}^*}{\log (1/r_j)}
\le c_4 \tilde p_{4r_j}.
$$
But we know that $r^{-1/4} \tilde p_r$ is bounded and
bounded away from zero. It therefore follows that 
for some absolute constants $c_1$ and $c_2$
and at least $98\%$ of the $j$'s,
$$
c_1 r_j^{1/4} \le p_{r_j} \le c_2 r_j^{1/4}.
$$
It then remains to get rid of the last $2\%$ of ``bad'' 
values of $j$. This can be done by pasting together 
``good'' configurations that are
``well-separated at the end'' 
of the 
annuli
$\{ 1 > | z| > r_{j_1} \}$ 
and $\{r_{j_1} > | z| > r_{j_1+j_2} \}$,
where $j_1$ and $j_2$ are ``good'' values such that $j_1+j_2 = j$.
See for instance \cite {LSWup} for more details.
\qed

\section {Other exponents}

The previous proofs need to be somewhat adjusted to 
show the corresponding result for the intersection exponent
$\xi$ (things are  more complicated due to the 
fact that there are two Brownian motions to take
care of, but no really new ideas are needed):
 
\begin {theorem}
\label {t.bmint}
One has $\xi = 5/4$. 
Furthermore, there exist two 
constants $c_1$ and $c_2$ such that for all $R>1$,
$$
c_1 R^{-5/4} \le \P [ Z^1 [0, \sigma_R^1] 
\cap Z^2 [ 0, \sigma_R^2 ] = \emptyset ] \le c_2 R^{-5/4} .
$$
\end {theorem}

Actually, it is possible to derive the value of many other 
exponents. For instance, 
suppose that 
$Z^1, \ldots , Z^k, \ldots$ are independent planar Brownian motions
started uniformly on the unit circle, and denote by 
$\sigma^1_R, \sigma_R^2, \ldots$ their respective hitting times
of the circle $R \partial \U$, then:

\begin {theorem}
\label {t.big}For all $k \ge 1$, there exist 
constants $c_1, c_2$ such that for all $R > 1$,
$$
c_1 R^{-\eta_k} \le \P [ {\cal D} ( Z^1 [ 0, \sigma_R^1] \cup
\cdots \cup Z^k [ 0, \sigma_R^k ] ) ] 
\le c_2 R^{-\eta_k}
$$
and
$$
c_1 R^{-\xi_k}
\le \P [ \hbox {The sets }
Z^1 [ 0, \sigma_R^1], \ldots , Z^k [0, \sigma_R^k] 
\hbox { are disjoint}]
\le c_2 R^{-\xi_k},
$$
where
$$
\eta_k=
\frac { ( \sqrt {24 k +1 } - 1 )^2 - 4 }{48}$$
and
$$
\xi_k = \frac { 4 k^2 -1 }{12}.
$$
\end {theorem}

The proof of these results is however more involved.
For other results and generalizations, see~\cite {LSW1,LSW2,LSW3}.
For instance, one can make sense of a continuum of exponents, 
or study intersection exponents for Brownian motion in a 
half-plane.

Let us mention that an instrumental role is also played in
the definition and determination of the exponents in Theorem \ref {t.big}
by the
critical exponents associated to non-intersection
events in a half-space. For instance, the half-space analog of 
the intersection exponent $\xi$ is:

\begin {theorem}
If $Z^1$ and $Z^2$ are defined as before.
Define
$$
q_R
:=
\P [
Z^1 [ 0, \sigma_R^1] \cap Z^2 [ 0, \sigma_R^2 ]
= \emptyset
\hbox { and }
Z^1 [ 0, \sigma_R^1] \cup Z^2 [ 0, \sigma_R^2 ]
\subset \H ].$$
There 
exist two constants $c_1$ and $c_2$ such that for all $R>1$,
$$
c_1 R^{-10/3}
\le
q_R 
\le c_2 R^{-10/3}.
$$
\end {theorem}

There is a close
relation between all these exponents (disconnection,
in the whole space, in the half-space), see \cite {LW1}.
The critical exponents in the half-space can be determined
in a similar way than the the whole-space exponents:
First one computes the ``derivative''  exponents associated
to chordal SLE. Then, using the identification between
chordal $SLE_6$ and reflected Brownian motion, one
transfers the SLE results into Brownian motion
results.
For the statements and proofs of all these
``half-space exponents'', see \cite {LSW1,LSW3}.
In order to get the value of all $\eta_k$ exponents, one then
uses the fact that a family of generalized exponents is
analytic, see \cite {LSWan} for more on this.

\medbreak
It has also been proved (using strong approximation of
simple planar random walks by Brownian motions) that
these exponents  describe the    
probabilities of the corresponding events for
planar simple random walks
 (see \cite {BL1, CM, LP1, LP2}).
For instance, if $S^1$ and $S^2$ denote two independent
simple random walks starting from neighbouring points,
then
$$
\P [ S^1 [0,n] \cap S^2 [0,n] = \emptyset ]
\approx
n^{-\xi/2} = n^{-5/8}
$$
when $n \to \infty$
(up-to-constants hold as well). The exponent is here
$\xi/2$ because we used here the parametrization in time and not
in space.
It is worthwhile stressing that it seems that to prove this
result that seems of combinatorial nature,
 one has to understand and use conformal invariance of
planar Brownian motion, its relation to $SLE_6$ as well
as the properties of $SLE_6$.

\section {Hausdorff dimensions}

In series of papers \cite {L2,L3,L4,L5} (before the
mathematical 
determination of the exponents in \cite {LSW1,LSW2,LSW3}),
 Lawler showed how to
use such up-to-constants estimates to estimate the 
Hausdorff dimension of various interesting 
random subsets of the planar Brownian curve in terms 
of the corresponding exponents.

More precisely, let $(Z_t, t \ge 0)$ denote a planar Brownian motion.
Then, we say that 
\begin {itemize}
\item 
The point $z=Z_t$ is a cut-point if $Z[0,t] \cap Z(t, 1] = \emptyset$.
\item
The point $z=Z_t$ is a boundary point if ${\cal D} ( Z[0,1] - z )$ i.e.
if $Z[0,1]$ does not disconnect $z$ from infinity.
\item
The point $z=Z_t$ is a pioneer point if ${\cal D} (Z [0,t] -z)$. 
\end {itemize}
Note that, loosely speaking, near $z=Z_t$, there are two 
independent Brownian paths starting at $z$: The future $Z^1 := 
(Z_{t+s}, s \in [0, 1-t])$ and the past $Z^2 := 
(Z_{t-s}, s \in [0,t] )$.
Furthermore, $z=Z_t$ is a cut-point if $Z^1 \cap Z^2 = \{z\}$,
$z$ is a boundary point if $Z^1 \cup Z^2$ do not disconnect $z$
from infinity and $z$ is a pioneer point if $Z^2$ does not disconnect
$z$ from infinity. 
Hence, the previous theorems enable us to estimate the 
probability that a given point $x \in \C$ is in the $\eps$-neighbourhood
of a cut-point (resp. boundary point, pioneer point).
Independence properties of planar Brownian paths then make it 
also possible to derive second moment estimates 
(i.e. the probability that two given points
$x$ and $x'$ are both in the $\eps$-neighbourhood
of such points) and to obtain the following result:

\begin {theorem}
\label {t.dim}
$\ $
\begin {itemize}
\item
The Hausdorff dimension of the set of cut-points is almost 
surely $2 - \xi$.
\item
The Hausdorff dimension of the set of boundary points is almost 
surely
$2- \eta_2$.
\item
The Hausdorff dimension of the set of pioneer points is almost surely
$2- \eta $.
\end {itemize}
\end {theorem}

Recall that $2-\xi = 3/4$, $2-\eta_2 = 4/3$, $2- \eta =7/4$.
Similar results hold for various other random subsets of the 
planar curve. 
We choose not to give the proof of this theorems in these
lectures since they are more using features of planar 
Brownian motion rather than $SLE_6$, but here
is a brief sketch in the case of the pioneer points.

\medbreak
\noindent
{\bf Sketch of the proof.}
Let ${\cal P}$ denote the set
of pioneer points on $Z[0,1]$.
Theorem~\ref {t.bmdisco}
roughly shows that for each $z$, the probability that 
$Z$ comes $\eps$-close to $z$ without disconnecting 
$z$ from infinity is comparable to $\eps^{1/4}$.
It follows that the expectation 
of the number $N_\eps$ of $\eps$-balls that 
are needed in order to cover ${\cal P}$ is 
comparable to (i.e. up-to-constants 
away from) $\eps^{-2+1/4} = \eps^{-7/4}$.
This in fact already shows that the 
Hausdorff dimension of ${\cal P}$ can a.s. not be
larger than $7/4$.

On the other hand, one has
good bounds on the second moment of $N_\eps$:
This is due to the fact that for two points $x$ and $x'$  
with $|x-x'|=r$ to be $\eps$-close
to pioneer points, then the following three events
must occur before time one:
\begin {itemize}
\item
$Z$ reaches ${\cal B} (x, 2r)$ without
disconnecting $x$
\item
$Z$ crosses the annulus 
$\{ z : \eps < | z-x| < r/2 \}$
without disconnecting $x$
\item
$Z$ crosses the annulus 
$\{ z : \eps < | z-x'| < r/2 \}$
without disconnecting $x'$.
\end {itemize}
Hence, it follows that
$\expect [ N_\eps^2 ] \le cst \times \eps^{-7/2}
\le cst E [ N_\eps]^2$.
Standard arguments can then be used   
to deduce from this that with positive
probability, the dimension of ${\cal P}$
is not smaller than $7/4$.
A zero-one law can finally be used to 
conclude that the dimension is a.s. equal
to $7/4$. See e.g. \cite {Lbuda} for details.
\qed

\bibc
The fact that 
one probably had to compute the value of the Brownian exponents
via an universality argument 
using another model (that should be closely related to 
critical percolation scaling limits) first appeared in \cite {LW2}.
The mathematical derivation of 
the value of the exponents was performed in the series of papers
\cite {LSW1,LSW2,LSW3,LSWan}. 
The properties of SLE that were later derived in \cite {LSWrest} 
enable to shorten some parts of some proofs, but it 
seems that analyticity of the family of generalized 
exponents derived in \cite {LSWan} can not be by-passed
for all exponents (for instance, it seems that it is 
needed to determine the exponent describing the 
probability that the union of three Brownian motions
does not disconnect a given point).
It can however be by-passed for those exponents that we have 
to focus on i.e.,  $\eta, \eta_2, \xi$.

Lemma \ref {l.separation} is a ``separation Lemma'' of the type
that had been derived by Lawler in the series of papers
relating the Hausdorff dimensions to the exponents
\cite {L1,L2,L3,L4,L5}. The proof presented here is 
adapted from the proof of the analogous but more
general results for the other exponents in \cite {LSWup}.
A good reference for the relation between Brownian exponents 
and Hausdorff dimensions is Lawler's review paper \cite {Lbuda}.
See also, Beffara \cite {Be1, Be2}.

Determining the Hausdorff dimensions of subsets of the SLE
processes is a difficult question.
Rohde-Schramm \cite {RS} have shown that the 
dimension of the SLE generating curve is not larger
than $1+\kappa /8$. It was conjectured to be a.s. equal to that
value (for $\kappa \le 8$). This has been proved to hold 
for the special values $\kappa = 8/3$ and $\kappa=6$, making use
of the locality and restriction properties (see 
\cite {LSWrest}, Beffara \cite {Be2}).
It now seems that Beffara \cite {Be3} managed to 
prove the general conjecture. 
 
The value of most of these exponents had been predicted/conjectured
before: Duplantier-Kwon \cite {DK} had predicted the 
values of $\xi_k$ using non-rigorous conformal field theory
considerations, Duplantier \cite {Dqg} more recently used also 
the so-called ``quantum gravity'' to predict the values
of all exponents. 
The fact that the dimension of the Brownian boundary was $4/3$ 
was first observed visually 
and conjectured by Mandelbrot \cite {Ma}.
Before the proof of this conjecture,  some rigorous bounds had been 
derived, for instance that the dimension of the 
Brownian boundary is strictly larger than $1$ and strictly smaller
than $3/2$ (see \cite {BJPP, BL2, Wbd}).

\chapter {SLE, UST and LERW}

\section {Introduction, LERW}

In the next two chapters, we will survey the rigorous results
that show that for some values of $\kappa$, $SLE_\kappa$ is 
indeed the scaling limit of discrete models.
There are at present only three values of $\kappa$ for
which this is the case: $\kappa=2$ is the scaling limit of LERW, 
$\kappa=6$ is the scaling limit of percolation cluster interfaces, 
and $\kappa=8$ is the scaling limit of the uniform spanning tree contour.

In all three cases, the convergence to SLE is derived 
as a  
consequence of three facts:
\begin {itemize}
\item
The ``Markovian'' property holds in the discrete case
(this is usually a trivial consequence of the 
definition of the microscopic model).
\item
Some macroscopic functionals of the model converge to  
conformally invariant quantities in the scaling limit (for a 
wide class of domains).
\item
One has ``a priori'' bounds on the regularity of the 
discrete paths.
\end {itemize}

Before going into more details, let us state the convergence 
theorem in the case of LERW that was presented in the 
introductory chapter:
Consider $\gamma^\delta$ the 
(time-reversal of the) loop-erasure of a simple random walk in 
$D \cap \delta \Z^2$, started from $0$ and stopped at the first
exit time of the simply connected (say, bounded) domain $D$.
Let $\gamma$ denote a radial $SLE_2$ in the unit disc started
uniformly on the unit circle (and aiming at $0$).
 Let $\Phi$ denote 
a conformal map from $\U$ onto $D$ that preserves $0$.
We endow the set of 
paths with the metric of uniform convergence modulo 
time-reparametrization:
$$
d ( \Gamma , \Gamma' ) 
= \inf_{\varphi} \sup_{t \ge 0} | \Gamma (t) - \Gamma'(\varphi(t)) |
$$
where the $\inf$ is over all increasing bijections
$\varphi$  from $[0,\infty)$ 
into itself.
Then,
\begin {theorem}
\label {t.lerw}
The law of $\gamma^\delta$ converges weakly when $\delta \to 0$
to the law of $\Phi (\gamma)$.
\end {theorem}

Actually, one can also use the convergence result to justify the fact that
$SLE_2$ is a simple path.
Instead of giving the
 basic ideas of the proof of this theorem, we will 
focus on a closely related problem:
 The uniform spanning trees scaling limit.

\section {Uniform spanning trees, Wilson's algorithm}

Suppose that a 
connected finite graph $G = ( V, E)$ is given ($V$ is the set of vertices and 
$E$ is the set of edges). We say that the subgraph $T \subset E$ is a 
spanning tree if it contains no loop,
 and if it has only one connected component.
We then define the uniform spanning tree as the 
uniform measure on the set of spanning trees.
For any two fixed points $a$ and $b$ in $G$, and any spanning tree $T$, there
exists a unique simple path in $T$ that joins $a$ to $b$ (it exists because 
$T$ has one connected component, it is unique because $T$ has
no loops). Hence, if $T$ is picked according to the UST measure,
this defines a random path $\gamma$ from $a$ to $b$. 
The following result had first been observed by Pemantle \cite {Pem}:

\begin {prop}
\label {p4}
The law of $\gamma$ is that of the loop-erasure of simple random walk
on $G$ started at $a$ and stopped at its first hitting of $b$.
\end {prop}

\begin{figure}
\centerline{\includegraphics*[height=3in]{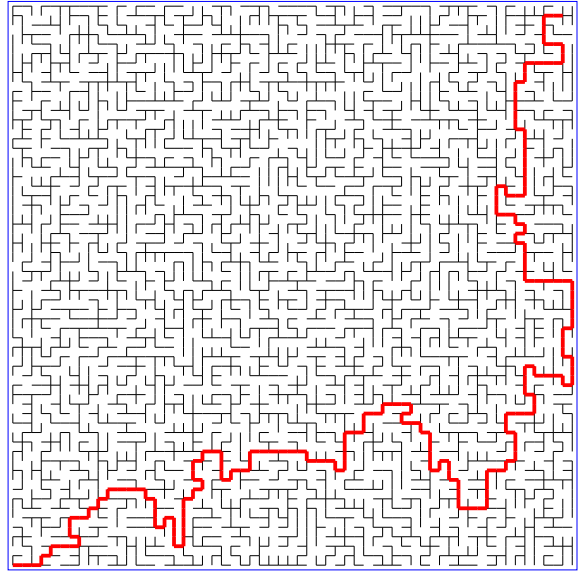}}
\caption{\label{f.lerw-in-ust}
A loop-erased walk as a subpath of the UST}
\end{figure}

This shows that LERW and UST are very closely related. Actually, it
turns out that an even  stronger relationship hold:
Suppose that an ordering of the vertices $v_0, v_1, \ldots, v_m$ 
of $G$ is given.
Define inductively the sets $A_m$ as follows: $A_0 = \{ v_0 \}$,
and for all $j \le m$,
$A_j = A_{j-1} \cup \gamma_j$
where $\gamma_j$ is the loop-erasure of a 
random walk started from $v_j$ and stopped at its first 
hitting of $A_{j-1}$.
Clearly, in this way, $A_m$ is a (random) tree that contains all 
vertices: It is a spanning tree.

\begin {prop}[Wilson's algorithm]
\label {p.Wilson}
The law of $A_m$ is the uniform spanning tree measure.
\end {prop}

Note that this algorithm yields a natural extension of uniform
spanning trees (or forests) in infinite graphs (see e.g. \cite {BLPS}
and the references therein for more on this subject).

\proof
One can derive this result using the explicit formulas that we 
derived in the introductory
 chapter for loop-erased random walks:
Indeed, it follows readily from the definition 
and the symmetry of the function 
$F$ that was defined there, and the fact that (since we are 
considering simple random walks), the transition probabilities
$p(x,y)$ are simply equal to $1/d_x$ where $d_x$ is the number of neighbours of
$x$), that for any possible spanning tree $T$,
$$
\P [ A_m = T ] 
=  F (v_1, v_2, \ldots, v_m ; \{ v_0 \} )
 \prod_{j=1}^m (1/ d_{v_j}) .
$$
This quantity is the same for all $T$: The law of $A_m$ is uniform.
\qed

\medbreak

Hence, 
if LERW has a conformally invariant scaling limit 
then UST also has a conformally invariant scaling limit (in a 
rather weak sense though, such as: for all $k$
given fixed points, the ``finite subtree
that go through these points'' converges in the scaling limit).

There is another way to encode planar trees that goes as follows.
 Suppose for instance that we are looking at a spanning tree of a bounded 
``simply connected'' 
graph $G \subset \Z^2$. Then, one can associate to each tree the contour 
of the tree which is a simple closed curve living on  a subset 
$G^\#$ of the lattice $(1/4 + \Z /2)^2$.
It is easy to see that (under mild assumptions on the domain), this curve visits every 
point of $(1/4 + \Z /2)^2$ that is close to the vertices of $G$. If the tree is chosen
according to the uniform measure on
 spanning trees, then the contour is chosen according
to the uniform measure on space-filling simple 
closed
curves in this graph $G^\#$. 

Hence, it is natural to study the behaviour of this space-filling curve in the
scaling limit.
In order to obtain  SLE (and not a 
closely related object that we would have to
define first)  it is (slightly) more convenient to consider a  
variant of the previously defined space-filling curve. 

\begin{figure}
\centerline{\includegraphics*[height=3in]{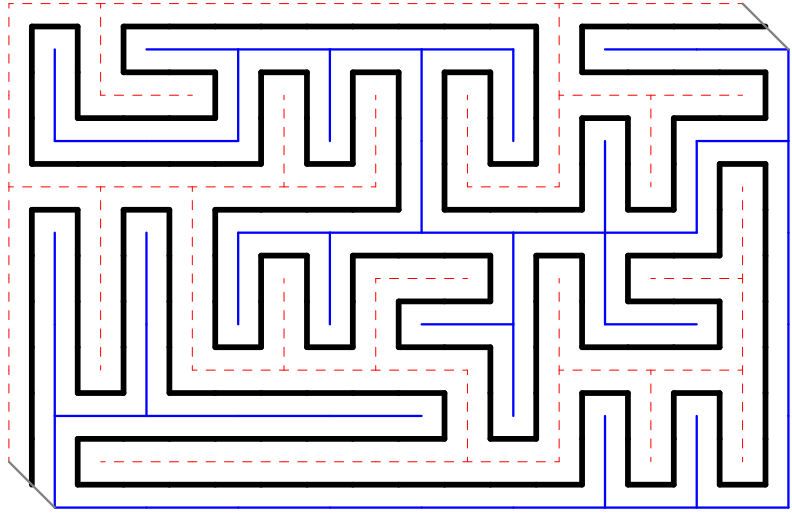}}
\caption{\label{f.dual}
The wired tree, the dual tree, the Peano curve.}
\end{figure}

More precisely, suppose that a certain connected graph of 
$(1/4 + \Z /2)^2$ is given together with two distinct 
``boundary
points'' $a$ and $b$. Then (for a suitable class
of ``admissible'' graphs), one is interested in the uniform measure
on simple space-filling curves $\eta$ from $a$ to $b$ in the graph (i.e. 
paths from $a$ to $b$ that visit all vertices exactly once).
An example of ``admissible'' graphs is given by the graph
obtained from removing from $G^\#$ a part of a 
simple closed space-filling curve $\gamma$.
This time, there is a one-to-one correspondence
between the family of simple space-filling curve 
$\eta$ (from $a$ to $b$) and the set of spanning trees in a certain
subgraph $G$ of $\Z^2$ obtained by wiring one part
of the boundary between $a$ and $b$ (i.e. by conditioning the
tree to contain this part of the boundary).
This is best seen on pictures, and not difficult to 
understand heuristically, but it is  somewhat messy to 
formulate precisely, so we will omit the precise statements  
here (see e.g., \cite {LSWlesl} for more details). 
 
Note that in this set-up, the Markovian type property for 
$\eta$ is immediate: If one conditions on the first step 
$\eta (1)$ of $\eta$, then the law of $\eta(1), \ldots, \eta(n)=b$
is simply the uniform measure on the space filling curves from
$\eta(1)$ to $b$ in the remaining graph.

\section {Convergence to chordal $SLE_8$}

Suppose that $D$ is a simply connected bounded planar domain with 
$C^1$ boundary and let $a,b$ denote two distinct points on 
$\partial D$.
For each $\delta$, we associate in  a ``suitable
approximation'' of $D \cap \delta Z^2$, denoted by $D_\delta$,
and the two boundary points $a_\delta$ and $b_\delta$ close
to $a$ and $b$.
We define $\eta^\delta$,  a uniformly chosen space-filling curve
from $a_\delta$ to $b_\delta$ in $D_\delta$.
 
\begin {theorem}
\label {th.ust}
When $\delta \to 0$, the law of $\eta^\delta$ converges 
weakly to that of a space-filling continuous 
 path $\eta$, such that the 
law of $(\eta [0,t], t \ge 0)$ is 
(up to time-change) that of chordal $SLE_8$ in $D$ from $a$ to $b$.
\end {theorem}

\begin{figure}
\centerline{\includegraphics*[height=2in]{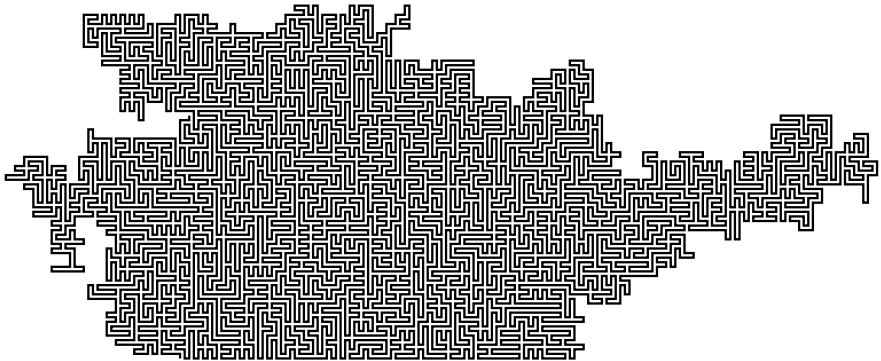}}
\caption{\label{f.peano}A sample of the beginning of the Peano curve.}
\end{figure}

\noindent
{\bf Some rough ideas from the proof.}
A first step is to obtain regularity estimates on the 
(discrete) random space-filling curve. This shows that the families
of probability measures defining $\eta^\delta$ is tight in 
an appropriate sense and therefore has subsequential limits.
These estimates have been derived in \cite {S1} (see 
also \cite {AB,ABNW}) and the basic
tools are Wilson's algorithm and estimates for simple random 
walks. This is not easy, and we refer to \cite {S1} for
details.  
Hence, one can work with a given decreasing sequence $\delta_n \to 0$ 
such that the law of $\eta_{\delta_n}$ converges towards that of
a random curve $\eta$, and one has to show that $\eta$ is in 
fact chordal $SLE_8$.

Let us first work on the discrete level. Suppose that 
$z_\delta$ is some discrete lattice approximation of $z \in D$ and that 
$c_\delta$ is some discrete lattice approximation of $c \in \partial D$ 
that is on the wired part of the boundary of $D$.
Let $P_1^\delta$ denote the part of the wired
boundary of $D_\delta$ 
 which is between $a_\delta$ and $c_\delta$, and let $P_2^\delta$ denote the 
part of the wired boundary which is between $c_\delta$ and $b_\delta$
(and $P_1$, $P_2$ are defined similarly in $D$).

We consider the event
$E^\delta$
that there exists a path in the corresponding tree that goes from 
$z_\delta$ to $P_2^\delta$ without touching $P_1^\delta$. 
By Wilson's algorithm, we see that $\P [ E^\delta (c,z)]$ 
is the probability that simple random walk on $D_\delta$ hits
$P_2^\delta$ before $P_1^\delta$.
One first key-observation is that when $\delta$ goes to 
$0$, the
probability of this event can be controlled in a rather uniform 
way: Uniformly over some suitable choices of $z$, $c$, $a$, $b$ and $D$,
it converges towards the probability that a Brownian motion in $D$
that is orthogonally reflected on the `free' part of $\partial D$, hits
$P_2$ before $P_1$.
This is a conformally invariant quantity.
Mapping $D$ onto the upper half-plane by some given
fixed mapping $g$
in such a way that 
$g(b)=\infty$, we see that 
$$
\lim_{\delta \to 0}
\P [ E^\delta  ] 
:=
h ( A, C, Z)  = F \left( \frac { Z-A}{C-A } \right)
,$$
where 
$$
F(r e^{i \theta})
= \frac {1}{\pi}
\tan^{-1} \left( \frac {1-r}{2 \sqrt {r} \sin \theta/2 } \right)
$$
and $A=g(a), Z=g(z), C=g(c)$.
This function $F$ can be computed for instance by first using reflection 
so that this probability is the probability that (non-reflected)
Brownian motion in the complex plane, started from $g(z)$ hits 
$[g(c), \infty)$ before $[g(a), g(z)]$, then to use the 
map 
 $x \mapsto (\sqrt x-\sqrt z)/(\sqrt x+\sqrt z)$ from $\C \setminus 
[g(a), + \infty)$
onto the unit disk and to look at the length of the image of $[g(c), + \infty)$
on the unit circle).

At each step $n$, define the conformal map
$\phi_n^\delta$ from  
a continuous approximation $D_n^\delta$ of $D^\delta \setminus
\eta[0,n]$ onto the upper half-plane  that is characterized 
by 
$\phi_n^\delta (x) - \phi_0^\delta (b)
= o( 1)$ when $x \to b$.
We then define $t_n^\delta$ to be the ``size''
$a(\phi_0^\delta (\eta^\delta [0,n]))$ of $\phi_0^\delta 
(\eta[0,n])$ and we put
$$
A_n^\delta
= \phi_n^\delta (\eta_n),
\
C_n^\delta = \phi_n^\delta (c_\delta)
\hbox { and }
Z_n^\delta = \phi_n^\delta (z_\delta).
$$ 

Suppose now that $\eps >0$ is small but fixed.
If one stops the 
uniform Peano curve at the first step $N$, at which 
either $|A_n^\delta-A_0^\delta|$ reaches $\eps$ or $t_n^\delta$
 reaches $\eps^2$,
(if $c$ and $z$ are not close to $a$),
then one does not yet know
whether $E^\delta$ holds or not. 
In fact the conditional
probability is just
equal to 
$$
\P [ E^\delta (c_\delta,z_\delta ,\eta_N, b_\delta, D_N^\delta)]
.$$
Hence, 
$$
\expect [\P [ E^\delta (c_\delta,z_\delta
,\eta_N, b_\delta, D^\delta_N)]
=
\P [ E^\delta  
] 
. $$
The right-hand side is close to 
$h(A,C,Z)$ and the right-hand side is close to 
$\expect [h (A_N^\delta, C_N^\delta, Z_N^\delta)]$
(in a uniform way as $\delta$ goes to $0$).
In fact, one can prove that 
$$
\expect [ h (A_N^\delta,C_N^\delta, Z_N^\delta)] = \expect [ h (A_0, C_0, Z_0)]
+ O (\eps^3).
$$
It turns in fact out, that the conformal map $\Phi_N^\delta$
is very close to the (properly normalized)  conformal map 
from $D \setminus \eta [0,N]$ onto $\H$ (i.e. removing the slit or 
the ``tube'' does not make much difference when $\delta$ is small).
In particular, when $\eps$ is small (and $\delta$ very small),
Loewner's equation shows that 
$$
(Z_N^\delta - Z_0)  = \frac {2t_N^\delta}{Z_0 -A_0 }
+ O (\eps^3) \hbox { and }
(C_N^\delta - C_0) =
 \frac {2t_N^\delta}{C_0-A_0}
+ O (\eps^3). $$
Hence, one can Taylor-expand $h$ in the previous estimate, so that
\begin {eqnarray*}
\lefteqn {\frac 12 
\expect [ (A_N^\delta-A_0)^2 ]
\partial_A^2 h(A_0, C_0, Z_0)
+
\expect [ A_N^\delta - A_0] 
\partial_A h (A_0, C_0, Z_0)
} \\
&&+ 2 \expect [t_N^\delta] 
\left( \frac {\partial_C h (A_0, C_0, Z_0) }{C_0 - A_0 } 
+ \frac {\partial_Z h (A_0, C_0, Z_0) } {Z_0 - A_0 }
\right)
= 
O ( \eps^3 ) .
\end {eqnarray*}
Using the explicit expression of $h$ as well as the fact that this
holds for various values $c$ and $z$ yields that in fact:
$$
\expect [ A_N^\delta- A_0 ] = O (\eps^3)
\hbox { and }
\expect [ (A_N^\delta-A_0)^2] 
= 8 \expect [t_N^\delta] + O (\eps^3).
$$
One can iterate this procedure using inductively 
defined stopping times
$N_2, N_3, \ldots$, and 
one can then use  this as a  seed to show that it is possible 
to find a Brownian motion $B$ such that $A_n^\delta$ remains close
to $B_{8 t_n^\delta}$, and then, after some additional work 
can be improved into the convergence theorem.
\qed

\medbreak

As the reader can see, this is only a very sketchy outline of a 
fairly long and technical proof. For details, see \cite {LSWlesl}.

\section {The loop-erased random walk}

The strategy of the proof of Theorem \ref {t.lerw} follows roughly the 
same lines. 
One has to identify  a conformal invariant quantity that appears in 
the scaling limit of LERW and that plays the role of the
probability of the  events $E$ in
the case of the uniform Peano curve.
The macroscopic quantities 
 that are used are related to 
the mean number of visits to a given point $z$ by the 
simple random walk started from $0$ and conditioned 
to leave the domain at the same point as the LERW.
See \cite {LSWlesl} for details.

\bibc
The convergence results presented in this chapter 
are proved in \cite {LSWlesl}, where the reader can
find more details.
For an introduction to LERW and UST, see for instance 
\cite {Lyons,LLERW}.
Rick Kenyon \cite {K0,K1} 
had proved that LERW (and UST's) have conformally
invariant features exploiting the relation between UST and 
dimer models (and some explicit computations). He also 
managed to determine directly 
(without using $SLE_2$ or $SLE_8$) \cite {K2,K3}
the
value of various critical exponents related to LERW and UST that
had been conjectured by Majumdar and Duplantier \cite {Maj,Dle}.
For instance, he showed that the expected length of a LERW
from $0$ to the boundary of the unit disc
 on the lattice $\delta \Z^2$ 
is of the order $\delta^{-5/4}$.
See also, Fomin's paper
 \cite {Fom} for another approach to some of these exponents.

In the recent preprint \cite {Ko}, Gady Kozma gives a completely different 
approach and justification to the existence of a scaling limit of LERW (that does 
not seem to use conformal invariance or SLE).

\chapter {SLE and critical percolation}

\section {Introduction}

\begin{figure}
\centerline{\includegraphics*[height=3in]{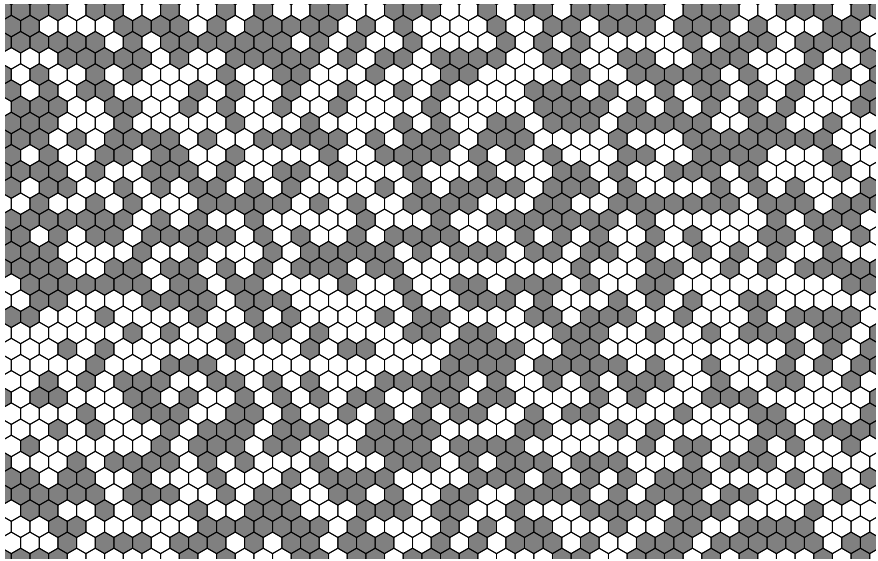}}
\caption{\label{f.crossing}Is there a left to right crossing
of white hexagons?}
\end{figure}

Consider a planar ``periodic'' lattice such that simple random 
walk on that lattice converges to planar Brownian motion. For
convenience, let us limit our discussion to the square lattice
and to the triangular lattice. 
Fix $p \in [0,1]$, and for each site of the lattice, decide
that with probability $p$, the site is open (with probability $1-p$,
it is therefore closed), and do that independently for all sites 
of the lattice.
One is interested in the properties of the connected components
(or ``clusters'') of open sites. It is now classical  (see e.g.,
\cite {Gr}
for an introduction to percolation) that there exists 
a critical value $p_c \in (0,1)$ such that:
\begin {itemize}
\item
If $p \le p_c$, there exists a.s. no infinite open cluster
(note that in dimension greater than $2$, the non-existence of an infinite 
open cluster at $p_c$ is still an open problem).
\item
If $p < p_c$, there exists a positive $\xi (p)$ such that when 
$n \to \infty$, the probability that $0$ is in the same 
connected component than $(n,0)$ decays exponentially fast, 
like $\approx \exp ( - n /\xi (p))$ (the positive quantity 
$\xi (p)$ is called the correlation length).
\item
If $p > p_c$, there exists almost surely no infinite open 
cluster.
\end {itemize}
The value of $p_c$ is lattice-dependent. In the case of the 
square lattice, it has been shown to be larger than 
$.556$ \cite {BE} (it is not expected to be any special 
number), while for    
the triangular lattice, it has been shown by 
Kesten and Wierman to 
be equal to $1/2$ (see e.g. \cite {Kesbook}).
This is not surprising because the triangular lattice
has a self-matching property: It is equivalent to 
say that the origin is in a finite open cluster or
to say that it is surrounded by a circuit (on the same
lattice, this is what makes the triangular
lattice so special) of closed sites. 
This property shows also that if $p=1/2$ on the triangular
lattice, the probability that there exists a left-to-right 
crossing of open sites of a square is exactly $1/2$ (otherwise,
there is a top-to-bottom crossing of closed sites). 
Russo, Seymour and Welsh \cite {Rus,SW}
have shown (this is sometimes known as the 
RSW theory) that this in fact implies that for any fixed 
$a$ and $b$, there exists a constant $c>0$, such that the 
probability $q(aN, bN)$
 of a left-to-right crossing of the $aN \times bN$ rectangle 
satisfies
$$ 1- c > q (aN, bN ) > c $$
for all large $N$.
This strongly suggests that when $N \to \infty$, $q(aN, bN)$
converges to a limit $F(b/a)$. 
A renormalizing group argument (loosely speaking, the rectangle
$2aN \times 2bN$ can be divided into four 
rectangles of size $aN \times bN$, which themselves can 
be divided into four rectangles etc.) also heuristically
suggests that not only the crossing probabilities 
converge but that in some sense, the information about 
``macroscopic connectivity properties'' should converge.
Note however that things are rather subtle.
Benjamini, Kalai and Schramm \cite {BKS} have for instance 
proved that 
if $A[N]$ denotes the event that there is a left-to-right crossing
of a $N \times N$ square say, and if one changes the status of a 
fixed proportion $\eps$ of the $N^2$ sites and looks at the event 
$\tilde A [N]$ that there exists a left-to-right crossing for
the new configuration, then the events $A[N]$ and $\tilde A [N]$ 
are asymptotically independent when $N \to \infty$. These events
are ``sensitive to noise''. 
When $N$ is large, it is not easy to ``see'' whether the 
crossing events occur or not
(in the Figure \ref {f.crossing}, each occupied 
site on the triangular lattice is represented by a white hexagon).

\begin{figure}
\centerline{\includegraphics*[height=3in]{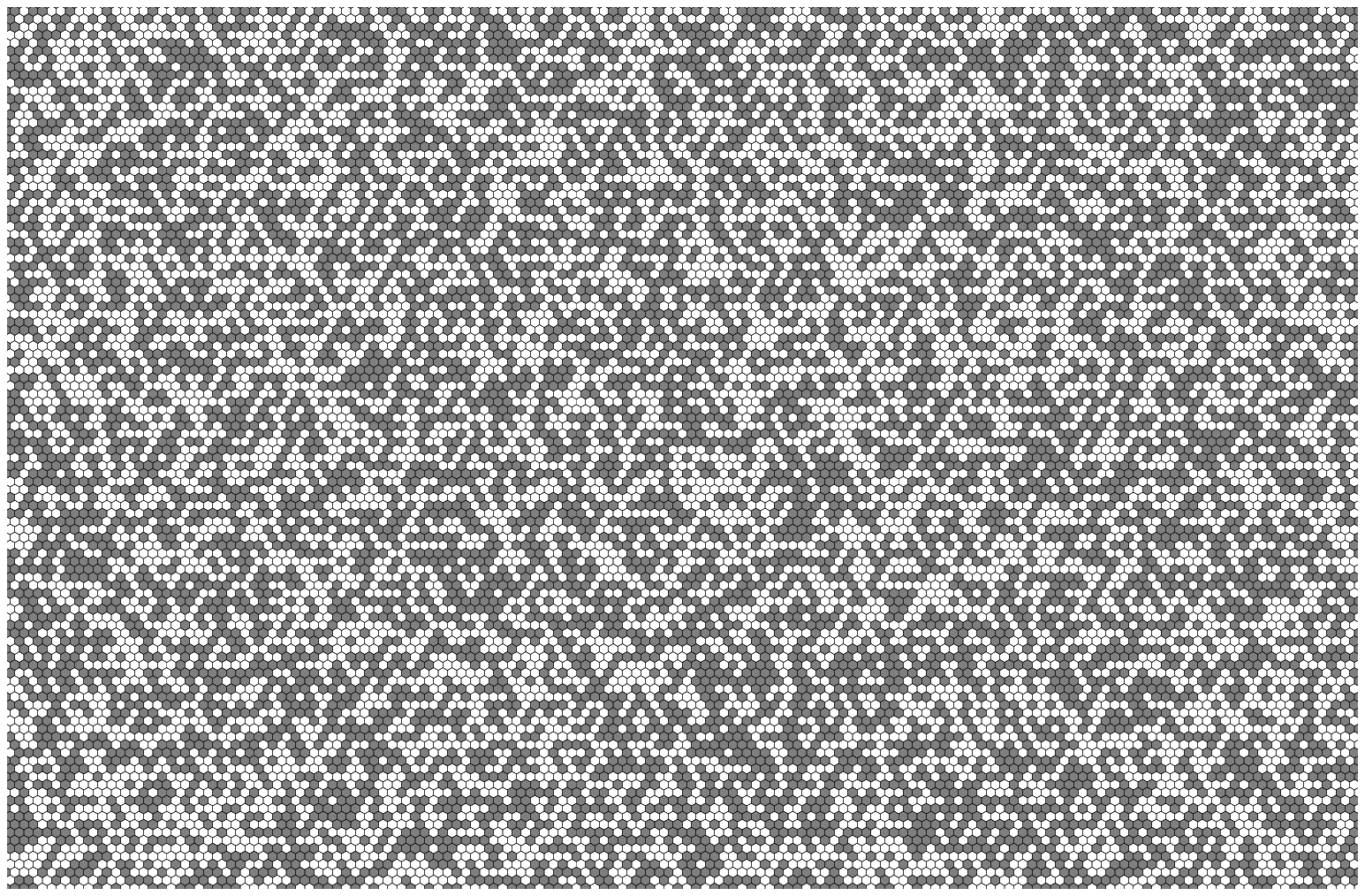}}
\caption{\label{f.crossing2}And now?}
\end{figure}

In fact, the renormalization argument suggests that even 
though the value of $p_c$ is lattice-dependent, on large scale,
what one sees at the value $p_c$ becomes lattice-independent. In 
other words, in the scaling limit, the behaviour of critical 
percolation should become lattice-independent (just as 
simple random walk converges to Brownian motion, for 
all ``regular'' lattices).  
Hence, the function $F(b/a)$ should be a universal
function describing the crossing-probabilities of a 
``continuous percolation process.''
In fact, this continuous percolation should be scale-invariant
(it is a scaling limit) as well as rotationally invariant (which 
would follow from lattice-independence). This leads to the 
stronger conjecture that it should be conformally invariant:
The connections in a domain $D$ and those in a domain $D'$ should have
the same law, modulo a conformal map from $D$ onto $D'$.

\section {The Cardy-Smirnov formula}

Using the conformal field theory 
ideas developed in \cite {BPZ,Ca1}, John Cardy \cite {Ca2}
 gave an exact prediction 
for the function $F$.
Extensive numerical work (e.g., \cite {LPS})
 did comfort these 
predictions.
Carleson noted that Cardy's function $F$ is 
closely related with the conformal maps 
from rectangles onto equilateral triangles, and that Cardy's
prediction could be rephrased as follows:

\begin {conjecture}[Cardy's formula]
If $D$ is conformally equivalent to the equilateral triangle $OAC$,
and if the four boundary points $a,o,c,x$  
are respectively mapped onto $A, O, C, X \in [CA]$, then (in the 
scaling limit when the mesh of the lattice goes to zero),
 the probability that there exists a 
crossing in $D$ from the part $(ao)$ of $\partial D$ to $(cx)$
is equal to 
$CX/CA$.
\end {conjecture}

We have seen that the SLE approach did 
provide a new justification to this formula.
Indeed, if the percolation exploration path
has a conformally invariant scaling limit,
it must be one of the chordal SLEs, as argued it the first Chapter.
Also, as the hitting probabilities computations in Chapter 3 show, 
$SLE_6$ is the unique SLE such that for all $X \in [CA]$, the two following 
probabilities are identical:
\begin {itemize}
\item
The SLE from $O$ to $A$ in the equilateral triangle hits $AX$ before $XC$
\item
The SLE from $O$ to $X$ in the equilateral triangle hits $AX$ before $XC$
\end {itemize}
This has to hold for the scaling limit of the critical exploration process. 
Hence, the unique possible conformally invariant scaling limit of the 
critical exploration process is $SLE_6$.
Another way to justify this is that this scaling limit has to satisfy locality and 
(cf. Chapter 4) that $SLE_6$ is the unique SLE that satisfies locality. Yet another (simpler) justification is that $SLE_6$ is the unique SLE for which the probability 
of the event corresponding to a left-right crossing of a square (or a rhombus) is $1/2$
(for an SLE starting from one corner and aiming at a neighbouring corner).

We have also seen that an $SLE_6$ from $O$ to $C$ in the equilateral 
triangle hits $XC$ before $AX$ with probability $CX/CA$.
Also, for the discrete exploration process, the corresponding event is 
precisely the event that there exists a crossing from $AO$ to $CX$.
Hence, we get a conditional result of the 
following type: If the scaling limit of critical 
percolation exists and is conformally invariant, then the scaling limit of 
the exploration process is $SLE_6$ and Cardy's formula holds. 
 
But in order to prove conformal invariance of critical percolation, one has to work 
with discrete percolation itself. In 2001, Stas Smirnov, proved that:

\begin {theorem}
Cardy's prediction is true in the case of critical site percolation
on the triangular lattice.
\end {theorem}

In fact, Smirnov's proof is a direct proof of Cardy's formula that does not 
rely at all on SLE. Then, with Smirnov's result, one can show that indeed 
the scaling limit of the percolation exploration process is $SLE_6$.

\medbreak
\noindent
{\bf Sketch of the proof.}
Suppose first for convenience that $AOC$ is an equilateral triangle 
and that the sides of the triangle have unit length
and are parallel to the
axis of the triangular grid (as we will see, this has in fact no other
influence on the proof than simplifying the notations).
For all $\delta = 1/n$, consider critical site percolation in 
$AOC$ on the triangular grid with mesh-size $1/n$.
For convenience, put $\tau= \exp ( 2 i \pi / 3)$  
and write $A_{1}= A$, $A_\tau = A_2= O$ and $A_{\tau^2} = A_3= C$.
For each face $z$ of the triangular grid  (i.e. for each site of
the dual hexagonal lattice), let $E_1(z)$ denote the event
that there exists a simple open (i.e. white) path from $A_1A_\tau$
to $A_1A_{\tau^2}$ that separates $z$ from $A_{\tau}A_{\tau^2}$.
Similarly, define the events $E_\tau (z)$ and $E_{\tau_2}(z)$ corresponding
to the existence of simple open paths separating $z$ from 
$A_1A_{\tau^2}$ and $A_1A_\tau$ respectively.
Define finally for $j = 1, \tau , \tau^2$,
$$
H_j( z) =  H_j^\delta (z) := \P [ E_j (z) ].
$$
The Russo-Seymour-Welsh theory ensures that the functions $H_j^\delta$
are uniformly  ``H\"older'' (actually, one 
first has to smooth out their 
discontinuities for instance in a linear way keeping only the values of 
$H_j^\delta$ at the center of the triangles). 
In particular, it shows that any for any sequence $\delta_n \to 0$,
the triplet of functions $(H_1^\delta, H_\tau^\delta, H_{\tau^2}^\delta)$
has a subsequential limit.
Our goal is now to identify the only possible 
such subsequential limit.

The Russo-Seymour-Welsh estimates also show that 
when $z \to A_{j \tau}A_{j \tau^2}$, the functions $H_j^\delta$
go uniformly to zero, and that when $z \to A_j$, the functions $H_j^\delta$
go uniformly to one.
Hence, for any subsequential limit $(H_1, H_\tau, H_{\tau^2})$, one has
 $H_j (z) \to  0 $ when  $z \to A_{j\tau} A_{j \tau^2}$, and 
 $ H_j (z) \to 1$ when $z \to A_j$.

Now comes the key-observation of combinatorial nature:
Suppose that $z$
is the center of a triangular face. Let $z_1, z_2,z_3$
denote the three (centers of the) neighbouring faces
 (with the same orientation as the triangle $A_1A_2A_3$)
and $s_1,s_2, s_3$ the three corners of the face containing
$z$ chosen in such a way that $s_j$ is the corner ``opposite''
to $z_j$.
We focus on the event
$E_1 (z_1) \setminus E_1 (z)$.
This is the event that there exists three disjoint 
paths $l_1$, $l_2$, $l_3$ such that 
\begin {itemize}
\item
The two paths $l_2$ and $l_3$ are open and join
the two sites $s_2$ and $s_3$ to $A_1A_3$ and $A_1A_2$
respectively.
\item
The path $l_1$ is closed (i.e., it consists only of closed sites), and joins $s_1$ to $A_2 A_3$.
\end {itemize}
One way to check whether this event holds is to 
start an exploration process from the corner $A_3$, say (leaving the 
open sites on the side of $A_1$ and the closed sites on 
the side of $A_2$).
If the event 
$E_1(z_1) \setminus E_1(z)$ is true, then the exploration process 
has to go through the face $z$, arriving into 
$z$ through the 
edge dual to $s_1s_2$. In this way, one has 
``discovered'' the simple paths $l_2$ and $l_1$ that 
are ``closest'' to $A_3$.
Then, in the remaining (unexplored domain), there must 
exist a simple open path from $s_3$ 
to $A_1 A_3$. 
But, the conditional probability of this event is the 
same as that of the existence of a simple closed 
path from $s_3$ to 
$A_1 A_3$ (interchanging open and closed in the 
unexplored domain does not change the probability 
measure). Changing all the colors once 
again, shows finally  that $E_1 (z_1) \setminus E_1 (z)$
has the same probability as the event that there exist
three disjoint 
paths $l_1$, $l_2$, $l_3$ such that 
\begin {itemize}
\item
The paths $l_1$ and $l_3$ are open and join
the two sites $s_1, s_3$ to $A_2A_3$ and $A_1A_2$
respectively.
\item
The path $l_2$ is closed, and joins $s_2$ to $A_1 A_3$.
\end {itemize}
This event is exactly 
$E_\tau ( z_2) \setminus E_\tau (z)$.
Hence, we get that,
$$
\P [ E_1 (z_1) \setminus E_1 (z)] 
=
\P [ E_\tau (z_2) \setminus E_\tau (z)]
= 
\P [ E_{\tau^2} (z_3) \setminus E_{\tau^2} (z)]
.$$
These identities can then be used to  show that 
for any equilateral contour $\Gamma$ (inside the 
equilateral triangle), 
the contour integrals of $H_j^\delta$ for $j=1,\tau, \tau^2$
are very closely related:
$$
\int_\Gamma dz H_1^\delta (z) =
\int_\Gamma dz H_\tau^\delta (z) / \tau + O(\delta^\eps)
= \int_\Gamma dz H_{\tau^2}^\delta (z) / \tau^2 
+ O (\delta^\eps)
$$ 
when $\delta \to 0$ for some $\eps >0$.
To see this, one has to 
expand the 
contour integrals as the sum of all
properly oriented contour integrals along all small 
triangles inside $\Gamma$. Then, the previous identities ensure that 
almost all terms cancel out. The remaining ``boundary''
terms are controlled with the
help of RSW estimates.

This result
then shows that for any subsequential limit $(H_1,H_\tau, H_{\tau^2})$,
the contour integrals of 
$H_1$, $H_\tau/\tau$ and of $H_{\tau^2}/\tau^2$ coincide.
It readily follows that the contour integrals of the functions 
$$
H_j + \frac {i}{\sqrt {3}} ( H_{j\tau} - H_{j \tau^2})
$$
for $j=1, \tau, \tau^2$
vanish. By Morera's theorem (see e.g. \cite {A1}), this ensures that these
functions are analytic. In particular, $H_1$ is harmonic.
The boundary conditions $H_j = 0$ on $A_{j \tau} A_{j \tau^2}$ for 
$j=1, \tau , \tau^2$ then ensure that 
$H_1 = 0$ on $A_2A_3$ and that the horizontal derivative of
$H_1$ on $A_1A_3 \cup A_2A_3$ vanishes.
Also, $H_1 (A_1) = 1$.
The only harmonic function 
in the equilateral triangle with these boundary conditions is 
the height 
$$H_1 (z) = \frac { d(z, BC) }{d(A, BC)}.
$$
This completes the proof of the Theorem when the domain is an equilateral triangle.

If $D$ now any simply connected domain, and $a=a_1$, $o=a_{\tau}$, $c=a_{\tau^2}$ are boundary points, the proof is almost identical. In its first part, the only difference is 
that one replaces the straight boundaries $A_{j} A_{j\tau}$
by approximations of the boundary of $D$ on the triangular lattice that is between the 
points $a_j a_{j \tau}$.
In exactly the same way, one obtains tightness and boundary estimates for the 
discrete functions $H_j^\delta$. Also, the argument leading to the fact that 
the contour integrals on equilateral triangles of 
$H_j + i ( H_{j \tau} - H_{j \tau^2})/\sqrt {3}$ for any subsequential limit 
vanish, remains unchanged. Hence, for any subsequential limit, one obtains a 
triplet of functions $(H_1, H_\tau, H_{\tau^2})$ such that for $j=1, \tau, \tau^2$:
\begin {itemize}
\item
The function $H_j + i ( H_{j \tau} - H_{j \tau^2})/ \sqrt {3}$ is analytic
\item
The function $H_j (x)$ tends to zero  when $x$ approaches the part of the 
boundary between $a_{j\tau}$ and $a_{j\tau^2}$. 
\item
The function $H_j (x)$ tends to one when $x \to a_j$.
\end {itemize}
The important feature is that this problem is conformally invariant: 
If $\Phi$ denotes a conformal map from $D$ onto the equilateral triangle such 
that $\Phi (a_j) = A_j$, and if $(H_1, H_\tau, H_{\tau^2})$ is such a triplet of 
functions, then the triplet $(H_1 \circ \Phi^{-1}, H_\tau \circ \Phi^{-1},
H_{\tau^2} \circ \Phi^{-1})$ solves the same problem in the equilateral
triangle. In the latter case, we have seen that the unique solution is given by 
$ H_j (x) = d( x, A_{j \tau} A_{j \tau^2}) / d (A_j, A_{j \tau} A_{j \tau^2})$.
Hence, the Theorem follows.
\qed

\medbreak

One should stress that this  proves much more than just the
asymptotic behaviour of the crossing probabilities.
It yields the asymptotic probability of the events $E_j (x)$ for $x$ inside the domain $D$ (and not only on its boundary).

\section {Convergence to $SLE_6$ and consequences}

One can use the previous result to prove that the discrete exploration 
process described in the introductory chapter indeed converges 
to chordal $SLE_6$. 

The regularity estimates are provided by the RSW theory and the 
discrete Markovian property is immediate. It remains to 
show that some macroscopic quantities converge to a conformally
invariant 
quantity in the scaling limit, but this is precisely what Smirnov's 
theorem shows. Hence, the method described in the previous 
chapter can be applied. Some adjustments are needed to 
take care of domains with rough boundary, though. In particular,
one can use the a priori bounds on the probability of having $5$
arms joining the vicinity of the origin to a large circle 
(the exponent $\alpha_5$ below) derived in \cite {KSZ}.

Exploiting this, one can therefore use the computations of critical 
exponents for $SLE_6$, to deduce asymptotic probabilities for 
discrete critical percolation on the triangular lattice:
For instance \cite {SmW, LSW5}, 
let $A_n [N]$ denote the event that there exists $n$ disjoint
open clusters  
joining the vicinity of the origin to the circle of radius $N$.
Then:
\begin {theorem}
When $N \to \infty$, 
one has 
$\P [ A_n [N]  ] \approx N^{- \alpha_n}$, 
where 
$\alpha_1 = 5/48$ and for all $n \ge 2$,
$\alpha_n = (4 n^2 - 1) /12$.
\end {theorem}
Note that the exponents $\alpha_n$ for $n \ge 2$ are the 
same than the Brownian intersection exponents
$\xi_n$ in Chapter 8. This  is 
not surprising because of the close relation between
$SLE_6$ and planar Brownian motion.
The exponent $\alpha_1$ corresponds to the event that radial 
$SLE_6$ winds only ``in one direction''  
around $0$ (see \cite {LSW5}.

Actually, Harry Kesten \cite {Kes} had shown that the previous 
result (for $n=1$ and $n=2$) would imply the following
description of the behaviour of percolation when the 
probability is near to the critical probability:

\begin {theorem}
If one performs site percolation on the triangular lattice with
probability $p$, then when $p \to 1/2+$,
the probability that the origin belongs to the infinite
cluster behaves like $(p-1/2)^{5/36 + o(1)}$.
When $p \to 1/2 - $, the correlation length
explodes like $(1/2 - p)^{-4/3+ o(1)}$.
\end {theorem}

See \cite {Kes, SmW} for more results as well 
as for the proofs...

\begin{figure}
\centerline{\includegraphics*[height=3in]{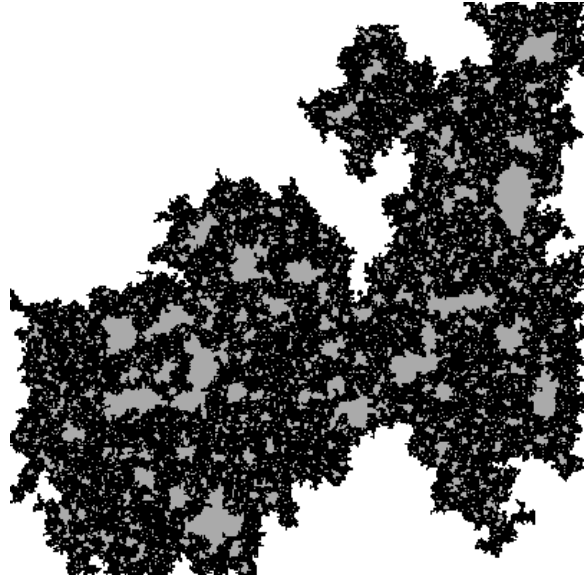}}
\caption{\label{f.cluster} 
Part of a (big) critical percolation cluster on the square lattice}
\end{figure}

\begin{figure}
\centerline{\includegraphics*[height=3in]{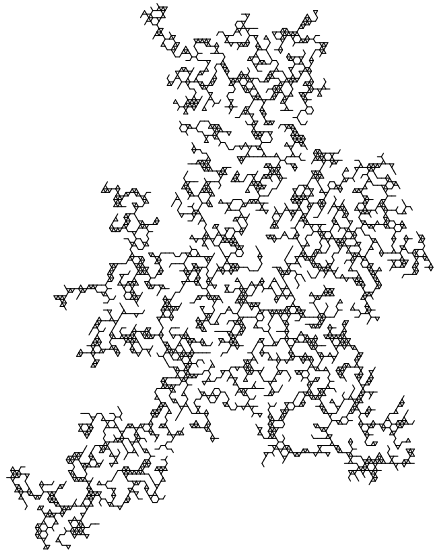}}
\caption{\label{f.cluster2} 
A critical percolation cluster on the triangular lattice}
\end{figure}

Let us conclude with the following combination of results
that we have mentioned in these lectures:
The following three curves are (locally) the same:
\begin {itemize}
\item
The outer boundary of the scaling limit of a large critical 
percolation cluster.
\item
The outer boundary of a planar Brownian motion.
\item
The scaling limit of long self-avoiding walks, provided this 
scaling limit exists and is conformally invariant.
\end {itemize}

\bibc
The value of critical exponents for percolation had been 
predicted by theoretical  
physicists \cite {DN,Pe,N,N2,SRG,SD,GA,Ca3,ADA}.
The conformal invariance conjecture for critical percolation  
had been discussed by Aizenman \cite {Ai,Ai2}.

Smirnov's complete detailed proof of 
Cardy's formula is contained in 
\cite {Sm1, Sm2}.
The actual detailed proof of the convergence of the discrete exploration 
process to $SLE_6$ (announced in \cite {Sm1}) should be 
written up \cite {Sm2} soon.
For the derivation of formulas and exponents for critical
percolation using $SLE_6$, see \cite {S2, LSW5, SmW}.

\chapter {What is missing}

\section {A list of ideas}

We have listed at the end of each chapter a list of references to papers that
develop ideas that are related to those presented in the corresponding chapter.
One aspect of SLE that we could have spent more time on is the actual computation of 
critical exponents. For simplicity, we have shown how to derive the Brownian exponents 
using radial $SLE_6$, but in general (for instance to derive the Hausdorff dimension 
of the SLE), one might as well work with chordal SLE.  Various exponents are derived in 
for instance in \cite {LSW1,LSW2,LSW3,LSW5,RS,Be2,Be3}. 

Before very briefly reviewing the results related to 
restriction properties, we would like to stress that the important
ideas underlying Rohde-Schramm's \cite {RS} proof of the
 existence and transience of the SLE
paths have not been presented in these lectures. The arguments \cite {RS}
require some non-trivial background in complex geometry. 
In two cases, the existence and/or transience of the SLE path is especially 
difficult to establish: For $\kappa=4$, because the domains generated by the SLE curve 
are not H\"older (see \cite {RS}). 
For $\kappa=8$, the only proof uses the fact that it is 
the scaling limit of the 
discrete uniform  Peano curves \cite {LSWlesl} described in Chapter 9.

One can also study geometric questions such as: Does the SLE have (local) cut points?
The answer is positive if and only if $\kappa < 8$ (see \cite {Be2}).

I plan to discuss
the following restriction properties 
in forthcoming lecture notes. The main reference is the
long recent paper \cite {LSWrest}. 
\begin {itemize}
\item
The full classification of the measures satisfying the restriction properties is one of the main goals of \cite {LSWrest}. 
These measures form a one-dimensional family indexed by a 
positive real-valued parameter $N$, that can be interpreted as the number of Brownian 
excursions that the measure is equivalent to. 
There exist two other important ways to
describe this one-dimensional family: The first one is via a variant of the 
$SLE_{8/3}$ process called $SLE(8/3, \rho)$. Loosely speaking, one replaces the driving 
Brownian motion by a Bessel process (see \cite {LSWrest} for all this), and the obtained 
simple random curve describes the outer boundary of the set satisfying the restriction 
property. 
The second description goes as follows: Consider an $SLE_\kappa$ with 
$\kappa \le 8/3$ and  add to this path a certain cloud Brownian loops
(this Poisson cloud of loops is also studied in \cite {LWbls}).
For a well-tuned density $d(\kappa)$ of the loops, one constructs the 
restriction measure corresponding to $N(\kappa)$ Brownian excursions. 
See also \cite {Dub2}. 
\item
This last description makes it possible to tie a link \cite {FW1,FW2}
with representation theory, and more precisely with highest-weight representations 
of the Lie Algebra of polynomial vector fields on the unit circle (the number $N(\kappa)$
is the highest-weight). This is related 
to considerations from conformal field theory. See also \cite {BB1,BB2,BB3} for 
the relation of SLE with ideas from conformal field theory. 
\item
The $SLE (\kappa, \rho)$ processes shed also some light on the computation of the 
(chordal) critical exponents. It turns out that they can be understood via the 
absolute continuity relations between Bessel processes (following from
Girsanov's Theorem); see \cite {Whid}. 
\end {itemize}
 
\section {A list of open problems}

Here is a list of open problems. Some of these were already mentioned in the previous 
chapters:

\subsection {Conformal invariance of discrete models}

So far, convergence of natural discrete models
towards SLE in the scaling limit has been proved only in the
two very special cases that we described in the last two chapters
(LERW-UST, and critical site percolation on the triangular lattice).
It is believed to hold for many other models:

\begin {itemize}
\item
The interface for a critical FK-percolation (see e.g. \cite {GLN} for 
an introduction to this dependent percolation model introduced by 
Fortuin and Kasteleyn) model for $q \le 4$ is conjectured to 
converge to chordal $SLE_\kappa$. Recall that the probability of a given 
realization is proportional to 
$$ p^{\# \hbox {open edges}} (1-p)^{\# \hbox {closed edges}}
q^{\# \hbox {connected components}}.$$
The relation between $q$ and $\kappa$ should be 
$$ \cos \frac {4\pi}\kappa = - \frac {\sqrt q }2,$$
where $q \in [0,4]$ and $\kappa \in [4,8]$. 
Here (as in the UST case and in some sense in the percolation case), the boundary conditions
have to be mixed (free on one part of the boundary, wired on the other -- this influences the way of counting the connected components).
See \cite {RS} for a more precise statement of this conjecture. Recall that for 
critical FK-percolation with parameter $q$ on the square lattice, the self-dual 
point is $p = \sqrt {q} / ( \sqrt {q} + 1)$ (proving that this self-dual point is the 
critical point is another open question, but it is not directly related to the 
SLE question; the question on the square grid is to prove that for this value of 
$p$, the interface converges to SLE).
Here self-dual means that the law of the dual graph of an 
FK percolation sample is also an FK percolation sample (in the dual lattice)
with the same parameters (see \cite {GLN}). 

Recall that when $q >4$, the FK percolation phase transition is conjectured to 
be a first-order transition (i.e. there can exist an infinite open cluster at the critical 
probability). The critical value $q=4$ corresponds to the special case $\kappa=4$.
Recall also (see e.g. \cite {GLN}) that the correlation functions of the critical 
$q$-Potts models are the same as those of the critical FK-percolation model.
Recall also that the usual percolation is the $q=1$ FK percolation model, and that 
the UST can be viewed as the $q=0$ critical FK percolation model (see e.g. \cite {H}).
For the critical FK percolation models, the Markovian property is clearly valid in the 
discrete case. The missing step is therefore the proof of conformal invariance.

It is interesting (and encouraging) to note that the integer values of $q$ correspond to the ``nice'' values of the angles $\alpha = \pi ( 1 - \kappa /4)$ of the isocele triangles for which 
hitting distributions are uniform
(Dub\'edat's observations \cite {Dub} mentioned at the end of Chapter 2):
 $\cos \alpha = \sqrt {q} / 2$. For $q=1$, it is the 
equilateral triangle, for $q=2$ (Ising), it is the isocele-rectangular triangle, and
for $q=3$, $\alpha = \pi /6$.
 
\item
Among all the critical percolation interfaces that are conjectured to converge to 
$SLE_6$ (this is the special case $q=1$ in the previous conjectures), it is 
worth stressing two cases, for which one has self-duality (and therefore some 
little hope to be able to prove something): The first one is 
bond-percolation on the square grid, and the second one is percolation on a 
Voronoi tessellation (see e.g. \cite {BS}).

\item
There exists a special model for which (as for the Ising model and for the 
uniform spanning tree model), the tools and arguments developed by Kenyon seem promising:
It is the so-called double-domino path, that is conjectured to converge to the 
special curve $SLE_4$ in the scaling limit.

\item
Note also that the Ising model itself (on the triangular lattice) has some self-duality
properties (this is due to the fact that for the Ising model, there are exactly two 
possible states for each site). Hence, Ising cluster interfaces (for appropriate boundary 
conditions, and on the triangular lattice) might converge to an SLE in the scaling limit. 

\item
For $\kappa < 4$, the relation with discrete models from statistical physics is not   
so clear. One relation is via the duality conjectures that we will discuss below.
The main open question is the convergence of the self-avoiding walk towards
the $SLE_{8/3}$ curve. Again, the main problem is to derive its conformal 
invariance. See \cite {LSWSAW} for a discussion. Let us insist that 
basically nothing is known rigorously on the asymptotic behaviour of the self-avoiding 
walk. For instance, to our knowledge, it has not even been disproved that the curve becomes
space-filling or a straight line in the scaling limit!

\item
It is likely that some discrete dynamic models can be shown to converge to SLE
(but their relation to models from statistical physics is unclear). For instance,
variations on the Laplacian random walk description of LERW that have some
conformally invariant features built in the model should in principle converge
to SLEs.
\end {itemize}

\subsection {Duality}
Another approach to the SLE curves when $\kappa < 4$ goes as follows:
It was conjectured (based on the computation of the dimensions) that in 
the scaling limit, the outer boundary of an $SLE_{\kappa'}$ hull for $\kappa' >4$
 at a given time looks (locally) like 
an $SLE_{16/\kappa'}$ curve.
Hence, the $SLE_\kappa$ curves for $\kappa<4$ correspond to the outer boundary of the 
scaling limit of critical FK-percolation clusters. 
The duality has been proved to hold in two cases: $\kappa=2$ 
(because of the relation between
LERW and UST that respectively converge to $SLE_2$ and $SLE_8$) and 
$\kappa=8/3$ 
(because of the restriction property considerations that allow to describe 
the outer boundary of conditioned $SLE_6$ processes in terms of $SLE_{8/3}$ 
processes (see \cite {LSWrest}).
In the general case, a weak form of duality has been identified by Dub\'edat \cite {Dub2},
that leads to conjecture a precise identity in law between the outer boundary
of an $SLE(\kappa', \rho_{\kappa'})$ process and the $SLE (16/\kappa', \rho_{\kappa'}')$
curve for well-chosen values of $\rho$ and $\rho'$. 

Proving this duality relation would be one way to settle the following open 
problem (it is only proved when $\kappa=6$ and $\kappa=8$):
Prove that the Hausdorff dimension of the boundary of $K_t$ is almost surely 
$1+2 / \kappa$ when $K_t$ is the hull of an $SLE_\kappa$ (chordal or radial)
for $\kappa >4$.
One would then combine duality with the computation of the dimension of the SLE
curves in \cite {Be3}.
There should however also exist a direct proof of this fact that does not 
rely on duality.  

\subsection {Reversibility}
The following conjecture follows very naturally from the fact that the SLEs are 
believed to be scaling limit of the previously described lattice models:
Suppose that 
$\kappa \le 8$ is given, and consider the chordal $SLE_\kappa$ curve $\gamma$ from 
$a$ to $b$ in a domain $D$ (where $a$ and $b$ are two boundary points).
One can time-reverse $\gamma$, and view it as a curve from $b$ to $a$ in $D$.
Then, the law of this time-reversal should be (modulo time-change) the law of an
$SLE_\kappa$ curve from $b$ to $a$ in $D$.
Another equivalent way of phrasing this is that if $\gamma$ is the chordal SLE path 
in the upper half-plane, the path $-1/ \gamma$ has the same law as $\gamma$ (modulo 
time-change).

This conjecture is very natural in terms of the lattice 
models, but on the other hand, it is not natural at all if one thinks of 
the actual definition of the SLE in terms of the Loewner chain (this is very 
non-reversible!).
In the special cases $\kappa=6$, $\kappa=8$ and $\kappa=2$, the 
result is a consequence of the convergence of the discrete reversible models to the SLEs.
So far, the reversibility of $\kappa=8/3$ is the only one that can be proved without reference to a reversible discrete model, and the tool here is the characterization of 
$SLE_{8/3}$ as the unique simple random curve that satisfies the restriction property.
In all other cases, the problem is to our knowledge open. This problem does not seem 
as out of reach as some of those that we just discussed.

Note that (as shown to me by Oded Schramm),
it is possible to show that reversibility of $SLE_\kappa$ fails to be true 
when $\kappa > 8$. This can seem surprising; more generally, the interpretation of 
$SLE_\kappa$ when $\kappa >8$ in terms of models from  statistical physics is 
not well-understood. Note that the asymptotic behaviour of $SLE_\kappa$ when $\kappa 
\to \infty$ is studied in \cite {Be4}.

\subsection {Quantum gravity and conformal field theory}

The arguments developed in conformal field theory under the name of 
quantum gravity suggest that some very interesting critical phenomena also occur for 
systems on certain random lattices. 
In particular, Duplantier \cite {Dqg,Dcif,D2003} showed that the value of the 
critical exponents in the plane (those exponents that can now be understood thanks to the 
SLE) can be predicted using the formula proposed by Knizhnik, Polyakov and Zamolodchikov
in \cite {KPZ}, that should relate the value of the critical exponents in the 
plane to the corresponding exponents on random lattices.

Recent progress has been made in the 
rigorous understanding of some of these random systems 
on these random graphs; see e.g. \cite {AnS,An,BMS,BDE} and the references 
therein. It seems that (as opposed to the rigid lattice case), the behaviour 
of some of these systems on random lattices might be
 accessible by ingenious combinatorial methods. 

Note \cite {Whid} that the KPZ formula seems to have a simple interpretation in terms 
of the $\rho$ in the $SLE(\kappa, \rho)$ processes. Maybe the combination of the 
determination of the exponents for SLE, and the results on random graphs will
provide in the end the rigorous justification to the KPZ relation.

More generally, the relation between SLE and conformal field theory 
(that has started to be investigated in \cite {BB1,BB2,BB3,FK,FW1,FW2}) 
and with the mathematical concepts used in conformal field theory
needs further understanding. It is not so clear whether this will be 
helpful to improve the knowledge on these critical two-dimensional systems 
(which was after all probably the initial motivation for the conformal field 
framework).
One related issue is
to manage to define SLE on general Riemann surfaces, see \cite {FK,Z,Dub4}.

\vfill
\eject



\begin {thebibliography}{999}
\bibitem {A1}
{L.V. Ahlfors,
{\em Complex analysis}, 3rd Ed., McGraw-Hill, New-York, 1978.}

\bibitem {A2}{
L.V. Ahlfors,
{\em Conformal Invariants, Topics in Geometric Function
Theory}, McGraw-Hill,  New-York, 1973.}

\bibitem {Ai}{
M. Aizenman (1996),
The geometry of critical percolation and conformal
invariance, Statphys19 (Xiamen, 1995), 104-120.
}

\bibitem {Ai2}
{M. Aizenman (1998),
Scaling limit for the incipient spanning clusters, in 
{\sl Mathematics of multiscale materials},
IMA Vol. Math. Appl. {\bf 99}, 
Springer, New York, 1-24.}
 
\bibitem {AB}
{M. Aizenman, A. Burchard (1999),
H\"older regularity and dimension bounds for random curves,
 Duke Math. J. {\bf 99}, 419--453. }

\bibitem {ABNW}
{M. Aizenman, A. Burchard, C.M. Newman, D.B. Wilson
(1999),
Scaling limits for minimal and random spanning trees in two dimensions. 
Random Structures Algorithms {\bf 15}, 319-367.}

\bibitem {ADA} { 
M. Aizenman, B. Duplantier, A. Aharony (1999),
Path crossing exponents and the external perimeter in 2D percolation.
Phys. Rev. Let. {\bf 83}, 1359-1362.}

\bibitem {An}
{O. Angel (2002),
Growth and Percolation on the Uniform Infinite Planar Triangulation,
preprint.}

\bibitem {AnS}
{O. Angel, O. Schramm (2002),
Uniform Infinite Planar Triangulations,
preprint.
}

\bibitem {BB1}
{M. Bauer, D. Bernard (2002),
$SLE_k$ growth processes and conformal field theories
 Phys. Lett. {\bf B543}, 135-138.}

\bibitem {BB2}
{M. Bauer, D. Bernard (2002),
Conformal Field Theories of Stochastic Loewner Evolutions,
preprint.}

\bibitem {BB3}
{M. Bauer, D. Bernard (2003),
SLE martingales and the Virasoro algebra, preprint.}


\bibitem {Be1}{V. Beffara (2001),
On some conformally invariant subsets of the planar Brownian curve,
Ann. Inst. Henri Poincar\'e, to appear}

\bibitem {Be2}
{V. Beffara (2002), Hausdorff dimensions for $SLE_6$, preprint.}

\bibitem {Be3}
{V. Beffara (2002), The dimension of the SLE curves, preprint.}

\bibitem {Be4}
{V. Beffara, in preparation}

\bibitem {BPZ0}
{A.A. Belavin, A.M. Polyakov, A.B. Zamolodchikov (1984),
 Infinite conformal symmetry of critical fluctuations in two dimensions, 
J. Statist. Phys. {\bf 34}, 763--774.}

\bibitem{BPZ}
{A.A. Belavin, A.M. Polyakov, A.B. Zamolodchikov (1984),
Infinite conformal symmetry in two-dimensional quantum field theory.
Nuclear Phys. B {\bf 241}, 333--380.}

\bibitem{BKS}
{I. Benjamini, G. Kalai, O. Schramm (1999), 
       Noise sensitivity of boolean functions and applications to percolation,
       Publ. Sci. IHES {\bf 90}, 5-43.
}

\bibitem {BLPS}
{I. Benjamini,   
R. Lyons, Y. Peres, O. Schramm (2001), 
       Uniform spanning forests.
       Ann. Probab. {\bf 29}, 1-65.
}

\bibitem {BS}
{I. Benjamini, O. Schramm  (1998),
Conformal invariance of Voronoi percolation, 
Comm. Math. Phys. {\bf 197}, 75-107.}

\bibitem {BE}
{J. van den Berg, A.  Ermakov (1996), 
A new lower bound for the critical probability 
of site percolation on the square lattice,
Random Structures Algorithms {\bf 8},199-212. }

\bibitem {BJ}
{R. van den Berg, A. Jarai (2001),
The lowest crossing in 2D critical percolation, preprint}

\bibitem {BJPP}
{C.J. Bishop, P.W. Jones, R. Pemantle, Y. Peres
(1997),
The dimension of the Brownian frontier is greater than $1$,  
J. Funct. Anal. {\bf 143}, 309--336. }

\bibitem {BMS}
{M. Bousquet-M\'elou, G. Schaeffer (2002),
The degree distribution in bipartite planar maps: applications to the Ising model,
preprint.}

\bibitem {BDE}
{J. Bouttier, B. Eynard, Ph. Di Francesco (2002),
Combinatorics of Hard Particles on Planar Graphs,
preprint.}

\bibitem {BL1}{
 K. Burdzy, G.F. Lawler (1990),
Non-intersection exponents for random walk and Brownian motion.
 I: Existence
 and an invariance principle, Probab. Theor. Rel. Fields {\bf 84},
393--410.}

\bibitem {BL2}
{ K. Burdzy, G.F. Lawler (1990),
Non-intersection exponents for random walk and Brownian motion.
 II: Estimates and applications to a random fractal,
Ann. Prob. {\bf 18}, 981-1009.
}

\bibitem {BP}
{R. Burton, R.  Pemantle 
(1993),
Local characteristics, entropy and limit theorems
 for spanning trees and domino tilings via transfer-impedances,
Ann. Probab. {\bf 21},  1329--1371. 
}

\bibitem {Ca1}
{J.L. Cardy (1984),
Conformal invariance and surface critical behavior,
Nucl. Phys. B {\bf 240} (FS12), 514--532.}

\bibitem{Ca2}
{J.L. Cardy (1992),
Critical percolation in finite geometries,
J. Phys. A, {\bf 25} L201--L206.}

\bibitem {Cabook}
{J.L. Cardy, 
{\em  Scaling and renormalization in statistical physics},
 Cambridge Lecture Notes in Physics {\bf 5},
 Cambridge University Press, 1996.
}

 \bibitem {Ca3}
{J.L. Cardy (1998),
The number of incipient spanning clusters in two-dimensional
percolation, J. Phys. A {\bf 31}, L105.}

\bibitem {CaChuo}
{J.L. Cardy (2001),
Lectures on Conformal Invariance and Percolation,
 Lectures delivered at Chuo University, Tokyo, preprint.}
 
\bibitem {CM1}
{L. Carleson, N.G. Makarov (2001),
Aggregation in the plane and Loewner's equation,
 Comm. Math. Phys. {\bf 216}, 583-607.}

\bibitem {CM2}
{L. Carleson, N.G. Makarov (2002),
Laplacian path models, preprint}

 \bibitem{CM}{
M. Cranston, T. Mountford (1991),
An extension of a result by Burdzy and Lawler,
Probab. Th. Relat. Fields {\bf 89}, 487--502.
}

\bibitem {DN} {M.P.M. Den Nijs (1979),
A relation between the temperature exponents of the eight-vertex 
and the $q$-state Potts model, J. Phys. A {\bf 12}, 1857-1868.
}

\bibitem {Dub} {J. Dub\'edat (2003), SLE and triangles,
El. Comm. Probab. {\bf 8}, 28-42.}

\bibitem {Dub2}
{J. Dub\'edat (2003), $SLE(\kappa, \rho)$ martingales and duality, 
preprint.}

\bibitem {Dub3}
{J. Dub\'edat (2003), 
Reflected planar Brownian motion, intertwining relations and crossing
probabilities, preprint.}

\bibitem {Dub4}
{J. Dub\'edat (2003),
preprint.}

\bibitem{Dle} 
{B. Duplantier (1992),
Loop-erased self-avoiding
walks in two dimensions: exact critical exponents and winding numbers,
Physica A {\bf 191}, 516--522.}

\bibitem {Dqg}
{B. Duplantier (1998),
Random walks and quantum gravity in two dimensions, Phys. Rev. Lett. {\bf
81},
5489--5492.}

\bibitem {Dp} 
{B. Duplantier (1999),
Harmonic measure exponents for two-dimensional percolation,
Phys. Rev. Lett. {\bf 82}, 3940-3943.}

\bibitem {Dcif}
{B. Duplantier (2000), 
 Conformally invariant fractals and potential theory, Phys. Rev. Lett. 
{\bf 84}, 1363-1367. 
}

\bibitem {D2003}
{B. Duplantier (2003),
Conformal Fractal Geometry and Boundary Quantum Gravity,
preprint}




\bibitem {DK}{
B. Duplantier, K.-H. Kwon (1988),
Conformal invariance and intersection of random walks, Phys. Rev. Let. {\bf
61},
 2514--2517.
}

\bibitem {Dur}{
P.L. Duren,
{\em Univalent functions}, Springer, 1983.}

\bibitem {Fom}
{S. Fomin (2001),
Loop-erased walks and total positivity, Trans. Amer. Math. Soc. 
{\bf 353}, 3563--3583.}


\bibitem {FK}
{R. Friedrich, J. Kalkkinen (2003),
preprint.}

\bibitem {FW1}
{R. Friedrich, W. Werner (2002),
Conformal fields, restriction properties, degenerate representations and SLE,
C.R. Ac. Sci. Paris Ser. I Math {\bf 335}, 947-952. 
}
 
\bibitem {FW2}
{R. Friedrich, W. Werner (2003),
Conformal restriction, highest-weight representations and SLE, preprint.}

\bibitem{GA}
{T. Grossman, A. Aharony (1987),
Accessible external perimeters of percolation clusters,
J.Physics A {\bf 20}, L1193-L1201}

\bibitem {Gr}
{G.R. Grimmett,
{\em Percolation}, Springer, New-York, 1989.
}

\bibitem {GLN}
{G.R. Grimmett (1997),
Percolation and disordered systems,
Ecole d'\'et\'e de Probabilit\'es de St-Flour XXVI, L.N. Math. {\bf 1665},
153-300}
 
 \bibitem {H}
 {O. H\"aggstr\"om (1995), Random-cluster Measures and Uniform Spanning Trees, Stoch. Proc. Appl. {\bf 59}, 267-275}
 
\bibitem {Hay}
{W.K. Hayman,
{\em Multivalent functions}, CUP, 1994 (second edition).}

\bibitem{IW}
{
N. Ikeda and S. Watanabe,
{\em Stochastic Differential Equations and Diffusion Processes},
Second edition, North-Holland, 1989.
}

\bibitem{Kennalgo}
{T. Kennedy (2002),
A faster implementation of the pivot algorithm 
for self-avoiding walks, J. Stat. Phys. {\bf 106},
 407-429.}

\bibitem {Kenn}
{T. Kennedy (2002),
Monte Carlo Tests of Stochastic Loewner Evolution Predictions 
for the 2D Self-Avoiding Walk, Phys. Rev. Lett. {\bf 88},
 130601. 
}

\bibitem {K0}
{R. Kenyon (1997),
Local statistics of lattice dimers, Ann. Inst. 
Henri Poincar\'e {\bf 33}, 591-618.}

\bibitem {Kdim}
{R. Kenyon (1999), 
Dim\`eres et arbres couvrants,
in {\sl Math\'ematique et Physique}, 
SMF Journ. Annu., 1-14.
}

 \bibitem {K1}
{
R. Kenyon (2000),
Conformal invariance of domino tiling, Ann. Probab.
{\bf 28}, 759-785.}

\bibitem{K2}
{
R. Kenyon (2000), The asymptotic determinant of the discrete
Laplacian, Acta Math.
{\bf 185}, 239-286.}

\bibitem {K3}
{R. Kenyon (2000),
Long-range properties of spanning trees in $\Z^2$,
J. Math. Phys. {\bf 41} 1338--1363.}

\bibitem {Kesbook}
{H. Kesten, {\em Percolation theory for mathematicians,}
Birh\"auser, Boston, 1982.}

\bibitem {Kes}
{H. Kesten (1987),
Scaling relations for 2D-percolation,
Comm. Math. Phys. {\bf 109}, 109-156.}

\bibitem {KSZ}
{H. Kesten, V. Sidoravicius, Yu. Zhang (2001),
Percolation of Arbitrary words on the Close-Packed Graph of $\Z^2$, 
Electr. J. Prob. {\bf 6}, paper no. 4.
} 


\bibitem {KPZ}
{V.G. Knizhnik, A.M. Polyakov, A.B. Zamolodchikov (1988),
Fractal structure of 2-D quantum gravity,
Mod. Phys. Lett. {\bf A3}, 819.
}

\bibitem {Ko}
{G. Kozma (2002),
Scaling limit of loop erased random walk - a naive approach,
preprint.}


\bibitem {Ku}
{P.P. Kufarev (1947),
A remark on integrals of the Loewner equation, 
Dokl. Akad. Nauk SSSR {\bf 57}, 655-656.}

 \bibitem {LPS}
{R. Langlands, Y. Pouillot, Y. Saint-Aubin (1994),
Conformal invariance in two-dimensional percolation,
Bull. A.M.S. {\bf 30}, 1--61}.

\bibitem {Lduke}
{G.F. Lawler (1980),
 A self-avoiding random walk, Duke Math. J. {\bf 47},
 655-694.
}
\bibitem {L1} {G.F. Lawler,
{\em Intersections of Random Walks,}
Birkh\"auser, Boston, 1991.}

\bibitem {L2}
{G.F. Lawler (1995),
 Nonintersecting planar Brownian motions, 
Mathematical Physics Electronic Journal {\bf 1},
paper no.1.}  

\bibitem {L3}{
G.F. Lawler (1996),
Hausdorff dimension of cut points for Brownian motion,
Electron. J. Probab. {\bf 1}, paper no.2.}

\bibitem {L4}{
G.F. Lawler (1996), The dimension of the frontier of planar Brownian
motion,
Electron. Comm. Prob. {\bf 1}, paper no.5.}

\bibitem {L5}{
G.F. Lawler (1997),
The frontier of a Brownian path is multifractal, preprint. }

\bibitem {Lsci}{
G.F. Lawler (1998),
Strict concavity of the 
intersection exponent for Brownian motion in two and three dimensions,
 Mathematical Physics Electronic Journal {\bf 5},
paper no. 5. }

\bibitem {LLERW}{G.F. Lawler (1999),
Loop-erased random walk, in {\sl Perplexing problems
in Probability}, Prog. Prob. {\bf 44}, Birkh\"auser,
 197-217.}

\bibitem {Lbuda}
{G.F. Lawler (1999),
Geometric and fractal properties of Brownian motion and random walk paths in two and three dimensions,
 Bolyai Mathematical Society Studies, {\bf 9}, 219-258.}

\bibitem {LLN}
{G.F. Lawler (2001), 
An introduction to the stochastic Loewner evolution,
preprint.}


\bibitem {LP1}
{ G.F. Lawler, E.E. Puckette (1997),
The disconnection exponent for simple random walk,
Israel J. Math. {\bf 99}, 109-122.}

\bibitem {LP2}
{ G.F. Lawler, E.E. Puckette (2000),
The intersection exponent for simple random walk,
Combin. Probab. Comput. {\bf 9}, 441-464.}

\bibitem {LSW1}
{G.F. Lawler, O. Schramm, W. Werner (2001),
Values of Brownian intersection exponents I: Half-plane exponents,
Acta Mathematica {\bf 187}, 237-273. }

\bibitem {LSW2}
{G.F. Lawler, O. Schramm, W. Werner (2001),
Values of Brownian intersection exponents II: Plane exponents,
Acta Mathematica {\bf 187}, 275-308.}

\bibitem{LSW3}
{G.F. Lawler, O. Schramm, W. Werner (2002),
Values of Brownian intersection exponents III: Two sided exponents,
Ann. Inst. Henri Poincar\'e {\bf 38}, 109-123.}

\bibitem {LSWan}
{G.F. Lawler, O. Schramm, W. Werner (2002),
Analyticity of planar Brownian intersection exponents,
Acta Mathematica {\bf 189}, to appear.}
 
\bibitem {LSW4/3}
{G.F. Lawler, O. Schramm, W. Werner (2001),
The dimension of the planar Brownian frontier is $4/3$, 
Math. Res. Lett. {\bf 8}, 401-411. 
}
 
\bibitem {LSWup}
{G.F. Lawler, O. Schramm, W. Werner (2001),
Sharp estimates for Brownian non-intersection probabilities,
in: {\sl In and Out of Equilbrium}, V. Sidoravicius Ed.,
 Prog. Probab. {\bf 51}, Birkh\"auser, 113-131.}

\bibitem {LSW5}
{G.F. Lawler, O. Schramm, W. Werner (2002),
One-arm exponent for critical 2D percolation,
Electronic J. Probab. {\bf 7}, paper no.2.}

\bibitem {LSWlesl}
{G.F. Lawler, O. Schramm, W. Werner (2001),
Conformal invariance of planar loop-erased random
walks and uniform spanning trees, preprint.}

\bibitem {LSWSAW}
{G.F. Lawler, O. Schramm, W. Werner (2002),
On the scaling limit of planar self-avoiding walks, 
preprint.}

\bibitem {LSWrest}
{G.F. Lawler, O. Schramm, W. Werner (2002),
Conformal restriction properties. The chordal case,
preprint.}

\bibitem {LW1}
{G.F. Lawler, W. Werner (1999),
Intersection exponents for planar Brownian motion,
Ann. Probab. {\bf 27}, 1601-1642.}

\bibitem {LW2}
{G.F. Lawler, W. Werner (2000),
Universality for conformally invariant intersection
exponents, J. Europ. Math. Soc. {\bf 2}, 291-328.}

\bibitem {LWbls}
{G.F. Lawler, W. Werner (2003),
The Brownian loop-soup, 
preprint.}

\bibitem {Le}
{N.N. Lebedev, {\em Special Functions and their Applications},
transl. from russian, Dover, 1972.}

\bibitem {LG}
{J.F. Le Gall (1992), Some properties of planar Brownian motion,
Ecole d'\'et\'e de Probabilit\'es de St-Flour XX, L.N. Math. {\bf 1527},
111-235.}

\bibitem {LV}
{O. Lehto, K.I. Virtanen, 
{\em Quasiconformal mappings in the plane}, second edition,
translated from German, Springer, New York, 1973.
}

\bibitem {Lev}
{P. L\'evy, 
{\em Processus Stochastiques et Mouvement Brownien},  
Gauthier-Villars, Paris, 1948.}
 
\bibitem {Lo}
{K. L\"owner (1923),
Untersuchungen \"uber schlichte konforme Abbildungen des
Einheitskreises I., Math. Ann. {\bf 89}, 103--121.}

\bibitem {Lyons}
{R. Lyons (1998),
A bird's-eye view of uniform spanning trees and forests,
 in {\sl Microsurveys in Discrete Probability},
 D. Aldous and J. Propp eds.,
 Amer. Math. Soc., Providence,  135--162.}

\bibitem {MS}{
N. Madras, G. Slade,
{\em The Self-Avoiding Walk}, Birkh\"auser, 1993.}

\bibitem {Maj}{S.N. Majumdar (1992),
Exact fractal dimension of the loop-erased random walk
in two dimensions,
Phys. Rev. Lett. {\bf 68}, 2329--2331.}

\bibitem {Ma}
{B.B. Mandelbrot,
{\em The Fractal Geometry of Nature},
Freeman, 1982.}

\bibitem{MR}
{D.E. Marshall, S. Rohde (2001), The Loewner differential 
equation and slit mappings, preprint.}


\bibitem {N}
{B. Nienhuis, E.K. Riedel, M. Schick (1980),
Magnetic exponents of the two-dimensional $q$-states Potts 
model, J. Phys A {\bf 13}, L. 189-192.}

\bibitem {N2}
{B. Nienhuis (1984),
Coulomb gas description of 2-D critical behaviour,
J. Stat. Phys. {\bf 34}, 731-761.}

\bibitem {N3}
{B. Nienhuis (1987),
Coulomb gas formulation of two-dimensional phase transitions, 
in {\sl Phase transitions and critical phenomena} {\bf 11},
Academic Press, 1--53. 
}

\bibitem{Pe}
{R.P. Pearson (1980),
Conjecture for the extended Potts model magnetic eigenvalue,
Phys. Rev. B {\bf 22}, 2579-2580.}


\bibitem {Pem}
{R. Pemantle (1991),
Choosing a spanning tree for the integer lattice uniformly,
 Ann. Probab. {\bf 19},  1559-1574. }

\bibitem {Po}
{A.M. Polyakov (1974),
A non-Hamiltonian approach to conformal field theory,
Sov. Phys. JETP {\bf 39}, 10-18.
}

\bibitem {P1}
{C. Pommerenke (1966),
On the L\"owner differential equation,
Michigan Math. J. {\bf 13}, 435--443.}

\bibitem {P2}
{C. Pommerenke,
{\em Boundary Behaviour of Conformal Maps},
Springer-Verlag, 1992.}

 \bibitem {RY}
{D. Revuz, M. Yor,
{\em Continuous Martingales and Brownian Motion}, Springer-Verlag, 1991.}

 \bibitem {RS}
{S. Rohde, O. Schramm (2001), 
Basic properties of SLE, preprint.}

\bibitem {Ru}
{W. Rudin,
{\em Real and Complex Analysis}, Third Ed., McGraw-Hill, 
1987.}

\bibitem {Rus}
{L. Russo (1978),
A note on percolation,
Z. Wahrscheinlichkeitsth. verw. Geb. {\bf 56}, 229-237.}


\bibitem {SD}
{H. Saleur, B. Duplantier (1987),
Exact determination of the percolation
hull exponent in two dimensions,
Phys. Rev. Lett. {\bf 58},
2325.}


\bibitem{SRG}
{B. Sapoval, M. Rosso, J. F. Gouyet (1985),
The fractal nature of a diffusion front 
and the relation to percolation,
J. Physique Lett. {\bf 46}, L149-L156}



\bibitem {S1}{
O. Schramm (2000), Scaling limits of loop-erased random walks and
uniform spanning trees, Israel J. Math. {\bf 118}, 221--288.}

\bibitem {S2}
{O. Schramm (2001), A percolation formula,
 Electr. Comm. Probab. {\bf 6}, 115-120.
}

\bibitem {SW}
{P.D. Seymour, D.J.A. Welsh (1978),
Percolation probabilities on the square lattice,
in {\sl Advances in Graph Theory},
 ann. Discr. Math. {\bf 3}, North-Holland, 227-245.
}

\bibitem {Sm1}
{S. Smirnov (2001),
Critical percolation in the plane: conformal invariance,
 Cardy's formula, scaling limits,
 C. R. Acad. Sci. Paris Sér. I Math. {\bf 333},  239--24}

\bibitem {Sm2}
{S. Smirnov, in preparation.}

\bibitem {SmW}
{S. Smirnov, W. Werner (2001),
Critical exponents for two-dimensional percolation, 
Math. Res. Lett. {\bf 8}, 729-744.}

\bibitem {V}
{B. Virag (2003), Brownian beads, preprint.}

\bibitem {VW}
{S.R.S. Varadhan, R.J. Williams 
(1985), 
Brownian motion in a wedge with oblique reflection. 
Comm. Pure Appl. Math. {\bf 38}, 405--443. 
}

\bibitem {W1}
{W. Werner (1994),
Sur la forme des composantes connexes du compl\'ementaire de la courbe 
brownienne plane,
 Probab. Theory Related Fields {\bf 98}, 307--337.}

\bibitem {Wbd}
{W. Werner (1996),
 Bounds for disconnection exponents, Electr. Comm. Probab.
{\bf 1}, 19-28.}

\bibitem {W3} 
{W. Werner (1997),
Asymptotic behaviour of disconnection and non-inter\-section exponents,
 Probab. Theory Related Fields {\bf 108},
131-152.}

\bibitem {W4} 
{W. Werner (2001), 
Critical exponents, conformal invariance and planar Brownian motion,
in {\sl Proceedings of the 4th ECM Barcelona 2000},
Prog. Math. {\bf 202}, Birkh\"auser, 87-103.}

\bibitem {Whid}
{W. Werner (2003),
Girsanov's theorem for $SLE(\kappa, \rho)$ processes, intersection exponents and
hiding exponents, preprint.}

\bibitem {WiD}
{D.B. Wilson (1996), 
Generating random spanning trees more quickly than the cover time,
 Proceedings of the Twenty-eighth Annual ACM Symposium on the Theory of
Computing (Philadelphia, PA, 1996), 296--303. }

\bibitem {Z}
{D. Zhan (2003),
preprint.}

\end {thebibliography}
\end {document}